\documentclass[]{imsart}

\usepackage{latexsym,amsmath,amssymb,amscd,amsfonts,amsthm,mathrsfs}
\usepackage[latin1]{inputenc}
\usepackage{natbib}
\usepackage{color}
\usepackage{geometry}
\usepackage{dsfont}

\theoremstyle{plain}

\newtheorem{thm}{Theorem}
\newtheorem{prop}{Proposition}
\newtheorem{cor}{Corollary}
\newtheorem{lem}{Lemma}
\newtheorem{Rem}{Remarks}

\def \P {\mathbb P}
\def \Pc {\mathcal P}
\def \S {\mathcal S}
\def \R {\mathbb R}
\def \N {\mathbb N}
\def \Z {\mathbb Z}
\def \E {\mathbb E}

\def \H {\mathcal H}

\def \cinf {\xi}
\def \xM {\alpha}

\DeclareMathOperator*{\KL}{KL}
\DeclareMathOperator*{\argmin}{argmin}

\DeclareMathOperator*{\pen}{pen}
\DeclareMathOperator*{\critp}{crit}
\def \M {\mathcal M}
\def \mus {\bar \mu}
\def \lambdai {\underline \lambda}
\def \lambdas {\bar \lambda}
\def \mel {\wp}
\def \B {\mathfrak{B}}

\begin{document}
\begin{frontmatter}

\title{Adaptive density estimation for clustering with Gaussian mixtures}
\runtitle{Adaptive density estimation using finite Gaussian mixtures}

\begin{aug}
  \author{Cathy Maugis-Rabusseau\ead[label=e1]{cathy.maugis@insa-toulouse.fr}}
  \address{Institut de Math\'ematiques de Toulouse, INSA de Toulouse, Universit\'e de Toulouse\\
           INSA de Toulouse,\\ 135, avenue de Rangueil,\\  31077 Toulouse Cedex 4, France.\\
           \printead{e1}\\
           \vspace*{0.2cm}
  }
  \and
  \author{Bertrand Michel\ead[label=e2]{bertrand.michel@upmc.fr}}
  \address{Laboratoire de Statistique Th\'eorique et Appliqu\'ee,\\Universit\'e Pierre et Marie Curie - Paris 6,\\  4 place Jussieu, 75252 Paris cedex 05 France\\
          \printead{e2}
  }
\runauthor{Maugis and Michel}

\end{aug}

\begin{abstract}
Gaussian mixture models are widely used to study clustering problems. These
model-based clustering methods require an accurate estimation of the unknown
data density by Gaussian mixtures. In Maugis and Michel
(2009)\nocite{MaugisMichel09}, a penalized maximum likelihood estimator is
proposed for automatically selecting the number of mixture components. In the
present paper, a collection of univariate densities whose logarithm is locally
$\beta$-Hölder with moment and tail conditions are considered. We show that
this penalized estimator is minimax adaptive to the $\beta$ regularity of such
densities in the Hellinger sense.
\end{abstract}




\end{frontmatter}

\section{Introduction}
Clustering methods consists of discovering clusters among observations. Many
cluster analysis methods have been proposed in statistics and learning theory,
roughly fall into three categories. The first one is based on similarity or
dissimilarity distances, the best-known are partitioned clustering methods as
k-means and the hierarchical clustering methods \citep[see for instance
Sections 14.3.6 and 14.3.12 in][]{Hastie09}. The second category consists of
density level set clustering methods which consider clusters as the connected
components of high density regions \citep[see][]{Hartigan75}. The third
category is composed of model-based clustering methods which define clusters as
observations having most likely the same distribution. In this last case, each subpopulation is
assumed to be distributed from a parametric density, like a Gaussian one and
thus the unknown data density is a mixture of these distributions \citep[see
for instance][]{McLachlan:Peel:2000}. The data clustering is then deduced
thanks to the maximum a posteriori (MAP) rule. The clustering problem being
based on data density estimation, it is then essential that this density be
efficiently estimated.

Because of their wide range flexibility, Gaussian mixture densities are widely
used to model the unknown distribution of continuous data for clustering
analysis \citep[see for instance][]{Lindsay95,McLachlan:Peel:2000}. By
recasting the clustering problem into a model selection problem, we have
proposed in \citet{MaugisMichel09} a non asymptotic penalized criterion. We
proved that the selected Gaussian mixture estimator fulfills an oracle
inequality. The aim of this new paper is to investigate the adaptive properties
of this estimator in order to justify the validity of our clustering method.
More precisely, adapting a recent approximation result, we show that our
estimator is minimax adaptive to the regularity parameter of a particular class
of Hölder spaces defined further. As far as we know, such a minimax adaptive
result has never been shown for a density estimator used for model-based
clustering methods.

We first recall the context of \citet{MaugisMichel09} in the unidimensional
case. Let us consider $n$ independent identically distributed random variables
$X_1,\ldots,X_n$ with values in $\R$. Their common unknown density $s$ belongs
to the set $\S$ of all density functions with respect to the Lebesgue measure
on $\R$. The considered unidimensional Gaussian mixtures are characterized by
their number of components $m$ and their means and variances parameters are
assumed to be bounded. These mixture densities are grouped into a model
collection $(\S_m)_{m\in\M_n}$, subsets of $\S$, defined by
\begin{equation}
   \S_{m}=\left\{x\in\R\mapsto \underset{u=1}{\stackrel{m}{\sum}}p_u \psi_{\sigma_u}(x-\mu_u); \mu_u\in[-\mus,\mus], \sigma_u^2\in[\lambdai,\lambdas]; p_u\in[0,1],  \underset{u=1}{\stackrel{m}{\sum}}p_u =1\right\}
   \label{defmodelSm}
\end{equation}
where $\psi$ is the Gaussian kernel defined by $\psi(x)=\pi^{-\frac 1 2}
\exp(-x^2)$ for all $x\in\R$ and  $\psi_\sigma(\cdot)=\sigma^{-1}
\psi\left(\frac{\cdot}{\sigma}\right)$ for all $\sigma>0$. The number of free
parameters, common to all the mixture densities of a given model $\S_m$ is
called dimension and is denoted $D(m)$. Considering a non asymptotic point of
view \citep[see for instance][]{Mas07}, the three bounds $\mus$, $\lambdai$ and
$\lambdas$ of each model $\S_m$ and also the maximum number of mixture
components in the collection may depend on $n$. Such mixtures are called sieves
according to the terminology introduced by \citet{Grenander}.

Over each model $\S_m$, a maximum likelihood estimator (MLE) $\hat s_m$ is
obtained by minimizing the empirical contrast
$$
    \gamma_n(t)=-\frac 1 n \sum_{i=1}^{n}\ln\left\{t(X_i)\right\}.
$$
The loss function associated to the likelihood contrast is the Kullback-Leibler
divergence: For two densities $s$ and $t$ in $\S$, the Kullback-Leibler
divergence is defined by
$$
    \KL(s,t)= \int \ln \left\{\frac{s(x)}{t(x)}\right\} s(x)\,dx
$$
if $s dx$ is absolutely continuous with respect to $t dx$ and $+\infty$
otherwise. The model $m^\star$ in the collection minimizing the
Kullback-Leibler risk
$$
    m^\star \in \underset{m\in\mathcal{M}_n}{\argmin}\ \E_s[\KL(s,\hat s_m)]
$$
is considered as the ``best'' model of the collection. Nevertheless this best
model $m^\star$ and also the associated density $\hat s_{m^\star}$ (called
oracle) are unknown since they depend on the true density $s$. A model $\hat m$
is then chosen by minimizing over $\M_n$ the following penalized criterion
$$
    \critp(m)=\gamma_n(\hat{s}_m)+\pen(m) .
$$
The penalty function $\pen:m \in\M_n\mapsto \pen(m)\in\mathbb{R}^+$ has to be
chosen such that the  Kullback-Leibler risk $\E_s[\KL(s,\hat s_{\hat m})]$ of
$\hat s_{\hat m}$ is close to the oracle risk $\E_s[\KL(s,\hat s_{m^\star})]$.
The construction of such penalties is proposed in Theorem 2.2 in
\citet{MaugisMichel09}. This theorem can be stated as follows in the univariate
context, where $d_H(g,h) = \frac{1}{\sqrt{2}}{\|\sqrt{g}-\sqrt{h}\|}_2$ denotes
the Hellinger distance between two densities $g$ and $h$ of $\S$:

\begin{thm}\label{ThMauMic09}
There exists two absolute constants $\kappa$ and $C$ such that, if
$$
    \pen(m)\geq\kappa\frac{D(m)}{n}\left\{1+2\,\mathcal{A}^2 + \ln\left(\frac{1}{1\wedge\frac{D(m)}{n}\,\mathcal{A}^2}\right)\right\}
$$
where
\begin{equation}
    \mathcal{A}=\sqrt{\ln(6\pi e^2)} + \sqrt{\pi} + \sqrt{\ln\left(\mus\sqrt{\frac{8}{c_1 \lambdai}}\right)} + \sqrt{\ln\left(\frac{144 \lambdas}{\lambdai}\right)},
\label{defA}
\end{equation}
then the model $\hat{m}$ minimizing
$$
    \critp(m)=\gamma_{n}(\hat s_m)+\pen(m)
$$
over $\M_n$ exists and
\begin{equation}
    \E\left[d_H^{\,2}(s,\hat{s}_{\hat{m}})\right]\leq \mathcal{C} \left[\inf_{m\in\M_n}\{\KL(s,\S_{m})+\pen(m)\}+\frac{1}{n}\right].
    \label{OracleIneq}
\end{equation}
\end{thm}

Note that a similar result can be found in \citet{MaugisMichel09} for
multivariate data clustering with variable selection. The method has been
successfully implemented and tested in practice \citep[see][]{MauMic10}.

Minimax adaptive estimation has been intensively studied in nonparametric
statistics, see for instance \citet{Tsybakov09}, and \citet{Mas07} for adaptive
minimax methods based on $l_0$ penalization. A natural optimality criterion is
the minimax risk, first introduced by \citet{Wolfowitz50}.  Let
$$
  \mathcal{R}(\tilde s_n, \H_\beta) = \underset{s\in\H_\beta}{\sup} \E_s[d_H^2(s,\tilde s_n)]
$$
be the maximal Hellinger risk of an estimator $\tilde s_n$ of $s$. The minimax
Hellinger risk on a density class $\H_\beta$ is then defined by
$$
  \mathcal{R}_n(\H_\beta)=\underset{\tilde s_n}{\inf}\ \mathcal{R}(\tilde s_n,\H_\beta)
$$
where the infimum is taken over all the possible estimators $\tilde s_n$ of
$s$. An estimator is said to be minimax on $\H_\beta$ if its maximal risk over
$\H_\beta$ reaches the minimax risk on this density class. Let us now consider
a collection $(\H_\beta)_{\beta \in \mathcal{B}}$ of density classes indexed by
a set $\mathcal{B}$ of regularity parameters $\beta$. An estimator is said to
be minimax adaptive if it reaches the minimax risk over $\H_\beta$ for all
$\beta$ of $\mathcal{B}$, without using the knowledge of $\beta$. In order to
motivate the clustering method based on Gaussian mixture estimator $\hat
s_{\hat m}$ proposed in \citet{MaugisMichel09}, we prove in this new paper that this
estimator is minimax adaptive over a particular collection of Hölder density
classes $(\H_\beta)_{\beta \in \mathcal{B}}$ defined further. Of course,
adaptive density estimation in one dimension is now a classical problem and
several adaptive estimators have been already proposed such as kernel
estimators or thresholding wavelet estimators. Nevertheless, although these
alternative methods maybe perform better than our penalized estimator
$\hat{s}_{\hat{m}}$ concerning density estimation in general, these have no
interest for clustering purposes.

The link between model selection and adaptive estimation is made through
approximation theory. Indeed, an adaptive estimation is possible only for
functional classes $\H_\beta$ that can be efficiently approximated by our
Gaussian mixture collection. Convolution is widely used in approximation theory
and many results are known on this topic. It is well known that the convolution
of a density $f$ with scaled versions $\psi_\sigma$ of the Gaussian kernel
$\psi$ converges to $f$ \citep[see for instance][chapter 20]{ChL09}. The
so-called quasi-interpolation method consists of replacing the functions $\psi
_\sigma * f$ by infinite linear combinations of scaled and translated Gaussian
kernels \citep[see for instance][chapter 36]{ChL09}. In a recent paper of
\citet{Hang10}, a nonlinear approximation algorithm based on
finite combinations of scaled and translated Gaussian kernels is defined to
give some approximation results in $L^{p} $ norm on some particular density
classes. Nevertheless, all these results cannot be straightly applied to study
the approximation capacities of Gaussian mixtures. Indeed, the coefficients in
these linear combinations are not necessary positive and their sum is not
constrained to be equal to one. Furthermore, the approximation results provided
by all these methods are not given for the Kullback-Leibler divergence as
required by our statistical context.

The approximation capacity of Gaussian mixtures has also been studied in non
parametric Bayesian works. Lemma 3.1 in \citet{GhosalVaart01} gives a
discretization result for Gaussian mixtures: assume that $s$ is a location or
location-scale mixture with a mixing distribution compactly supported or with
sub-Gaussian tails, $s$ can be approximated by a finite Gaussian mixture with a
small number of components, the error being controlled in $L_1$ and $L_\infty$
norms. In \citet{GhosalVaart07b}, these authors take advantage of this method
for approximating by finite Gaussian mixtures some twice continuously
differentiable functions with additional regularity conditions. More recently,
\citet{KRV} prove an approximation result by finite Gaussian
mixtures for densities whose logarithm is locally Hölder. Their approximation
result is given for the  Kullback-Leibler divergence. This last result can be
successfully  adapted in our context to control the bias term in the right side
term of the oracle inequality (\ref{OracleIneq}) on these particular density
classes. Concerning approximation, the contribution of our work consists of
checking that the non explicit constants of the approximation bounds given in
\citet{KRV} are actually uniform over a density class $\H_{\beta}$ we define.
For easier reading, all the approximation results are given and proved in this
preprint version although a large part of them can be found in \citet{KRV}.

The paper is organized as follows: The main results are presented in
Section~\ref{sect:MainResults}. The density classes $\H_{\beta}$ are introduced
in Section~\ref{subs:Spaces} and an approximation result, adapted of
\cite{KRV}, is given in Section~\ref{subs:ThApprox}. Next, a lower bound of the
minimax risk is given in Section~\ref{subs:LowerBound} and the adaptive
property of our penalized Gaussian mixture estimator on these density classes
$\H_{\beta}$ is addressed in Section~\ref{subs:adaptative}. The approximation
result, the lower bound and the adaptive result are respectively proved in
Sections \ref{sect:ProofApprox}, \ref{ProofLowB} and \ref{ProofAdapt}. Finally,
some technical results are developed in Appendices~\ref{sect:AppendixApprox}
and \ref{sect:AppendixMinor}.

\section{Main results}\label{sect:MainResults}
\subsection{The density classes $\mathcal{H}(\beta,\Pc)$}\label{subs:Spaces}
The adaptation result given further requires a slightly modified version of the approximation result by finite Gaussian mixtures proved in \cite{KRV}. This approximation result concerns densities whose logarithm is locally $\beta$-Hölder and that fulfills additional tail, moments and monotonicity conditions. More precisely, let $\beta >0$, $r = \lfloor \beta \rfloor$ be the largest integer less than $\beta$ and $k \in \N$ such that $\beta \in  (2k, 2k + 2]$. Let also $\Pc$ be the set of parameters $\left\{\gamma,l^+,L,\varepsilon,C,\xM,\cinf,M\right\}$ where $L$ is a polynomial function on $\R$ and the other parameters are positive constants. We then define the density class $\mathcal{H}(\beta,\Pc)$ of all densities $f$ satisfying the following conditions:
\begin{enumerate}
\item {\bf Smoothness.} $\ln f$ is assumed to be locally $\beta$-H\"older: for all $x$ and $y$ such that $|y-x| \leq \gamma$,
    \begin{equation}
    \left|(\ln f)^{(r)}(x) - (\ln f)^{(r)}(y) \right|  \leq r! \, L(x)  |y-x |^{\beta-r}.
    \label{Cond1}
    \end{equation}
Furthermore for all $j \in \{0,\ldots,r \}$,
    \begin{equation}
    | (\ln f)^{(j)}(0) | \leq l^+ .
    \label{Cond1b}
    \end{equation}

\item {\bf Moments.} The derivative functions $(\ln f)^{(j)}$ for $j=1,\ldots,r$ and the polynomial function $L$ fulfill
    \begin{equation}
    \int_{\mathbb{R}} \left| (\ln f)^{(j)}(x) \right| ^{\frac{2 \beta + \varepsilon}{j}} f(x) dx \leq C \ ,   \quad  \quad \int_{\mathbb{R}} \left| L(x) \right|^{2  + \frac{\varepsilon}{\beta}} f(x) dx \leq C.
    \label{Cond2a}
    \end{equation}

\item {\bf Tail.} For all $x\in \R$,
    \begin{equation}
    f(x) \leq M \psi(x) .
    \label{Cond2b}
    \end{equation}

\item {\bf Monotonicity.} $f$ is strictly positive, $f$ is nondecreasing on $(-\infty,-\xM)$ and nonincreasing on $(\xM, \infty)$, and $f(x) \geq \cinf$ for all $x\in [-\xM,\xM]$.
\end{enumerate}

\begin{Rem} The monotonicity assumption can be relaxed by assuming that there exist two constants $c>0 $ and $\bar \sigma>0$ such that $\forall 0<\sigma<\bar{\sigma},\ \forall x\in\mathbb{R}$, $$\frac{K_{\sigma} f(x)}{f(x)} \geq c. $$ This condition corresponds to the first point given in Lemma~\ref{BoiteOutilBornes0} in Appendix~\ref{sect:AppendixApprox} which is a key point to prove the approximation result. In the following, the strong monotonicity condition is assumed in the definition of the density class $\mathcal{H}(\beta,\Pc)$ to simplify the proofs of the lower bound. \end{Rem}

\begin{Rem}
For easier reading, the monotonicity assumption is stated on a symmetric interval but it is possible to consider this assumption on a general interval $[\alpha_1,\alpha_2]$ with $\alpha_1<\alpha_2$. This monotonicity assumption allows us to lower bound the convolution $f\ast \psi_\sigma$ by $f$ up to a multiplicative constant  according to Remark 3 in \citet{GGR99}.
\end{Rem}

\begin{Rem}
These density classes are more restrictive than those considered in \citet{KRV}: Indeed the upper bounds in (\ref{Cond2a}) have to be uniform on the density class $\H(\beta,\Pc)$ and we also need the additional Condition (\ref{Cond1b}). These restrictions allow us to control the Kullback-Leibler divergence between a density
of $\H(\beta,\Pc)$ and a convenient finite Gaussian mixture, uniformly over $\H(\beta,\Pc)$. Note that Condition (\ref{Cond2b}) is here assumed on $\R$ but it could be assumed only outside an interval as in \citet{KRV}.
\end{Rem}

\begin{Rem}\label{RemLarger}
In the sequel, $\Pc'$ is said to be "larger than" $\Pc$ if at least one of the following conditions is fulfilled:
\begin{itemize}
\item at least one constant among $M$, $C$ or $l^+$ of $\Pc'$ is larger than the corresponding one of $\Pc$,
\item the constant $\gamma$ of $\Pc'$ is smaller than the corresponding one of $\Pc$,
\item   for all $x\in\mathbb{R},\ L(x)\leq L'(x)$ where $L$ (resp. $L'$) belongs to $\Pc$ (resp. $\Pc'$)\\
\end{itemize}
\end{Rem}

\subsection{Approximation result}\label{subs:ThApprox}
For any function $f$, $K_\sigma f$ denotes the convolution $f\ast\psi_\sigma$ and $\Delta_\sigma f$ is the error term $K_\sigma f-f$. As explained in \cite{KRV}, for a $\beta$-smooth density $f$ with $\beta\leq 2$ and under reasonable regularity assumptions, it is possible to define a finite location-scale Gaussian mixture $\mel_\sigma$ such that $\KL(f,\mel_\sigma)=O(\sigma^{2\beta})$. The usual approach consists of discretizing the continuous mixture $K_\sigma f$. But as $\|f - K_\sigma f\|_\infty$ remains of order $\sigma^2$ when $\beta>2$, this approach appears to be inefficient for smoother densities. An alternative strategy is proposed in Kruijer et al. \cite{KRV}, based on the following successive convolutions of $f$: $f_0=f$ and for all $j\geq 0$, $f_{j+1}=f - \Delta_\sigma f_j$. In their paper, the density is approximated by a discretized version of the continuous mixture $K_\sigma f_k$ where  $k \in \N$ is such that $\beta \in  (2k, 2k + 2]$.

In our framework, Lemma 4 in \cite{KRV} cannot be directly used since the upper bound over the Kullback-Leibler divergence between $f$ and the finite Gaussian mixture is not uniform over $\H(\beta,\Pc)$. Thus some additional work is necessary in order to prove an uniform version of this approximation result. Another reason for revisiting the approximation results given in \citet{KRV} is that these ones are stated for $\sigma \leq \bar{\sigma}$ where $\bar{\sigma}$ depends on the approximated density $f$. Thus we also need to check that  it is possible to choose the same $\bar{\sigma}$ for all the densities of $\mathcal{H}\left( \beta,\Pc \right)$. The proof of Theorem~\ref{ThApproximation} consists of carefully following the method of \cite{KRV} in order to obtain this uniform version. A sketch of the proof is given below and a self-contained proof is detailed in Section~\ref{sect:ProofApprox}.

\begin{thm}~\label{ThApproximation}
There exists a positive constant $\bar \sigma(\beta)<1$ such that for all $f\in\H(\beta,\Pc)$ and for all $\sigma < \bar \sigma(\beta)$, there exists a finite Gaussian mixture of density $\mel_\sigma$ with less than $G_\beta \sigma^{-1} |\ln \sigma|^{\frac 3 2}$ support points, with the same variance $\sigma$ for each component and with means belonging to $[-\mu_\sigma,\mu_\sigma]$
where
$$
    \mu_\sigma \leq \tilde G_\beta |\ln \sigma|^{\frac 1 2}
$$
such that
\begin{equation}
\label{ctrapproxKL}
    \KL(f,\mel_\sigma)=\int_{\mathbb{R}} f(x) \ln\left(\frac{f(x)}{\mel_\sigma(x)}\right)dx  \leq c_\beta\ \sigma^{2\beta}
\end{equation}
where  $c_{\beta}$ is uniform on $\mathcal{H}\left( \beta,\Pc \right)$ and continuous on $\beta$ . The constant $\bar\sigma(\beta)$ only depends on $\H(\beta,\Pc)$ and is a continuous function of $\beta$. Moreover, $G_\beta$ and $\tilde G_\beta$ are two positive constants that only depend on $\H(\beta,\Pc)$, and are both increasing functions of $\beta$.
\end{thm}
The two constants $G_\beta$ and $\tilde G_\beta$ are explicitly defined by Equations (\ref{Gtildebeta}) and (\ref{Gbeta}) in the proof of Theorem~\ref{ThApproximation} in Section~\ref{subsect:ProofThApproximation}.

\begin{proof}[Sketch of the proof]
Let $f$ be a density in a given class $\H(\beta,\Pc)$. First, the convolution $K_\sigma f_k$ is shown  to be close to $f$ on a subspace of $\R$ where the derivative functions of $\ln f$ and $L$ are efficiently controlled (see Lemma~\ref{Lemma1}). On this subspace, the difference $K_\sigma f_k -f$ is controlled by $f(x) R_f(x) O(\sigma^\beta)$, apart from a term $\sigma^H$ where $H$ can be arbitrarily large. The term $O(\sigma^\beta)$ is uniform on $\H(\beta,\Pc)$ and $R_f$ is a polynomial function of $L$ and the derivative functions of $\ln f$. Next, since $f_k$ is not necessarily a positive function, a density function $h_k$ is defined from $f_k$. The previous result is then adapted for controlling $K_\sigma h_k -f$ on a more restrictive subspace of $\R$ (see Lemma~\ref{Lemme2}). Based on  this result, a control of the Kullback-Leibler divergence between $f$ and the continuous Gaussian mixture $K_\sigma h_k$ is obtained in Proposition~\ref{PropKLFKsigmahk}: $\KL(f,K_\sigma h_k)\leq c_\beta \sigma^{2\beta}$ where $c_\beta$ is a multiplicative constant uniform on $\H(\beta,\Pc)$. Finally, using a discretization result, a similar control is obtained for $\KL(f,\mel_\sigma)$ where $\mel_\sigma$ is a finite Gaussian mixture fulfilling conditions given in Theorem~\ref{ThApproximation}.
\end{proof}

\subsection{Lower bound}\label{subs:LowerBound}

In order to show that the MLE penalized estimator $\hat s _{\hat m}$ is adaptive to the smoothness parameter $\beta$, a lower bound of the  minimax risk $\mathcal{R}_n(\H(\beta,\Pc))$ is required. For all $0 < \underline \beta < \bar \beta$, a ``large enough'' parameter set $\Pc(\underline \beta, \bar \beta)$ is found  such that for all $\beta \in [\underline \beta, \bar \beta]$ , $\mathcal{H}(\beta,\Pc(\underline \beta, \bar \beta))$ is well defined and a lower bound is given for  the density classes $\mathcal{H}\left(\beta, \Pc(\underline \beta, \bar \beta)\right)$. Note that in Theorem~\ref{ThApproximation}, the constants $c_{\beta}$, $\bar{\sigma}(\beta)$, $G_\beta$ and $\tilde G_\beta$ cannot be bounded uniformly for all $ \beta \in \R^+$. Nevertheless, it can be proved that $\hat s _{\hat m}$ is minimax adaptive on a range of regularity $[\underline \beta,\bar \beta]$.

First, the parameter set $\Pc(\underline \beta, \bar \beta)$ has to be defined rigorously. Its definition is rather technical since it depends on the way the lower bound is proved. The proof is based on the construction of some oscillating functions, this standard method is presented for instance in \citet[][see Section 7.5]{Mas07}. Let us take some infinitely differentiable function $\varphi: \R \rightarrow \R$ with compact support included into $(\frac 1 4,\frac 3 4)$ such that
$$
    \int_{\R} \varphi(x) dx =0 \textrm{ and } \int_{\R} \varphi(x)^2 dx=1.
$$
We set $A = \underset{0\leq k \leq r+1}{\max}\|\varphi^{(k)}\|_{\infty} > 1$ and
let $D$ be some positive even integer. For any positive integer $j\in\{1,\ldots,D\}$, we consider the function
$$
    \begin{array}{r c l}
    \varphi_{j}: \R &\rightarrow& \R\\
                  x & \mapsto & \frac{\cinf D^{-\beta}}{A}\varphi \left( \frac{D}{\xM}(x+\frac{\xM}{2}) - (j-1)\right).
    \end{array}
$$
Moreover, let $\mathcal{T}(\xM,\cinf)$ be the space of functions $\omega:\R\rightarrow \R^+$ such that
$w$ is nondecreasing on $(-\infty,-\frac{\xM}{2})$, nonincreasing on $(\frac{\xM}{2},+\infty)$,
$\omega(x) = 2\cinf$ for all $x\in \left[-\frac{3\xM}{4},\frac{3\xM}{4}\right]$, and $\omega(-\xM)=\omega(\xM)=\cinf$.

Next, let $\tilde\Pc=\left\{\frac{\xM}{4},\ln(2\cinf),\tilde L, \tilde \varepsilon, \tilde C, \xM,\cinf, \tilde M\right\}$ be a parameter set such that $\mathcal{T}(\xM,\cinf) \bigcap \H(\beta,\tilde\Pc)$ is nonempty. Based on a function $\omega\in \mathcal{T}(\xM,\cinf) \bigcap \H(\beta,\tilde\Pc)$ and the functions $\varphi_j$, we consider the functional space
$\mathcal{J}(\beta,D)=\left\{f_{\theta};\ \theta\in\{0,1\}^D\right\}$ where
for all $\theta \in \{0,1\}^D$ and for all $x\in\R$,
\begin{equation}
    f_{\theta}(x) =\omega(x) + \underset{j=1}{\stackrel{D}{\sum}} (2\theta_j -1) \varphi_{j} (x).
    \label{deffvtheta}
\end{equation}

\begin{prop} \label{prop:exist}
There exists a parameter set $\Pc(\underline\beta,\bar\beta)$ such that for all $D \in \N^{*}$ and for all $\beta \in [\underline \beta, \bar \beta]$,
\begin{equation*}
\label{inclusionJH}
\mathcal{J}(\beta,D) \subset  \mathcal{H}\left(\beta, \Pc(\underline \beta, \bar \beta)\right) .
\end{equation*}
\end{prop}

\begin{Rem} Note that if such a parameter set exists, Proposition~\ref{prop:exist} is also true for all the parameter sets larger than it (in the sense given in Remark~\ref{RemLarger}).
A key point to prove the lower bound stated in the next theorem is that the parameter set $ \Pc(\underline \beta, \bar \beta)$ does not depend on $D$.
\end{Rem}

\begin{thm}~\label{ThMinoration}
Suppose that one observes independent random variables $X_1,\ldots,X_n$ with common density $s$ with respect the Lebesgue
measure on $\R$. For any  $\beta\in[\underline\beta,\bar\beta]$ and any parameter set $\Pc(\underline\beta,\bar\beta)$ given by Proposition~\ref{prop:exist}, there exists a positive
constant $\kappa_\beta$ such that
$$
    \mathcal{R}_n(\H(\beta,\Pc(\underline\beta,\bar\beta))):=\underset{\tilde s}{\inf}\ \underset{s}{\sup}\ \E[d_H^2(s,\tilde{s})] \geq \kappa_\beta\ n^{-\frac{2\beta}{2\beta +1}}
$$
where the supremum (resp. the infimum) is taken over all densities $s$ in $\H(\beta,\Pc(\underline\beta,\bar\beta))$ (resp. over all possible estimators $\tilde s$ of $s$).
\end{thm}

Proposition~\ref{prop:exist} and Theorem~\ref{ThMinoration} are proved in Section~\ref{subs:propexist} and Section~\ref{subs:proofThMinor} respectively.
After establishing Proposition~\ref{prop:exist}, the Hellinger distance and the Kullback-Leibler divergence between two functions of $\mathcal{J}(\beta,D)$ are controlled in Lemma~\ref{propHellinger} and Lemma~\ref{propKullback} respectively. These controls are required to combine a corollary of a Birgé's Lemma \citep[see][]{Birge05} and the so-called Varshamov-Gilbert's Lemma. These last two results can be found in \citet[][see Corollary 2.19 and Lemma 4.7]{Mas07} and are reminded in Appendix~\ref{sect:AppendixMinor}.

\subsection{Adaptive density estimation}  \label{subs:adaptative}
In a non asymptotic model selection approach, the model collection may increase with the sample size $n$, leading to an adaptive procedure.
As it was already explained, the adaptive properties of $\hat s _{\hat m}$ are studied on a range of regularity $[\underline \beta,\bar \beta]$. Preliminary, we fix $ 0 < \underline\beta  < \bar\beta $ and we also choose $a_{\bar \beta} > 1  $ large enough such that
\begin{equation}
    \label{Cdttechn}
    \frac{G_{\bar \beta}  }{a_{\bar \beta}}  \left(  \frac{\ln a_{\bar \beta}} {\ln 2}   + 3  \right)^{3/2}  \leq 1,
\end{equation}
where $G_{\bar \beta}$ is defined in Theorem~\ref{ThApproximation}. The parameters of the Gaussian mixture models $(\S_m)_{m \in \mathcal{M}_n}$ are now specified in order to apply the approximation results provided by Theorem~\ref{ThApproximation}:
\begin{equation*}
    \S_{m} =\left\{x\in\R\mapsto \underset{u=1}{\stackrel{m}{\sum}}p_u \psi_{\sigma_u}(x-\mu_u); \mu_u \in [-\mus(m),\mus(m)], \sigma_u^2 \in [\lambdai(m) , \lambdas(m)], p_u\in[0,1],  \underset{u=1}{\stackrel{m}{\sum}}p_u =1\right\}
    \label{defmodelSmbis}
\end{equation*}
where $\sqrt{\lambdai(m)} := a_{\bar \beta} m^{-1} (\ln m)^{3/2}$,  $\mus(m) = \tilde G _{\bar \beta} | \ln  \sqrt{ \lambdai (m) }|^{1/2}$ and $\lambdas(m)>\lambdai(m)$ for all $m$. Note that the last parameter $\lambdas(m)$ can be taken the same for all $m$ and is denoted $\lambdas$ in the sequel. Since a $n$-sample is observed, it is natural to suppose that the number of mixture components $m$ is less than $n$ and we also assume that the mixtures have at least two components: $\mathcal{M}_n = \{2, \dots,n\}$. Note that when the sample size $n$ increases, mixtures with small component variances and many components $m$ are available in the model collection. This obviously improves the approximation capacity of the Gaussian mixtures.

\begin{thm}~\label{ThAdapt}
Assume that $n \geq 3$ and let $\hat s_{\hat m}$ be the penalized maximum likelihood estimator minimizing the penalized criterion defined in Theorem~\ref{ThMauMic09}.
Then there exists a constant $c_{\underline\beta,\bar\beta}$ such that for all $\beta\in[\underline\beta,\bar\beta]$ and for all $s\in\H(\beta,\Pc(\underline\beta,\bar\beta))$,
$$
\E\left[d_H^{\,2}(s,\hat{s}_{\hat{m}})\right] \leq c_{\underline\beta,\bar\beta}\ (\ln n)^{\frac{5\beta}{2\beta+1}}\ n^{\frac{-2\beta}{2\beta+1}}.
$$
\end{thm}

Theorem~\ref{ThAdapt} shows that the penalized estimator $\hat s_{\hat m}$ is adaptive on the regularity $\beta$ of the density classes defined in Section~\ref{subs:Spaces}, up to a power of $\ln(n)$. This logarithm term is due to the penalty shape given in Theorem~\ref{ThMauMic09}. It is not detected in practice as shown in \citet{MauMic10} and we suspect that it could be removed from the penalty shape. Note that the non parametric Bayesian estimator defined in \cite{KRV} has a similar rate of convergence with a greater power of the logarithm term.

\section{Proof of the approximation result}\label{sect:ProofApprox}
In this section, the density functional space $\H(\beta,\Pc)$ is fixed. To make the proofs and the results easier to read, we use the notation $c_{\beta}$  (resp. $\bar{\sigma}(\beta)$) for denoting constants (resp. upper bound on $\sigma$) that only depends on  $\beta$ and $\Pc $. We also use the notation $c_{\beta,p}$ (resp. $\bar{\sigma}(\beta,p)$) if it also depends on an other parameter $p$. Moreover, we introduce the following notation: For any nonnegative integer $h$, the $h$-fold convolution of the Gaussian kernel $\psi$ is denoted $\psi^{\ast h}$ and for any nonnegative $t$, the $t$-th moment of $\psi^{\ast h}$ is defined by $\nu_{h,t} = \int x^{t} \psi^{\ast h}(x) dx.$
We also denote as $l_j(.)$ the $j$th derivative $\frac{d^j}{dx^j} \ln f (x)$ of $\ln f$ and we consider a subset $A_\sigma$ defined by
\begin{equation*}
    A_{\sigma} := \left\{ x\in\R; \ |l_j(x)|  \leq \B \sigma^{-j} | \ln \sigma |^{-j/2} , j=1\ldots r ,   \ L(x)  \leq \B \sigma^{-\beta} | \ln \sigma |^{-\beta/2} \right\}
\label{Asigmadef}
\end{equation*}
if $\beta>1$ and $A_{\sigma} := \left\{ x\in\R;   \ L(x)  \leq \B \sigma^{-\beta} | \ln \sigma |^{-\beta/2} \right\}$ otherwise.

\subsection{Approximation by a continuous mixture}\label{subs:ApproxP1}
\begin{lem} \label{Lemma1}
Let $\beta>0$ and $k\in\mathbb{N}$ such that $\beta\in(2k,2k+2]$. For all $H>0$, there exists $\bar{\sigma}(\beta,H)>0$ such that for all $\sigma < \bar{\sigma}(\beta,H)$,
for all $f \in \H\left( \beta,\Pc \right)$ and for all $x \in A_{\sigma} $ we have
\begin{equation*}
    (K_{\sigma} f_k) (x) = f(x) \left[1 + R_f(x) O_{\beta,H}( \sigma^{\beta})\right] + O_{\beta,H}(\sigma^{H})
    \label{ControleL1}
\end{equation*}
with $R_f(x) = a_{r+1} L(x)$ if $\beta \leq 1$, and
$$
    R_f(x) = a_{r+1} L(x) + \sum_{j=1}^r   a_{j} |l_{j}(x)| ^{\frac \beta j}
$$
otherwise. In both cases, the $a_j$'s are nonnegative constants that are uniform on $\H\left( \beta,\Pc \right)$. Furthermore, $\bar{\sigma}(\beta,H)$ is a continuous function of $\beta$ and $H$.
\end{lem}
\begin{proof}
Let $H >0$ and  $f \in \H\left( \beta,\Pc \right)$. If $\beta>1$, for all $x$ and $y$ such that $|y-x|\leq \gamma$, there exists $\rho$ such that $|x-\rho|\leq |x-y|$ and
$$
\ln f (y)  =\ln f(x) + \sum_{j=1}^{r} \frac{l_j(x)}{j!} (y-x)^j + \frac{l_r(\rho)-l_r(x)}{r!} (y-x)^{r}.
$$
Then, the smoothness condition (\ref{Cond1}) implies, since $|\rho - x|\leq |y-x|\leq \gamma$, that
\begin{eqnarray*}
\left|\ln f (y)  - \ln f(x) - \sum_{j=1}^{r} \frac{l_j(x)}{j!} (y-x)^j\right| &\leq &  L(x) |\rho-x|^{\beta-r} |y-x|^{r}\\
&\leq &  L(x) |y-x|^{\beta}.
\end{eqnarray*}
Thus we have
\begin{equation}
\ln f (y) \leq \ln f(x) + B(x,y) \label{borne1}
\end{equation}
and
\begin{equation}
\ln f (y)\geq \ln f(x) + \tilde B(x,y) \label{borne2}
\end{equation}
with $ B(x,y) = \underset{j=1}{\stackrel{r}{\sum}} \frac{l_j(x)}{j!} (y-x)^j + L(x) |y-x|^{\beta} $ and $\tilde B(x,y)=\underset{j=1}{\stackrel{r}{\sum}} \frac{l_j(x)}{j!} (y-x)^j - L(x) |y-x|^{\beta}$.
Note that for $\beta\leq 1$, (\ref{borne1}) and (\ref{borne2}) are valid with $B(x,y)=-\tilde B(x,y)=L(x)|x-y|^\beta$.

Let $x \in A_{\sigma} $ and $y \in D_x:= \{ y\in\R; |y-x| \leq k' \sigma | \ln \sigma |^{1/2} \}$ where $k'$, chosen below, has to be identical for all $f \in  \H\left( \beta,\Pc \right)$ and we  also assume that $\sigma$ is small enough to satisfy
\begin{equation}
\label{sig:kprime}
k' \sigma |\ln \sigma |^{1/2} \leq \gamma .
\end{equation}
Then (\ref{borne1}) gives that for all $y \in D_x$, $f(y) \leq f(x) \exp[B(x,y)]$ and thus
\begin{equation}
K_{\sigma} f(x) \leq f(x) \int_{D_x} e^ {B(x,y)} \psi_{\sigma} (y-x) dy +  \int_{D_x^c} f(y) \psi_{\sigma} (y-x) dy .
\label{DxDcx}
\end{equation}
For the sequel, note that, for $x \in A_{\sigma} $ and $y \in D_x$, if $\beta\geq 1$,
\begin{eqnarray*}
\left| B(x,y) \right| & \leq& \sum_{j=1}^r \frac{k'^j}{j!} \sigma^j |\ln \sigma |^{j/2} \B \sigma^{-j} |\ln \sigma |^{-j/2} + k'^{\beta} \sigma^{\beta }  |\ln \sigma |^{\beta /2} \B \sigma^{-\beta } |\ln \sigma |^{-\beta /2 } \\
& \leq& \B \sum_{j=1}^r \frac{k'^j}{j!} + k'^{\beta} := d_1(\beta,k')
\end{eqnarray*}
and thus
\begin{eqnarray}
e^{B(x,y )} & \leq &  \sum_{j=0}^r \frac{1}{j !} B^j(x,y) + \left| B^{r+1}(x,y)  \right| \, \left|\sum_{j \geq 0} \frac{1}{ (j+r+1) !} B^j(x,y) \right| \notag \\
 & \leq &  \sum_{j=0}^r \frac{1}{j !} B^j(x,y) + \left| B^{r+1}(x,y)  \right| \, \left|\sum_{j \geq 0} \frac{1}{ j !} d_1(\beta,k')^j \right| \notag \\
 & \leq &  \sum_{j=0}^r \frac{1}{j !} B^j(x,y) + d_2(\beta,k')  \left| B^{r+1}(x,y) \right| \label{MajoreB}
\end{eqnarray}
with $d_2(\beta,k') = \exp[d_1(\beta,k')]$.

\paragraph{Case k=0 :} We consider that $\beta\in(1,2]$ thus $r=1$. The case $\beta\in(0,1]$ is discussed hereafter.
We have $ B(x,y) =l_1(x) (y-x) + L(x) |y-x|^{\beta} $ and (\ref{MajoreB}) yields
\begin{eqnarray*}
e^{B(x,y) } & \leq &  1 + B (x,y) + d_2(\beta,k') B^2(x,y)   \\
  & \leq & 1+ l_1(x)(y-x) + L(x)  |y-x|^{\beta}  \\
  & & + d_2(\beta,k') \left[  l_1(x)^2  (y-x)^2 + 2 L(x) l_1(x) |y-x |^{ \beta}  (y-x)+ L^2(x) | y-x |^{2 \beta}  \right] \\
  & \leq & 1+ l_1(x)(y-x) + L(x)  |y-x|^{\beta}  \\
  & & + d_2(\beta,k') \left[ (\B k')^{2-\beta} | l_1(x) | ^{\beta}  |  y-x |^{\beta} + 2 \B k'^{\beta} l_1(x)  (y-x)+ \B k'^{\beta} L(x) | y-x |^{\beta}  \right]
  \end{eqnarray*}
since $|l_1(x) (y-x) | \leq \B k'$ and $| L(x)(y-x)^{\beta} | \leq \B k'^{\beta} $. Since $\psi_\sigma$ is symmetric, ${\int_{D_x} (y-x)\psi_{\sigma}(y-x) dy=0}$ and thus,
\begin{eqnarray}
\int_{D_x} e^ {B(x,y)} \psi_{\sigma} (y-x) dy
& \leq  & 1 + L(x) \int_{D_x}  \psi_{\sigma} (y-x) |y-x|^{\beta} dy\notag\\
      & & +  d_2(\beta,k') \left[ (\B k')^{2-\beta} | l_1(x) | ^{\beta} + \B k'^{\beta}L(x) \right] \int_{D_x}  \psi_{\sigma} (y-x)  | y-x |^{\beta} dy\notag\\
& \leq  & 1 + \left\{ d_2(\beta,k') (\B k')^{2-\beta} |l_1(x)| ^{\beta} + [1+d_2(\beta,k') \B k'^{\beta}] L(x) \right\} \nu_{1,\beta} \sigma^{\beta} \notag\\
& \leq  & 1 +   d_3(\beta,k') \left[2 L(x)+| l_1(x) | ^{\beta}  \right]  \sigma^{\beta} . \label{borneDx}
\end{eqnarray}
Let  $k'=k'(2,0,1,H)$ given by Lemma~\ref{Lemma10}, and $\bar \sigma(\beta,H)$ such that (\ref{sig:kprime}) is satisfied. Note that $\bar \sigma(\beta,H)$ can be taken as a continuous function of  $\beta $ and $H$. Then, for the second integral in the  right hand of (\ref{DxDcx}), using (\ref{Cond2b}), it gives
\begin{eqnarray}
\int_{D_x^c} f(y) \psi_{\sigma} (y-x) dy &=& \int_{\R} f(y) \psi_{\sigma} (y-x) \mathds{1}_{ |y-x | \geq k' \sigma | \ln \sigma |^{\frac 1 2} }(y) dy \notag \\
                                         &\leq& M \int_{|u| \geq k' | \ln \sigma |^{\frac 1 2}} \psi (u) du  \notag \\
                                         & \leq & c_{\beta,H} \sigma^H \label{b2}.
\end{eqnarray}
Furthermore, since $k'$ depends on $H$, $d_3(\beta,k')$ can be also rewritten as a constant $c_{\beta,H}$. Finally, (\ref{DxDcx}), (\ref{borneDx}) and (\ref{b2}) give that
$$
(K_{\sigma}f) (x) \leq  f(x) \left\{  1 +   c_{\beta,H} R_f(x)  \sigma^{\beta}  \right\} + c_{\beta,H} \sigma^H
$$
with $R_f(x) =  2 L(x)  + | l_{1}(x)|^\beta$.\\
\\
For $ \beta \in (0,1] $, (\ref{DxDcx}) is still valid with $B(x,y) = L(x) |y-x|^{\beta}$. For $x \in A_{\sigma}$, the first integral in (\ref{DxDcx}) can be treated as in the case $\beta \in (1,2]$: it yields
$$
    \int_{D_x}  \exp(L(x) |x-y|^{\beta} )\psi_{\sigma}(x-y)  dy \leq 1 + c_{\beta}  L(x) \nu_{1,\beta} \sigma^{\beta} .
$$
Using Lemma (\ref{Lemma10}) as before, it gives
$$
    \int_{D_x^c} f(y)  \psi_{\sigma}(x-y)  \leq  c_{\beta,H} \sigma^H
$$
and finally, for all $x \in A_{\sigma}$,
$$
    K_{\sigma} f(x)  \leq f(x) \left[ 1 + c_{\beta,H} L(x) \sigma^{\beta} \right] + c_{\beta,H} \sigma^{\beta} .
$$
A similar lower bound can be shown in the same way, it is proved in the general case further.

\paragraph{Case k=1:} We consider $\beta \in (3,4]$, a similar proof gives the result for $\beta\in(2,3]$. According to (\ref{MajoreB}), for $x \in A_{\sigma} $ and $y \in D_x$,
\begin{equation}
\label{decexpB}
    e^{B(x,y) } \leq 1 + B(x,y) + \frac 1 2 B(x,y)^2 + \frac 1 6 B(x,y)^3 + d_2(\beta,k') B(x,y)^4
\end{equation}
with $ B(x,y) = l_1(x)(y-x) + \frac 1 2 l_2(x) (y-x)^2 +\frac 1 6 l_3(x) (y-x)^3  + L(x)|y-x|^{\beta}$. Thus $\exp[B(x,y)]$ is upper bounded by a linear combination of terms of the form
$$
    \left[ L(x) |x-y|^{\beta} \right] ^{\eta_4}   \prod_{j=1}^3\left[ l_j(x) (y-x)^{ j} \right]^{\eta_j}
$$
with $\sum_{j=1}^4 \eta_j \leq 4$. Let $\mathcal{A}_1(x,y) $ be the sum of such terms for which $\eta_4 \beta + \sum_{j=1}^3 j \eta_j  < \beta$ and let $\mathcal{A}_2(x,y)$ be the others terms and thus $ e^{B(x,y) } \leq \mathcal{A}_1(x,y)+ \mathcal{A}_2(x,y)$. Note that the constant $d_2(k',\beta)$ only appears in the terms of $\mathcal{A}_2(x,y)$. By removing inside $\mathcal{A}_1(x,y) $ all the terms for which the power of $(y-x)$ is an odd integer (since $\int u^t \psi_\sigma(u) du =0$ if $t$ is an odd integer), it yields
\begin{eqnarray}
\int_{D_x}  \mathcal{A}_1(x,y) \psi_{\sigma}(y-x) dy
&= & \int_{D_x} \psi_{\sigma}(y-x)  \left\{ 1 + \frac 1 2  \left[l_1^2(x) +  l_2(x) \right] (y-x)^2 \right\} dy \notag\\
& \leq & 1 +   \frac {\nu_{1,2}}{2} \left\{l_1^2(x) + l_2(x)  \right\} \sigma^2.
\label{DxA1}
\end{eqnarray}
Next, for each term of $\mathcal{A}_2(x,y)$, we have for all $x\in A_{\sigma}$ and all $y \in D_x$
\begin{eqnarray*}
\left[ L(x) |x-y|^{\beta} \right] ^{\eta_4}  \prod_{j=1}^3\left| l_j(x) (y-x)^{ j} \right|^{\eta_j}
& = & \left\{  \left[   L(x)|x-y|^{\beta} \right]^{\eta_4}  \, \prod_{j=1}^3 \left| l_j(x) (x-y)^{j} \right| ^{\eta_j}  \right\}^{ \frac{\eta_4 \beta + \sum_{j=1}^3 j \eta_j - \beta}{\eta_4 \beta + \sum_{j=1}^3 j \eta_j}}   \\
& & \quad \times \,   \left\{   L(x)^{\eta_4} \, \prod_{j=1}^3 \left| l_j(x)\right|^{\eta_j}  \right\}^{ { \frac{\beta}{\eta_4 \beta + \sum_{j=1}^3 j \eta_j}} } \,     |x-y|^{\beta} \\
\end{eqnarray*}
and finally
\begin{eqnarray*}
\left[ L(x) |x-y|^{\beta} \right] ^{\eta_4}  \prod_{j=1}^3\left| l_j(x) (y-x)^{ j} \right|^{\eta_j}
& \leq & c_{\beta} \left\{   L(x)^{\eta_4} \, \prod_{j=1}^3 \left| l_j(x)\right|^{\eta_j}  \right\}^{ { \frac{1}{\eta_4  + \sum_{j=1}^3 \frac j {\beta} \eta_j}} } \,     |x-y|^{\beta}
\end{eqnarray*}
according to Lemma~\ref{LemMajlu}. Note that
\begin{eqnarray*}
\frac{1}{\eta_4  + \sum_{j=1}^3 \frac j {\beta} \eta_j}\ln \left\{   L(x)^{\eta_4} \, \prod_{j=1}^3 \left| l_j(x)\right|^{\eta_j}  \right\}
& = & \frac{1}{\eta_4  + \sum_{j=1}^3 \frac j \beta \eta_j} \left\{  \eta_4 \ln L(x) +  \sum_{j=1}^3  \eta_j  \frac j \beta \ln | l _j  (x)|^{\beta / j} \right\} \\
& \leq  &   \ln\left\{\frac{1}{\eta_4  + \sum_{j=1}^3 \frac j {\beta} \eta_j} \left( \eta_4 L(x) +  \sum_{j=1}^3  \eta_j  \frac j {\beta}  | l _j  (x)|^{\beta / j } \right) \right\}.
\end{eqnarray*}
since the logarithm function is concave.
According to Lemma \ref{LemMajlu}, for each term of $\mathcal{A}_2(x,y)$ we thus have
\begin{eqnarray*}
\left[ L(x) |x-y|^{\beta} \right] ^{\eta_4}  \prod_{j=1}^3\left| l_j(x) (y-x)^{ j} \right|^{\eta_j} & \leq &
\frac{c_{\beta}}{\eta_4  + \sum_{j=1}^3 \frac j {\beta} \eta_j}
|x-y|^{\beta}    \left\{  \eta_4 L(x) +  \sum_{j=1}^3  \eta_j  \frac j {\beta}   |l _{j}(x)|^{\beta / j } \right\}
 \end{eqnarray*}
and then
$$
    | \mathcal{A}_2(x,y) | \leq  c_{\beta,k'} \left\{  a_4 L(x) +  \sum_{j=1}^3  a_j  | l _j  (x)|^{\beta / j }  \right\} |x-y|^{\beta}
$$
where $c_{\beta,k'}$ comes from $d_2(\beta,k')$ in (\ref{decexpB}) and where the $a_j$'s are positive constants that only depend on $\H(\beta,\Pc)$. It leads to
\begin{equation} \label{DxA2}
    \int_{D_x} \left|  \mathcal{A}_2(x,y)  \right| \psi_{\sigma}(y-x) dy  \leq   c_{\beta,k'} R_f(x) \sigma^{\beta}
\end{equation}
where $R_f(x) = a_{4} L(x) + \sum_{j=1}^3   a_{j} |l_{j}(x)|^{\beta/j}$. 
Finally,  (\ref{DxA1}) and (\ref{DxA2}) together yield
$$
    \int_{D_x} e^ {B(x,y)} \psi_{\sigma} (y-x) dy \leq  1 +   \frac {\nu_{2}}{2} \left\{l_1^2(x) + l_2(x)  \right\} \sigma^2 + c_{\beta,k'}  R_f(x)\sigma^{\beta}.
$$
Let  $k'=k'(2,0,1,H)$ given by Lemma~\ref{Lemma10}, and $\bar \sigma(\beta,H)$ such that (\ref{sig:kprime}) is satisfied and where $\bar \sigma(\beta,H)$ can be taken as a continuous function of  $\beta $ and $H$. Next, for all $f \in \H\left( \beta,\Pc \right)$, and all $\sigma < \bar \sigma(\beta,H)$,
$$
    \int_{D_x^c} f(y) \psi_{\sigma} (y-x) dy \leq  c_{\beta,H}\sigma^H .
$$
Finally, for $\sigma < \bar \sigma(\beta,H)$,
$$
    (K_{\sigma}f) (x) \leq   f(x) \left[  1 +  \frac{ \nu_{1,2}}{2}  \left\{l_1^2(x) + l_2(x)  \right\}   \sigma^{2}  +  c_{\beta,H} R_f(x) \sigma^{\beta}\right] +  c_{\beta,H} \sigma^H
$$
and the similar lower bound is obtained in the same way, see the general case further. Thus,
\begin{equation} \label{56}
    (K_{\sigma}f) (x) =   f(x) \left[  1 +  \frac{ \nu_{1,2}}{2}  \left\{l_1^2(x) + l_2(x)  \right\}   \sigma^{2}  +  R_f(x) O_{\beta,H}(\sigma^{\beta})  \right] +  O_{\beta,H}(\sigma^H).
\end{equation}
Now, we need a similar result for $f_1$ instead of $f$. Equation (\ref{56}) depends on the kernel $\psi$
through the values of $\nu_{1,2}$. In fact, it holds for any symmetric kernel $\phi$ such that $\int \phi (x) x^t dx= \nu_{1,t}  < \infty$ and $\int_{|x| > k' |\ln \sigma|^{1/2} } \phi(x)  |x|^{t} dx=  O_{\beta,H}( \sigma^{H})$ when $k'$ is large enough. For $ \psi^{\ast 2}$, these properties follow from Lemma~\ref{Lemma10} : let $k'$ and $\bar \sigma(\beta,H)$ such that (\ref{sig:kprime}) is satisfied and where $\bar \sigma(\beta,H)$ is a continuous function of  $\beta $ and $H$. Thus, denoting $\nu_{2,u}=\int x^u \psi^{\ast 2}(x) dx$, for all $\sigma < \bar \sigma(\beta,H)$,
$$
    (K^2_{\sigma}f) (x) =   f(x) \left[  1 +  \frac{ \nu_{2,2}}{2} \left\{l_1^2(x) + l_2(x)  \right\}  \sigma^{2}  +  R_f(x)O_{\beta,H} (\sigma^\beta)\right] + O_{\beta,H}(\sigma^H).
$$
Now, since $f_1 = 2f -  K_{\sigma} f$ and $ \nu_{2,2} = 2 \nu_{1,2}$, it yields for all $\sigma < \bar \sigma(\beta,H)$ that
$$
    (K_{\sigma}  f_1)(x)   = f(x) \left[  1  +  R_f(x)O_{\beta,H}(\sigma^{\beta} ) \right] +  O_{\beta,H}(\sigma^H).
$$

\paragraph{General case:} Let $ \beta \in (2k,2k+2] $. We give the main ideas of the proof in the general case. According to (\ref{MajoreB}), for $x \in A_{\sigma} $ and $y \in D_x$, $\exp[B(x,y)]$ is upper bounded by a linear combination of terms of the form
$$
    \left[ L(x) |x-y|^{\beta} \right] ^{\eta_{r+1}}   \prod_{j=1}^r\left[ l_j(x) (y-x)^{ j} \right]^{\eta_j}
$$
with $\underset{j=0}{\stackrel{r}{\sum}} \eta_j \leq r+1$. We then decompose $ e^{B(x,y) }$ into $\mathcal{A}_1(x,y) $ and $\mathcal{A}_2(x,y)$ as before. By removing inside $\mathcal{A}_1(x,y) $ all the terms for which the power of $(y-x)$ is an odd integer, it yields
\begin{equation} \label{DxA1gen}
 \int_{D_x}  \mathcal{A}_1(x,y) \psi_{\sigma}(y-x) dy \leq 1 +  \sum_{u=1}^k  \nu_{1,2u} Q_u(x)  \sigma^{2u}
 \end{equation}
where the $Q_u$'s are positive functions that can be expressed in function of $L$ and the $l_{u}$'s. Following the same method as for $\beta \in (3,4]$, it yields
\begin{equation}
    | \mathcal{A}_2(x,y) | \leq  c_{\beta,k'} \left\{  a_{r+1} L(x) +  \sum_{j=1}^r  a_j  | l _j  (x)|^{\beta / j }  \right\} |x-y|^{\beta}
\label{majA2}
\end{equation}
and
\begin{equation}
\label{DxA2gen}
    \int_{D_x} \left|  \mathcal{A}_2(x,y)  \right| \psi_{\sigma}(y-x) dy  \leq   c_{\beta,k'} R_f(x) \sigma^{\beta}
\end{equation}
where $R_f(x) = a_{r+1} L(x) + \underset{j=1}{\stackrel{r}{\sum}} a_{j} |l_{j}(x)|^{\frac \beta j}$. Finally, (\ref{DxA1gen}) and (\ref{DxA2gen}) together yield
$$
    \int_{D_x} e^ {B(x,y)} \psi_{\sigma} (y-x) dy \leq 1 +  \sum_{u=1}^k  \nu_{1,2u} Q_u(x)  \sigma^{2u} + c_{\beta,k'}  R_f(x)\sigma^{\beta}.
$$
Using Lemma \ref{Lemma10}, let $k'$ depending on $H$ and $\bar \sigma(\beta,H)>0$ such that (\ref{sig:kprime}) is satisfied and where $\bar \sigma(\beta,H)$ is a continuous function of  $\beta $ and $H$. For all $\sigma < \bar \sigma(\beta,H)$,
$$
    \int_{D_x^c} f(y) \psi_{\sigma} (y-x) dy \leq  c_{\beta,H} \sigma^H .
$$
Finally, since $k'$ depends on $H$,
$$
    (K_{\sigma}f) (x) \leq   f(x) \left[  1 +  \sum_{u=1}^k  \nu_{1,2u} Q_u(x)  \sigma^{2u}  + c_{\beta,H} \,  R_f(x) \sigma^{\beta}  \right] +  c_{\beta,H}  \sigma^H
$$
and the similar lower bound is obtained in the same way (see further). Thus,
$$
    (K_{\sigma}f) (x) =   f(x) \left[  1 +  \sum_{u=1}^k  \nu_{1,2u} Q_u(x)  \sigma^{2u}  +  R_f(x) O_{\beta,H}(\sigma^{\beta} ) \right] + O_{\beta,H}(\sigma^H).
$$
Now, we need a similar result for $f_k$ instead of $f$. According to Lemma~\ref{Lemfk},
\begin{equation}
    f_k = \sum_{i=0}^k (-1)^i \left(\substack{k+1 \\ i+1} \right) K_{\sigma}^i f .
    \label{decfk}
\end{equation}
For all $h\leq k$, the same method can be applied with $\psi^{\ast h}$ instead of $\psi$ and it yields with the same functions $Q_u$ and $R_f$:  for $\sigma < \bar \sigma (\beta,H)$,
$$
    (K^h_{\sigma}f) (x) =   f(x) \left[  1 +  \sum_{u=1}^k  \nu_{h,2u} Q_u(x)  \sigma^{2u}  +  R_f(x) O_{\beta,H}(\sigma^{\beta} ) \right] + O_{\beta,H}(\sigma^H).
$$
According to (\ref{decfk}), for $\sigma < \bar \sigma (\beta,H)$
\begin{eqnarray*}
K_{\sigma} f_k(x) &=& \sum_{i=0}^k (-1)^i \left(\substack{k+1 \\ i+1} \right) K_{\sigma}^{i+1} f(x)  \\
                  &=& \sum_{j=1}^{k+1} (-1)^{j+1} \left(\substack{k+1\\j} \right) K_{\sigma}^{j} f(x)  \\
                  &=& \sum_{j=1}^{k+1} (-1)^{j+1} \left(\substack{k+1 \\ j} \right) f(x) \left[  1 +  \sum_{u=1}^{k}  \nu_{j,2u} Q_u(x)  \sigma^{2u}  +  R_f(x) O_{\beta,H}(\sigma^{\beta} ) \right] +     O_{\beta,H}(\sigma^H) \\                
 &=& f(x) \left\{ 1 +   \sum_{u=1}^{k+1}  \left[ \sum_{j=1}^{k+1} (-1)^{j+1} \left(\substack{k+1 \\j} \right)  \nu_{j,2u}  \right] Q_u(x)  \sigma^{2u}  +  R_f(x) O_{\beta,H}(\sigma^{\beta} )   \right\} +  O_{\beta,H}(\sigma^H)\\
\end{eqnarray*}     
and then             
$$ K_{\sigma} f_k(x) = f(x) \left\{ 1 +  R_f(x) O_{\beta,H}(\sigma^{\beta} )   \right\} +  O_{\beta,H}(\sigma^H) $$
since $\sum_{i=0}^k (-1)^i \left(\substack{k+1 \\ i+1} \right) = 1$ and $\sum_{j=1}^{k+1} (-1)^{j+1} \left(\substack{k+1 \\j} \right)  \nu_{j,2u}=0$ according to Lemma~\ref{LemmeSomme}.

To complete this proof, we give the method for obtaining the lower bound in the general case.
Using (\ref{borne2}) and proceeding in the same way as for the upper bound, it yields
\begin{eqnarray}
K_{\sigma} f(x)
&\geq& f(x) \int_{D_x} e^ {\tilde B (x,y)} \psi_{\sigma} (y-x) dy \notag \\
&\geq&  f(x) \int_{D_x}  \left[ 1 + \sum_{j=1}^r \frac{1}{j !} \tilde B ^j (x,y) - d_2(\beta,k') \left| \tilde B ^{r+1} (x,y) \right| \right] \psi_{\sigma} (y-x) dy \notag \\
&\geq&  f(x) \int_{D_x} \mathcal{A}_1(x,y)   \psi_{\sigma} (y-x)  dy  + f(x) \int_{D_x} \mathcal{A}_2(x,y)   \psi_{\sigma} (y-x)    dy \label{A1A2}
\end{eqnarray}
where $\mathcal{A}_1(x,y) $ (resp. $\mathcal{A}_2(x,y)$) contains the terms which powers are less than $\beta$ (resp. larger than $\beta$). For the first integral,
\begin{eqnarray*}
\int_{D_x} \mathcal{A}_1(x,y)   \psi_{\sigma} (y-x)  dy
&=&    \int_{\R} \mathcal{A}_1(x,y)   \psi_{\sigma} (y-x)  dy  -  \int_{D_x^c} \mathcal{A}_1(x,y)   \psi_{\sigma} (y-x)  dy  \\
&\geq& 1 + \sum_{u=1}^k \nu_{1,2u} Q_u(x) \sigma^{2u} - \left| \int_{D_x^c} \mathcal{A}_1(x,y)   \psi_{\sigma} (y-x)  dy \right|.
\end{eqnarray*}
Note that $\mathcal{A}_1 (x,y)$ is  a linear combination of terms of the form  $\prod_{j=1}^r\left[ l_j(x) (y-x)^{ j} \right]^{\eta_j}$, where $\sum _{j=1}^r j \eta_j  $ is even.
Since $x \in A_{\sigma}$, and since $|x-y| \leq 1$, we can find $h >0 $ and a constant $c_{\beta}$ that only depends on $\beta$ such that
$$
    \left| \mathcal{A}_1(x,y)  \right| \leq c_{\beta} \sigma^{-h}  (x-y)^{2}.
$$
Finally,
$$
    \left| \int_{D_x^c} \mathcal{A}_1(x,y)   \psi_{\sigma} (y-x)  dy \right|  \leq  c_{\beta} \sigma^{-h} \int_{D_x^c} (x-y)^{2}  \psi_{\sigma}(y-x)  dy.
$$
Then we apply Lemma \ref{Lemma10} with $H' > h + H$  and it gives that
$$
    \left| \int_{D_x^c} \mathcal{A}_1(x,y)   \psi_{\sigma} (y-x)  dy \right|  \leq  c_{\beta,H}\ \sigma^{H}.
$$
To find a lower bound for the second integral in (\ref{A1A2}), we note that according to (\ref{majA2}),
$$
    | \mathcal{A}_2(x,y) | \geq -   c_{\beta,H}  \left\{  a_{r+1} L(x) +  \sum_{j=1}^r  a_j  | l _j  (x)|^{\beta / j }  \right\} |x-y|^{\beta}
$$
and thus
\begin{equation*}
    \int_{D_x} \left|  \mathcal{A}_2(x,y)  \right| \psi_{\sigma}(y-x) dy  \leq    c_{\beta,H}  R_f(x) \sigma^{\beta}
\end{equation*}
where $R_f$ is defined as before. We finally obtain that for $\sigma < \bar \sigma (\beta,H)$,
$$
    K_{\sigma}f (x) \geq   f(x) \left[  1 +  \sum_{u=1}^k  \nu_{1,2u} Q_u(x)  \sigma^{2u}  +    c_{\beta,H}  \,  R_f(x) \sigma^{\beta}  \right] +  c_{\beta,H}  \sigma^H.
$$
\end{proof}

For a density $f$ belonging to $\H(\beta,\Pc)$, Lemma~\ref{Lemma1} shows that the convolution $K_\sigma f_k$ is close to $f$ on a subspace of $\R$ where the derivative functions of $\ln f$ and $L$ are efficiently controlled. Furthermore, the control on the difference $K_\sigma f_k - f$ is uniform over $\H(\beta,\Pc)$, which is required to upper bound the Kullback-Leibler divergence between $f$ and $K_\sigma f_k$. Thus $K_\sigma f_k$ seems to be a good candidate to approximate the density function $f$. Nevertheless, the function $f_k$ is not a density function: Its integral over $\R$ is equal to 1 (see Lemma~\ref{Lemfk}) but it can take negative values. To remedy this problem, \citet{KRV} define a density function $h_k$ as follows:
Considering the subspace
$$
    J_{\sigma,k}  =   \left\{ x\in\R;\ f_k(x) > \frac 1 2 f(x) \right\},  
$$
the following positive function is defined
$$
    \forall x\in\R,\ g_k(x) =  f_k(x) \mathds{1}_{J_{\sigma,k}}(x) + \frac 1 2 f(x) \mathds{1}_{J_{\sigma,k}^c}(x) 
$$
and it is normalized to obtain a density function
\begin{equation}
    \forall x\in\R,\ h_k(x) =  \frac{g_k(x)}{\int g_k (u) du}.
    \label{W20b}
\end{equation}
Note that the constant $1/2$ is arbitrary in the definition of $J_{\sigma,k}$, any other number of $(0,1)$ could be used.

Now, the result of Lemma~\ref{Lemma1} has to be extended for the convolution $K_\sigma h_k$. For this purpose, the integral of $K_\sigma^t f$ for all nonnegative integers $t\leq k$ is controlled over $A_\sigma^c$ and $E_\sigma^c$ where $A_\sigma$ is defined by (\ref{Asigmadef}) and $E_{\sigma} = \{x\in\R;\ f(x) \geq \sigma^{H_1}\}$ with $H_1> 4 \beta$.

\medskip
\begin{lem}\label{Lemme2}
Let $\beta>0$ and $k\in\mathbb{N}$ such that $\beta\in(2k,2k+2]$. There exists $\bar \sigma (\beta,H_1) >0 $ such that for all $\sigma \leq \bar \sigma (\beta,H_1) $, for all $f \in \H(\beta,\Pc)$ and for all nonnegative integers $t \leq k$,
\begin{equation}\label{W23g}
    \int_{A_{\sigma}^c} \left( K_{\sigma}^t f \right)(x) dx = O_{\beta}(\sigma^{2 \beta})
\end{equation}
and
\begin{equation}\label{W23d}
    \int_{E_{\sigma}^c} \left( K_{\sigma}^t f \right)(x) dx = O_{\beta,H_1}(\sigma^{2 \beta}) .
\end{equation}
Furthermore, for $\sigma \leq \bar \sigma (\beta,H_1) $, $A_{\sigma} \cap E_{\sigma} \subset J_{\sigma,k}$ and
\begin{equation}\label{W24}
    \int_{\mathbb{R}} g_k(x)  dx = 1+ O_{\beta,H_1} (\sigma^{2 \beta}).
\end{equation}
Thus, for all $H>0$, there exists $\bar \sigma (\beta,H_1,H) >0 $ such that for all $\sigma \leq \bar \sigma (\beta,H_1,H) $ and for all $x \in A_{\sigma} \cap E_{\sigma}$,
\begin{equation}\label{W25}
    \left|(K_{\sigma} h_k)(x)-f(x)\right| = f(x) R_f(x) O_{\beta,H_1,H}( \sigma^{\beta}) +  O_{\beta,H_1,H}(\sigma^{H}).
\end{equation}
Furthermore,  $\bar \sigma (\beta,H_1) $ and  $\bar \sigma (\beta,H_1,H)$ are both continuous functions of $\beta$, $H_1$ and $H$ for the last one.
\end{lem}

\begin{Rem} \label{Rem:W24bis} The left term in (\ref{W24}) does not depend on $H_1$ whereas the right term does. Indeed, the presence of $H_1$ here is only technical and by choosing for instance $H_1 = 4\beta +1$, it gives that there exists a  positive constant $\bar \sigma (\beta)$, continuous in $\beta$ such that for all $\sigma \leq \bar \sigma (\beta) $,
\begin{equation}
\label{W24bis}
   \int_{\mathbb{R}} g_k(x) dx = 1+ O_{\beta} (\sigma^{2 \beta}).
\end{equation}
\end{Rem}

\begin{proof}
For $\delta \leq 1 $ to be chosen further, let
$$
    A_{\sigma,\delta} := \{ x \in \R, | l_j(x) | \leq \delta \B \sigma^{-j} | \ln \sigma|^{-j/2} ,  \forall j \in 1 \dots r ,  L(x) \leq \delta\B\sigma^{-\beta} |\ln \sigma|^{- \beta / 2} \}
$$
if $\beta>1$ and let
$$
    A_{\sigma,\delta} := \{ x \in \R,  L(x) \leq \delta\B\sigma^{-\beta} |\ln \sigma|^{- \beta / 2} \}
$$
otherwise. Note that for all $\delta \leq 1 $, $A_{\sigma,\delta} \subset A_{\sigma}$. In the sequel we assume that $\beta>1$, the proof being easily adapted for $\beta<1$.

\noindent \textbf{Proof of (\ref{W23g}):}\\
\noindent $\bullet$ Case $t = 0$: If $X \sim f$, then
\begin{eqnarray*}
\int_{A_{\sigma}^c} \left( K_{\sigma}^0 f \right)(x) dx
&=&    \int_{A_{\sigma}^c} f (x) dx \\
&\leq& \sum_{j=1}^{r} P(\ |l_j(X)| >   (\delta\B) \sigma^{-j} |\ln \sigma |^{-j/2}  )  + P( |L(X) | >  (\delta\B) \sigma^{-\beta}  |\ln \sigma |^{-\beta/2})\\
&\leq& \sum_{j=1}^{r} P(\ |l_j(X)|^{\frac{2\beta+\varepsilon}{j} } >  (\delta\B)^{\frac{2\beta+\varepsilon}{j}} \sigma^{-2\beta-\varepsilon} |\ln \sigma |^{- \frac{2\beta+\varepsilon}{2} } ) \\
 & &     +   P(|L(X)|^{\frac{2 \beta + \varepsilon}{\beta}} >  (\delta\B)^{\frac{2\beta+\varepsilon}{\beta}} \sigma^{-2 \beta- \varepsilon}  |\ln \sigma |^{-\frac{2\beta+\varepsilon}{2}})\\
&\leq& \sum_{j=1}^{r} P(\ |l_j(X)|^{\frac{2\beta+\varepsilon}{j} } >  (\delta\B)^{\frac{2\beta+\varepsilon}{j}} \sigma^{-2\beta} )	 +  P( |L(X)|^{\frac{2 \beta + \varepsilon}{\beta}} >  (\delta\B)^{\frac{2\beta+\varepsilon}{\beta}} \sigma^{-2 \beta})
\end{eqnarray*}
since $\sigma^{-\varepsilon} |\ln \sigma|^{-\frac{2\beta + \varepsilon}{2}} > 1$ for $\sigma$ small enough (say $\sigma \leq \bar \sigma (\beta)$). Then, Markov Inequality together with (\ref{Cond2a}) gives
$$
    P(\ |l_j(X)|^{\frac{2\beta+\varepsilon}{j} } > (\delta\B)^{\frac{2\beta+\varepsilon}{j}} \sigma^{-2\beta} )\leq (\delta\B)^{-\frac{2\beta+\varepsilon}{j}}\sigma^{2 \beta} \mathbb{E}[|l_j(X)|^{\frac{2\beta + \varepsilon}{j}}]
    \leq c_\beta  \sigma^{2 \beta}	
$$
and
$$
    P( |L(X)|^{\frac{2 \beta + \varepsilon}{\beta}} > (\delta\B)^{\frac{2\beta+\varepsilon}{\beta}} \sigma^{-2 \beta} ) \leq c_\beta  \sigma^{2\beta}.
$$
Finally,
$$
    \int_{A_{\sigma}^c} \left( K_{\sigma}^0 f \right)(x) dx \leq c_\beta \sigma^{2 \beta}.
$$

\noindent $\bullet$ Case $t = 1$: Let $X \sim f$ and $ U \sim \psi$, then $X + \sigma U \sim K_{\sigma} f$. By applying Lemma~\ref{Lemma10} with $H = 2\beta$, let $k'$ depending on $\beta$ and $\bar \sigma (\beta)$ such that for all $\sigma \leq \bar \sigma (\beta)$, $ k'  \sigma | \ln \sigma |^{1/2} \leq \gamma  $ and $P( |U| > k' | \ln \sigma |^{1/2} ) \leq c_{\beta} \sigma^{2 \beta}$. Then,
\begin{eqnarray}
\int_{ A_{\sigma}^c } K_{\sigma} f (x) dx
&=&    P(X + \sigma U \in A_{\sigma}^c)   \notag \\
&=&    P(X + \sigma U \in A_{\sigma}^c  \cap |U|  \leq  k' | \ln \sigma|^{1/2} ) +P(X + \sigma U \in A_{\sigma}^c  \cap |U|  > k' | \ln \sigma|^{1/2} )  \notag \\
\end{eqnarray}
and
\begin{eqnarray}
\int_{ A_{\sigma}^c } K_{\sigma} f (x) dx &\leq& P(X + \sigma U \in A_{\sigma}^c  \cap  X \in A_{\sigma,\delta} \cap  |U|  \leq  k' | \ln \sigma|^{1/2} )  \notag  \\
& &  \; + \; P(X + \sigma U \in A_{\sigma}^c  \cap  X \in A_{\sigma,\delta}^c \cap  |U|  \leq k' | \ln \sigma|^{1/2} )  +  P( |U|  >  k' | \ln \sigma|^{1/2} ) \notag  \quad \quad\\
&\leq& P(X + \sigma U \in A_{\sigma}^c  \cap  X \in A_{\sigma,\delta} \cap  |U|  \leq  k' | \ln \sigma|^{1/2} ) + P(X \in A_{\sigma,\delta}^c )+   c_{\beta} \sigma^{2 \beta}  \label{decm1} .
\end{eqnarray}
The second term in (\ref{decm1}) can be shown to be bounded by a multiple of $\sigma^{2 \beta}$ in the same manner as for $t=0$ for $\sigma \leq \bar \sigma(\beta)$. We now show that for $\sigma$ small enough, the first term in (\ref{decm1}) is zero for every function $f \in \H(\beta,\Pc)$. On the one hand, according to (\ref{Cond1}) there exists $y \in [X,X + \sigma U]$ such that
$$  l_j(X + \sigma U)  = \sum_{u= 0} ^{r-j-1} \frac{l_{j+u}( X ) }{u !} (\sigma U) ^{u}  + \frac{l_{r}(y ) }{(r-j)!} (\sigma U) ^{r-j} . $$
If $  X \in A_{\sigma,\delta} $ and $ |U|  \leq  k' | \ln \sigma|^{1/2} $ it yields
\begin{eqnarray*}
| l_j(X + \sigma U) |& \leq & \sum_{u= 0} ^{r-j-1} \left| \frac{l_{j+u}( X ) }{u !} \right| | \sigma U |^{u}  + \left| \frac{l_{r}(y ) - l_{r}(X )}{(r-j)!}\right| |\sigma U| ^{r-j}  +  \left| \frac{l_{r}(X )}{(r-j)!} \right| | \sigma U| ^{r-j}  \\
& \leq & \sum_{u= 0} ^{r-j} \left| \frac{l_{j+u}( X ) }{u !} \right|  | \sigma U|  ^{u} + \frac{r!}{(r-j)!}   L(X) |y-X|^{\beta-r} \left| \sigma U \right|^{r-j}  \\
& \leq & \sum_{u= 0} ^{r-j} \left| \frac{l_{j+u}( X ) }{u !} \right|  | \sigma U|  ^{u} + \frac{r!}{(r-j)!}   L(X)\left| \sigma U \right|^{\beta-j} \\
& \leq & \sum_{u= 0} ^{r-j}  \frac{1}{u !} \delta \B \left( \sigma | \ln \sigma |^{1/2} \right)^{-(u+j)} \left(  \sigma k' | \ln \sigma |^{1/2} \right)^u  + \frac{r!}{(r-j)!}    \delta \B \left( \sigma | \ln \sigma |^{1/2} \right)^{-\beta}  \left(  \sigma k' | \ln \sigma |^{1/2} \right)^{\beta-j} .
\end{eqnarray*}
And thus for $\delta$ small enough, $| l_j(X + \sigma U) | \leq  \B \left(  \sigma  | \ln \sigma |^{1/2} \right)^{-j}$ for all $j \in \{1\dots r \}$. Since ${X + \sigma U \in A_{\sigma}^c}$, this means that
\begin{equation}
    \label{Contrad1}
    L(X+\sigma U) >\B(\sigma |\ln \sigma |^{1/2})^{-\beta}.
\end{equation}
On the other hand, let $\eta  = \max |z_i |$ where the $z_i$'s are the roots of $L$. Suppose that $\mbox{deg}(L)=q$, for $j = 1, \dots,q,$  $| L^{(j)}(x) | / | L(x)| \rightarrow 0$ when $|x|$ tends to infinity. Consequently, since $L$ does not vanish out of $[-\eta,\eta]$, there exists $c>0$ only depending on $L$ such that if $|x| > \eta +1$, then $| L^{(j)}(x)| \leq c |L(x)|$. If $|X| > \eta +1 $, then
\begin{eqnarray}
| L(X + \sigma U) |
& \leq & L(X) + \sum_{j= 1} ^{q} \left| \frac{L^{(j)}(X)}{j !} \right| | \sigma U |^{j}  \label{devL} \\
& \leq &  \delta \B \left( \sigma | \ln \sigma |^{1/2} \right)^{-\beta}
          + c |L(X)| \sum_{j= 1} ^{q}  \frac{1}{j !}  \left( \sigma k' | \ln \sigma  |^{1/2} \right)^{j}  \notag \\
& \leq & \delta \B \left( \sigma | \ln \sigma |^{1/2} \right)^{-\beta}  + c_{\beta} \delta \B \left( \sigma | \ln \sigma  |^{1/2} \right)^{-\beta+1}  \notag \\
& \leq & 2 \delta \B \left( \sigma | \ln \sigma |^{1/2} \right)^{-\beta}    \notag
\end{eqnarray}
for $\sigma \leq \bar \sigma (\beta) $ where $\bar \sigma (\beta) $ can be chosen as a continuous function of $\beta$. It then leads to a contradiction with (\ref{Contrad1}) for $\delta$ chosen small enough and thus $\P( X + \sigma U \in A_{\sigma}^c  \cap  X \in A_{\sigma,\delta} \cap  |U|  \leq  k' | \ln \sigma|^{1/2}  \cap |X| > \eta +1 ) = 0 . $ Next, let $\bar{L} : = \underset{j =0 ... r}{\max} \underset{|x| \leq \eta +1  }{\sup} | L^{(j)}(x) | $. If  $|X| \leq \eta +1 $, (\ref{devL}) implies that
\begin{eqnarray*}
| L(X + \sigma U) |
& \leq &  \delta \B \left( \sigma | \ln \sigma |^{1/2} \right)^{-\beta}
          + c \bar{L} \sum_{j= 1} ^{q}  \frac{1}{j !}  \left( \sigma k'  | \ln \sigma  |^{-1/2} \right)^{j}\\
& \leq & \delta \B \left( \sigma | \ln \sigma |^{1/2} \right)^{-\beta}  + c_{\beta} \sigma | \ln \sigma  |^{-1/2}\\
& \leq & 2 \delta \B \left( \sigma | \ln \sigma |^{1/2} \right)^{-\beta}
\end{eqnarray*}
for $\sigma \leq \bar \sigma (\beta) $ where $\bar \sigma (\beta) $ can be chosen as a continuous function of $\beta$, which also gives a contradiction with (\ref{Contrad1}). Thus $\P( X + \sigma U \in A_{\sigma}^c  \cap  X \in A_{\sigma,\delta} \cap  |U|  \leq  k' | \ln \sigma|^{1/2}  \cap |X| \leq  \eta +1 ) = 0 $. Finally, for $\sigma  \leq \bar \sigma (\beta) $, $\P( X + \sigma U \in A_{\sigma}^c  \cap  X \in A_{\sigma,\delta} \cap  |U|  \leq  k' | \ln \sigma|^{1/2})=0$ and (\ref{decm1}) gives that
$$
    \int_{A_\sigma^c} K_\sigma f(x) dx = O_\beta(\sigma^{2\beta}).
$$

\noindent $\bullet$ Case $t \geq 2$. The same method as before can be applied by assuming $X \sim K^{t-1}_{\sigma} f$ and $ U \sim \psi$. Similarly, $\int_{ A_{\sigma}^c } K_{\sigma}^{t} f (x) dx$ can be decomposed into three terms as in (\ref{decm1}). Two of them are $O_{\beta}(\sigma^{2\beta})$ and the remaining term is zero for $\delta$ small enough.

\paragraph{Proof of (\ref{W23d}):~\\}
\noindent $\bullet$ Case $t = 0$: According to Condition (\ref{Cond2b}),
\begin{eqnarray*}
\int_{E_{\sigma}^c}  f (x) dx & \leq & \sigma^{H_1/2} \int_{E_{\sigma}^c} \sqrt{f(x)} dx\\
                              & \leq & \sigma^{H_1/2} \int_{\mathbb{R}} \sqrt{M \pi^{-\frac 1 2}} \exp\left(-x^2\right) dx \\
                              & \leq & \sigma^{2\beta} M^{\frac 1 2} \pi^{\frac 1 4}
\end{eqnarray*}
since $H_1 > 4 \beta$.

\noindent $\bullet$ Case $t = 1$: We have
$$
    \int_{E_{\sigma}^c}  K_{\sigma} f (x) dx = \int_{E_{\sigma}^c \cap A_{\sigma}}  K_{\sigma} f (x) dx + \int_{E_{\sigma}^c \cap A_{\sigma}^c} K_{\sigma} f (x) dx
$$
where the second integral is less than $ \int_{A_{\sigma}^c} K_{\sigma} f (x) dx $ which is $O_{\beta}(\sigma^{2\beta})$ for $\sigma \leq \bar \sigma (\beta)$ according to (\ref{W23g}). For $\delta \leq 1$ to be chosen further uniformly on $\H(\beta,\Pc)$, consider the set
$ E_{\sigma,\delta} = \{ x\in\mathbb{R}; f(x)\geq \sigma^{\delta H_1}\}$. Let $X \sim f$ and $U \sim \psi$. By applying Lemma~\ref{Lemma10} as before with $H = 2\beta$, let $k'$ depending on $\beta$ and $\bar \sigma (\beta)$ such that for all $\sigma \leq \bar \sigma (\beta)$, $ k'  \sigma | \ln \sigma |^{1/2} \leq \gamma  $ and
$P( |U| > k' | \ln \sigma |^{1/2} ) \leq c_{\beta} \sigma^{2 \beta}$. Then,
\begin{eqnarray}
\int_{E_{\sigma}^c \cap A_{\sigma}}  K_{\sigma} f (x) dx
& = &   \P \left( X + \sigma U \in  E_{\sigma}^c \cap A_{\sigma} \right) \notag \\
&\leq&  \P \left( X + \sigma U \in  E_{\sigma}^c \cap A_{\sigma}  \; ; \; |U| \leq k' |\ln \sigma|^{1/2} \right) + \P \left(  |U| >  k' |\ln \sigma|^{1/2} \right)  \notag \\
\end{eqnarray}
\begin{eqnarray}
\int_{E_{\sigma}^c \cap A_{\sigma}}  K_{\sigma} f (x) dx & \leq  &  \P \left( X + \sigma U \in  E_{\sigma}^c \cap A_{\sigma}  \; ; \; |U| \leq k' |\ln \sigma|^{1/2} \; ; \; X \in A_{\sigma} \right) +  \P \left(X \in A_{\sigma}^c \right) \notag \\
& & \; + \; \P \left(  |U| >  k' |\ln \sigma|^{1/2} \right) \notag \\
& \leq  &  \P \left( X + \sigma U \in  E_{\sigma}^c  \; ; \; |U| \leq k' |\ln \sigma|^{1/2} \; ; \; X \in A_{\sigma} \cap E_{\sigma,\delta} \right) \label{Ligne1} \\
& & \; + \; \P \left(X \in  E_{\sigma,\delta}^c \right)  + \P \left(X \in A_{\sigma}^c \right) + c_{\beta} \sigma^{2 \beta}. \label{Ligne2}
\end{eqnarray}
According to (\ref{W23g}), the second term in (\ref{Ligne2}) is $O_{\beta}(\sigma^{2\beta})$ for $\sigma \leq \bar \sigma (\beta)$. The first term in (\ref{Ligne2}) can be bounded as previously for $t=0$, leading to the condition
\begin{equation*}
    \delta H_1 \geq 4 \beta
\end{equation*}
and thus we choose $\delta \in (0,1)$ to satisfy this last condition. It remains to control the probability given in (\ref{Ligne1}). On the one hand, since $X + \sigma U \in  E_{\sigma}^c$ and $X \in E_{\sigma,\delta}$, 
\begin{equation}
    \label{Contardic1}
    | \ln f(X + \sigma U)  - \ln f(X ) | >  (1-\delta) H_1 |\ln\sigma | .
\end{equation}
On the other hand, since $X \in A_{\sigma} $ and $|U| \leq k' |\ln \sigma|^{1/2}$,
\begin{eqnarray*}
| \ln f(X + \sigma U)  - \ln f(X ) | &  \leq & \sum_{j=1}^r  \left| \frac{l_j(X)}{j!} \sigma^j U^j \right| + L(X)  \left| \sigma U \right|^{\beta } \\
                                 &  \leq &   \B e^{k'} + \B k'^\beta := d_1(\beta,k').
\end{eqnarray*}
This is in contradiction with (\ref{Contardic1}) for $\sigma \leq \exp \left(- \frac{d_1(\beta,k')}{(1-\delta)H_1}\right)$ and then (\ref{Ligne1}) is zero for $\sigma \leq \bar \sigma (\beta,H_1)$ where $\bar \sigma (\beta,H_1)$ can be chosen continuous.

\noindent $\bullet$ Case $t \geq 2$. We follow the same proof as before:
$$
    \int_{E_{\sigma}^c} K_{\sigma}^t f (x)  dx = \int_{E_{\sigma}^c \cap A_{\sigma}}  K_{\sigma}^t f (x) dx + \int_{E_{\sigma}^c \cap A_{\sigma}^c} K_{\sigma}^t  f (x)dx
$$
where the second integral is less than $ \int_{A_{\sigma}^c} K_{\sigma}^t f (x) dx $ which is $O_{\beta}(\sigma^{2\beta})$ according to (\ref{W23g}). Let $X \sim f$ and $U_1, \dots,U_t \sim \psi$. By applying Lemma~\ref{Lemma10} as before with $H = 2\beta$, let $k'$ depending on $\beta$ and $\bar \sigma (\beta)$ such that for all $\sigma \leq \bar \sigma (\beta)$, $ k'  \sigma | \ln \sigma |^{1/2} \leq \gamma  $ and $ P( |U| > k' | \ln \sigma |^{1/2} ) \leq c_{\beta} \sigma^{2 \beta} $. Then,
\begin{eqnarray}
\int_{E_{\sigma}^c \cap A_{\sigma}}  \hskip-0.5cm  K_{\sigma}^t f (x) dx
&=&     \P \left( X + \sigma U_1 + \dots + \sigma U_t \in  E_{\sigma}^c \cap A_{\sigma} \right) \notag \\
&\leq&  \P \left( X + \sigma U_1+ \dots + \sigma U_t \in  E_{\sigma}^c \cap A_{\sigma}  \; ; \;  \forall j \ |U_j| \leq k' |\ln \sigma|^{1/2} \right) + \sum_{j=1}^t \P \left(  |U_j| >  k' |\ln \sigma|^{1/2} \right)  \notag \\
&\leq&  \P \left( X + \sigma U_1+ \dots + \sigma U_t \in  E_{\sigma}^c \cap A_{\sigma}  \; ; \;  \forall j \ |U_j| \leq k' |\ln \sigma|^{1/2} \; ; \; X \in A_{\sigma} \cap E_{\sigma,\delta} \right) \label{Ligne1m} \\
& & + \quad  \P \left(X \in  E_{\sigma,\delta}^c \right)  + \P \left(X \in A_{\sigma}^c \right)  +\ c_{\beta}.
\label{Ligne2m}
\end{eqnarray}
According to (\ref{W23g}), the second term in (\ref{Ligne2m}) is $O_{\beta}(\sigma^{2\beta})$ for $\sigma \leq \bar \sigma (\beta)$, as well as the first term if $\delta H_1 \geq 4 \beta $.  We thus choose $\delta \in (0,1)$ to satisfy this last condition. As before, we check that the probability given in  (\ref{Ligne1m}) is $0$. On the one hand, since $X + \sigma U_1 + \dots + \sigma U_t \in  E_{\sigma}^c$ and $X \in E_{\sigma,\delta}$, 
\begin{equation}
    \label{Contardic1gen}
    | \ln f(X + \sigma U_1 + \dots + \sigma U_t)  - \ln f(X ) | >  (1-\delta) H_1 |\ln\sigma| .
\end{equation}
On the other hand, since $X \in A_{\sigma} $ and for all $j$, $|U_j| \leq k' |\ln \sigma|^{1/2}$, then
\begin{eqnarray*}
 | \ln f(X  + \sigma U_1 + \dots + \sigma U_t)  - \ln f(X ) | &  \leq & \sum_{j=1}^r  \left| \frac{l_j(X)}{j!} \sigma^j  \left( \sum_{i=1}^t U_i\right)^j \right| + L(X)   \left| \sigma \sum_{i=1}^t U_i \right|^{\beta} \\
   &  \leq & \sum_{j=1}^r  \left| \frac{l_j(X)}{j!}  \left( \sigma k'  t | \ln \sigma | ^{1/2} \right)^j \right| +  L(X)  \left( \sigma t k'  | \ln \sigma | ^{1/2} \right)^{\beta} \\ \\
  &  \leq & \B e^{tk'} +\B(tk')^\beta .
  \end{eqnarray*}
This in contradiction with (\ref{Contardic1gen}) for $\sigma \leq \bar \sigma (\beta,H_1)$ and finally (\ref{Ligne1m}) is zero.

\paragraph{Proof of $\mathbf{E_{\sigma} \cap A_{\sigma} \subset J_{\sigma,k}}$:~\\}
For $\beta \leq 2$, the inclusion  is obvious since $f_0 = f $ and thus $J_{\sigma,k} = \R$. To prove the case $\beta >2$, we show by induction on $u \in \N$, $1 \leq u \leq k $ that for every $h \in (0,1)$, there exists a continuous function $\bar \sigma(\beta,H_1,h)$ such that for all $\sigma \leq \bar \sigma(\beta,H_1,h)$, for all $f \in \H(\beta,\Pc)$ and all $x \in E_{\sigma} \cap A_{\sigma}$,
\begin{equation}\label{fugeqf}
    f_u(x) \geq \left(1-\frac h 2\right)f(x) .
\end{equation}
\noindent $\bullet$ Let $u= 1$ and $L^{(1)}$ defined by
$$
    L^{(1)}(x) = \sum_{j=1}^r \frac{|l_j(x) |}{j!}  \left(k' \sigma | \ln \sigma|^{1/2} \right)^{j-1} + L(x)\left(k' \sigma | \ln \sigma|^{1/2} \right)^{\beta-1} .
$$
Then, for all $x$ and all $y \in D_x$, we have
\begin{equation}
    \label{borneu1}
    \ln f(x) - L^{(1)}(x) |y-x|\leq \ln f (y) \leq \ln f(x) + L^{(1)}(x) |y-x|.
\end{equation}
Note that for all $x$ in $ A_{\sigma}$, $L^{(1)}(x) \leq \B e^{k'}  \sigma^{-1} | \ln \sigma |^{-1/2}$. Starting from (\ref{borneu1}) and following the proof of Lemma~\ref{Lemma1} for the case $\beta=1$, it yields for all $x$ in $ A_{\sigma}$ and for all $H' >0$ and all $\sigma < \bar \sigma(\beta,H')$
\begin{equation}\label{KsigmaR2}
    (K_{\sigma} f) (x) = f(x) \left[1 + L^{(1)}(x) \sigma O_{\beta}(1)\right] + O_{\beta,H'}(\sigma^{H'}).
\end{equation}
For every $x \in E_{\sigma} \cap A_{\sigma}$, taking $H' = H_1+1$ in (\ref{KsigmaR2}) it yields
\begin{eqnarray*}
\frac{f_1(x)}{f(x)} & = & 2  - \frac{K_{\sigma}f(x)}{f(x)} \\
                    & = & 1 - L^{(1)}(x) \sigma O_{\beta}(1) - \sigma \frac{ \sigma^{H_1}}{f(x)} O_{\beta,H_1+1}(1).
\end{eqnarray*}
Next, $\left|\frac{ \sigma^{H_1}}{f(x)}\right| \leq 1 $ since $x \in E_{\sigma}$. Thus, for every $h\in (0,1)$, there exists $\bar \sigma(\beta,H_1,h)$ such that for all $\sigma \leq \bar{\sigma(\beta,H_1,h)}$, for all $f \in \H(\beta,\Pc)$ and every $x \in E_{\sigma} \cap A_{\sigma}$, $f_1(x) > (1 - h) f(x)$.

\noindent $\bullet$ The previous point is sufficient for $\beta \leq 4$ since $k=1$ in this case. We now also suppose that $\beta >4$ and thus that $k \geq 2$. Suppose that the integer $2 \leq u \leq  k$ is such that (\ref{fugeqf}) is true for the integer $u-1$. Let $h \in (0,1)$, there exists $\bar \sigma(\beta,H_1,h)$ such that for all $\sigma < \bar \sigma(\beta,H_1,h)$, for all $f \in \H(\beta,\Pc)$ and every $x \in E_{\sigma} \cap A_{\sigma}$, $f_{u-1}(x) > \left(1 - \frac h 2 \right) f(x)$. Note that since $2u \leq \beta$
we find that for all $x$ and all $y \in D_x$,
\begin{equation*}
    \ln f(x) +  \sum_{j=1}^{2u-1 } \frac{l_j(x)}{j!} (y-x)^j - (y-x)^{2u} L^{(u)}(x)\leq \ln f(y) \leq \ln f(x) + \sum_{j=1}^{2u -1 } \frac{l_j(x)}{j!} (y-x)^j + (y-x)^{2u} L^{(u)}(x)
    \label{borne1beta2}
\end{equation*}
with $L^{(u)}(x) =  \sum_{j=2u}^{r} \frac{|l_j(x)|}{j!} \left(k' \sigma | \ln \sigma|^{1/2} \right) ^{j-2u}  + L(x)  \left(k' \sigma | \ln \sigma|^{1/2} \right)^{\beta - 2u}$. Thus for all $x$ in $ A_{\sigma}$, $\left| L^{(u)}(x) \right| \leq \B e^{k'}  \sigma^{-2u} | \ln \sigma |^{-u}$. Following the proof of Lemma \ref{Lemma1} for the case $\beta = 2u$, it yields for all $x \in  A_{\sigma}$ and for all $H' >0$ and all $\sigma < \bar \sigma(\beta,H')$,
\begin{equation} \label{Ksigfup1}
    (K_{\sigma} f_{u-1}) (x) = f(x) \left[1 + R^{(u)}(x) O_{\beta}( \sigma^{2u})\right] +  O_{\beta,H'}(\sigma^{H'})
\end{equation}
with $ R^{(u)} = a_{2u+2} L^{(u)}(x) + \sum_{j=1}^{2u}   a_{j} |l_{j}(x)| ^{\frac {2u} j} $ and we have  $ \sup_{f \in \H(\beta,\Pc)} \sup_{x \in E_{\sigma} \cap A_{\sigma}}  \left| \sigma^{2u} R^{(u)}(x) \right| \leq \frac{c_{\beta}}{|\ln \sigma |^{u}} $. Then, using (\ref{Ksigfup1}) with $H' =H_1+1$, it yields for all $x  \in  A_{\sigma} \cap E_{\sigma} $
\begin{eqnarray*}
\frac{f_{u}(x)}{f(x)} & = & 1 - \frac{K_{\sigma} f_{u-1}(x)  - f_{u-1}}{f(x)} \\
& = & -\left[R^{(u)}(x) O_{\beta}(\sigma^{2u}) + \frac{\sigma^{H_1+1}}{f(x)} c_{\beta,H_1+1} \right]+ \frac {f_{u-1}(x)}{f(x)}.
\end{eqnarray*}
There exists $\bar{\sigma} (\beta,H_1,h) $ such that for all $\sigma \leq \bar{\sigma} (\beta,H_1,h) $, $\left|  R^{(u)}(x) O_{\beta}(\sigma^{2u}) + \frac{\sigma^{H_1+1}}{f(x)} c_{\beta,H_1+1} \right| \leq h /2 $ and the induction is complete. By choosing $u = k$ and $h = 1/2$, it finally gives that  $E_{\sigma} \cap A_{\sigma} \subset J_{\sigma,k}$ for $\sigma \leq \bar{\sigma} (\beta,H_1) $.

\paragraph{Proof of (\ref{W24}):}
We have
\begin{eqnarray*}
\int_\R g_k(x)  dx  &=& \int_{J_{\sigma,k}} f_k(x) dx + \int_{J_{\sigma,k}^c} \frac 1 2 f(x) dx\\
                    &=& 1 + \int_{J^c_{\sigma,k}} \left[\frac 1 2 f(x) - f_k(x)\right] dx
\end{eqnarray*}
since $\int_\R f_k(x) dx=1$ (see Lemma~\ref{Lemfk}). Moreover, $ f_k$  is a linear combination of $K_{\sigma}^t f$, $t = 0, \dots,k$, according to Lemma~\ref{Lemfk}. Thus it yields $ \int_\R g_k(x)  dx =  1+ O_{\beta,H_1}( \sigma^{2 \beta}) $ for $\sigma \leq \bar{\sigma} (\beta,H_1) $ thanks to (\ref{W23g}), (\ref{W23d}) and that $J^c_{\sigma,k}\subset A_\sigma^c \cup E_\sigma^c$ for $\sigma \leq \bar{\sigma} (\beta,H_1)$.

\paragraph{Proof of (\ref{W25}):}
Let $H >0$. According to Lemma \ref{Lemma1} and (\ref{W24}), for all $\sigma \leq \bar{\sigma} (\beta,H,H_1) $, for all $x \in E_{\sigma} \cap A_{\sigma} \subset J_{\sigma,k}$, we have
\begin{eqnarray*}
\left|K_{\sigma} h_k(x) -f(x)\right| & \leq & \left(\int_\R g_k(y) dy\right)^{-1} \,  \left| K_{\sigma} f_k (x) -f(x)\right| + \left|\left(\int_\R g_k(y) dy\right)^{-1} - 1\right| f(x)\\
                                     &      & + \left(\int_\R g_k(y) dy\right)^{-1} \int_{J^c_{\sigma,k}} \left\{\frac 1 2 f(u) - f_k(u)\right\} \psi_\sigma (x-u) du\\
                                     &\leq & c_\beta f(x) R_f(x) \sigma^\beta + c_\beta \sigma^H + \left(\int_\R g_k(y) dy\right)^{-1} (\diamondsuit)
\end{eqnarray*}
where $(\diamondsuit)=\int_{J^c_{\sigma,k}} \{\frac 1 2 f(u) - f_k(u)\} \psi_\sigma (x-u) du$. Let $D_x=\{u\in\R; |x-u|\leq k' \sigma |\ln \sigma|^{\frac 1 2}\}$ such that $k' \sigma |\ln \sigma|^{\frac 1 2}\leq \gamma$. According to the third result of Lemma~\ref{Lemfk} and Lemma~\ref{Lemma10},
$$
  0 \leq  \int_{J^c_{\sigma,k}\cap D_x^c} \left\{\frac 1 2 f(u) - f_k(u)\right\} \psi_\sigma (x-u) du \leq 2^{k+2}\frac{M}{\sqrt{\pi}} \int_{D_x^c} \psi_\sigma(x-u) du \leq c_\beta \sigma^H.
$$
Next, if $x\in A_{\sigma}\cap E_{\sigma}$, there exists $t>1$ such that for all $u\in D_x$, $u\in A_{\sigma,t}\cap E_{\sigma,t}$. This result can be proved by adapting some parts of the proof of (\ref{W23g}) and (\ref{W23d}). Moreover, by changing $\B$ into $t \B$, it can be shown that there exists $\bar\sigma(\beta)$ such that for all $\sigma<\bar\sigma(\beta)$, $A_{\sigma,t}\cap E_{\sigma,t}\subset J_{\sigma,k}$. Thus for $\sigma$ small enough, $\int_{J^c_{\sigma,k}\cap D_x^c} \{\frac 1 2 f(u) - f_k(u)\} \psi_\sigma (x-u) du =0.$
Finally, $(\diamondsuit)\leq c_\beta \sigma^H$ and $\left|K_{\sigma} h_k(x) -f(x)\right|\leq c_\beta f(x) R_f(x) \sigma^\beta + c'_\beta \sigma^H$.
\end{proof}

\begin{prop}\label{PropKLFKsigmahk}~\\
Let $\beta>0$ and $k\in\mathbb{N}$ such that $\beta\in(2k,2k+2]$. There exists a positive constant $\bar \sigma (\beta)$ such that for all $ f \in \H(\beta,\Pc)$ and all $\sigma<  \bar \sigma (\beta)$,
$$
    \mbox{KL}(f,K_\sigma h_k) = \int_{\mathbb{R}} f(x) \ln\left(\frac{f(x)}{K_\sigma h_k(x)}\right) dx = O_{\beta}(\sigma^{2\beta})
$$
where $h_k$ is defined by (\ref{W20b}) and where $\bar \sigma (\beta)$ can be chosen as a continuous function of $\beta$.
\end{prop}
\begin{proof}
Preliminary, we remark that if $p$ and $q$ are two densities and $S$ is a set, then
$$
    \int_S p \ln\left(\frac p q\right) \leq \int_S p\ \frac{p-q}{q} = \int_S  \frac{(p-q)^2}{q} + \frac{q(p-q)}{q} = \int_S  \frac{(p-q)^2}{q} + \int_{S^c}(q-p)
$$
since $\int_S p = 1 - \int_{S^c} p$, $\int_S q = 1 - \int_{S^c} q$ and $\int_S (p-q) = \int_{S^c} (q-p)$. We use this inequality with the densities $f$ and $K_\sigma h_k$, and the sets $A_\sigma$ and $ E_\sigma$, where $E_\sigma$ is defined with $H_1=4 \beta+1$, to obtain the following control of $\KL(f,K_\sigma h_k)$:
\begin{eqnarray}
\int_{\mathbb{R}} f(x) \ln\left(\frac{f(x)}{K_\sigma h_k(x)}\right) dx
&=& \int_{A_\sigma \cap E_\sigma} f(x) \ln\left(\frac{f(x)}{K_\sigma h_k(x)}\right) dx + \int_{A_\sigma^c \cup E_\sigma^c} f(x) \ln\left(\frac{f(x)}{K_\sigma h_k(x)}\right) dx\notag\\
&\leq& \int_{A_\sigma \cap E_\sigma} \frac{[f(x)-K_\sigma h_k(x) ]^2}{K_\sigma h_k(x)}dx\label{PropKLINT1}\\
&+& \int_{A_\sigma^c \cup E_\sigma^c} [K_\sigma h_k (x) -f(x)] dx\label{PropKLINT2}\\
&+& \int_{A_\sigma^c \cup E_\sigma^c} f(x) \ln\left(\frac{f(x)}{K_\sigma h_k(x)}\right) dx.\label{PropKLINT3}
\end{eqnarray}

\noindent $\bullet$ \textbf{Control of (\ref{PropKLINT1})}:\\
Let $H>0$. According to Lemma~\ref{Lemme2} with $H_1=4\beta+1$, there exists $\bar\sigma (\beta,H)>0$ such that for all $x \in A_\sigma \cap E_\sigma$ and for all $\sigma < \bar \sigma (\beta,H)$,
$\left[K_\sigma h_k(x) - f(x)\right]^2 \leq \left[\Lambda_{\beta,H}f(x) R_f(x)\sigma^\beta + \Omega_{\beta,H} \sigma^H\right]^2$
where $\Lambda_{\beta,H} $ and $\Omega_{\beta,H}$ are two constants.
Moreover, according to Lemma~\ref{BoiteOutilBornes0}, there exists $\bar{\sigma}(\beta)>0$ such that for all $\sigma<\bar{\sigma}(\beta)$,
$$
    K_\sigma h_k(x) \geq \frac{D}{1+A_{\beta} \sigma^{2\beta}} f(x)
$$
with $D=\frac{\cinf\sqrt{\pi}}{6 M}$. Thus for all $\sigma < \bar \sigma (\beta,H) \wedge \bar\sigma (\beta)$,
\begin{eqnarray*}
\frac{\left[f(x) - K_\sigma h_k(x)\right]^2}{K_\sigma h_k(x)}
&\leq & \frac{\Lambda_{\beta,H}^2}{D}  (1+A_{\beta}\sigma^{2\beta})\sigma^{2\beta} R_f(x)^2 f(x) + \frac{\Omega_{\beta,H}^2}{D}  (1+A_{\beta}\sigma^{2\beta})\sigma^{2H}  \frac{1}{f(x)}\\
&+    & \frac{2 \Lambda_{\beta,H} \Omega_{\beta,H}}{D} (1+A_{\beta}\sigma^{2\beta})\sigma^{\beta + H} R_f(x).
\end{eqnarray*}
Then,
\begin{eqnarray}
\int_{A_\sigma \cap E_\sigma} \frac{[f(x)-K_\sigma h_k(x) ]^2}{K_\sigma h_k(x)}dx
&\leq & \frac{\Lambda_{\beta,H}^2}{D} (1+A_{\beta}\sigma^{2\beta})\sigma^{2\beta} \int_{A_\sigma \cap E_\sigma} R_f(x)^2 f(x)dx\notag\\
&+    & \frac{\Omega_{\beta,H}^2}{D} (1+A_{\beta}\sigma^{2\beta})\sigma^{2H-2(4\beta+1)} \int_{A_\sigma \cap E_\sigma}   f(x)dx  \notag\\
&+    & \frac{2 \Lambda_{\beta,H} \Omega_{\beta,H}}{D} (1+A_{\beta}\sigma^{2\beta})\sigma^{\beta + H-4\beta-1} \int_{A_\sigma \cap E_\sigma} R_f(x) f(x)dx . \hskip 1cm \label{Prop1ExpAux}
\end{eqnarray}
Thus the two integrals
$\int_{A_\sigma \cap E_\sigma} R_f(x)^2 f(x)dx$ and $\int_{A_\sigma \cap E_\sigma} R_f(x) f(x)dx$ have to be controlled.

The first integral can be decomposed into
\begin{eqnarray*}
\int_{A_\sigma \cap E_\sigma} R_f(x)^2 f(x)dx
&=& \int_{A_\sigma \cap E_\sigma} \left[a_{r+1} L(x) + \sum_{j=1}^r a_j |l_{j}(x)|^{\frac \beta j}\right]^2f(x) dx\\
&=& a_{r+1}^2 \int_{A_\sigma \cap E_\sigma} L(x)^2 f(x) dx + \sum_{j=1}^r a_j^2 \int_{A_\sigma \cap E_\sigma} |l_{j}(x)|^{\frac{2\beta}{j}} f(x) dx\\
& & + 2\sum_{j=1}^r a_{r+1}a_j \int_{A_\sigma \cap E_\sigma} |l_{j}(x)|^{\frac \beta j} L(x) f(x) dx\\
& & + \sum_{\underset{j\neq j'}{j,j'=1}}^r a_{j'}a_j \int_{A_\sigma \cap E_\sigma} |l_{j}(x)|^{\frac \beta j} |l_{j'}(x)|^{\frac{\beta}{j'}} f(x) dx.
\end{eqnarray*}
Using the Hölder inequality and Condition (\ref{Cond2a}), for all  $ j=1,\dots,r$,
\begin{equation}
\int_{A_\sigma \cap E_\sigma} |l_{j}(x)|^{\frac{2\beta}{j}} f(x) dx
\leq  \left[\int_{\mathbb{R}} |l_{j}(x)|^{\frac{2\beta+\varepsilon}{j}} f(x) dx\right]^{\frac{2\beta}{2\beta+\varepsilon}} \left[\int_{\mathbb{R}} f(x) dx\right]^{\frac{\varepsilon}{2\beta+\varepsilon}}\leq  C^{\frac{2 \beta}{2\beta+\varepsilon}}
\label{intlj}
\end{equation}
and
$\int_{A_\sigma \cap E_\sigma} L(x)^2 f(x) dx \leq \left[\int_{\mathbb{R}} L(x)^{2 +\frac{\varepsilon}{\beta}}(x) f(x) dx\right]^{\frac{2\beta}{2\beta+\varepsilon}} \left[\int_{\mathbb{R}} f(x) dx\right]^{\frac{\varepsilon}{2\beta+\varepsilon}}\leq C^{\frac{2 \beta}{2\beta+\varepsilon}}.$
Next, using the Cauchy-Schwarz inequality and (\ref{intlj}),
for all $j,j'\in\{1,\dots,r\}, j\neq j'$,
$$
\int_{A_\sigma \cap E_\sigma} |l_{j}(x)|^{\frac \beta j} |l_{j'}(x)|^{\frac{\beta}{j'}}f(x) dx \leq  \left[\int_{\mathbb{R}} |l_{j}(x)|^{\frac{2\beta}{j}} f(x) dx\right]^{\frac 1 2} \left[\int_{\mathbb{R}} |l_{j'}(x)|^{\frac{2\beta}{j'}} f(x)dx\right]^{\frac 1 2}\leq  C^{\frac{2 \beta}{2\beta+\varepsilon}}
$$
and for all $j\in\{1,\dots,r\}$,
$$
\int_{A_\sigma \cap E_\sigma} |l_{j}(x)|^{\frac \beta j} L(x) f(x) dx \leq  \left[\int_{\mathbb{R}} |l_{j}(x)|^{\frac{2\beta}{j}} f(x) dx\right]^{\frac 1 2} \left[\int_{\mathbb{R}} L(x)^2 f(x) dx\right]^{\frac 1 2} \leq  C^{\frac{2 \beta}{2\beta+\varepsilon}}.
$$
Finally, $\int_{A_\sigma \cap E_\sigma} R_f(x)^2 f(x)dx \leq \left(\underset{j=1}{\stackrel{r+1}{\sum}}a_j\right)^2 C^{\frac{2 \beta}{2\beta+\varepsilon}}$.

For the second integral,
\begin{eqnarray*}
\int_{A_\sigma \cap E_\sigma} \hskip-0.5cm R_f(x) f(x)dx &=& \int_{A_\sigma \cap E_\sigma} \left[a_{r+1} L(x) + \sum_{j=1}^r a_j |l_{j}(x)|^{\frac \beta j}\right] f(x) dx\\
&\leq & a_{r+1} \sqrt{\int_{\mathbb{R}}L(x)^2 f(x) dx}\sqrt{\int_{\mathbb{R}} f(x) dx} + \sum_{j=1}^r a_j \sqrt{\int_{\mathbb{R}}|l_{j}(x)|^{\frac{2\beta}{j}} f(x) dx} \sqrt{\int_{\mathbb{R}} f(x) dx}\\
&\leq & \underset{j=1}{\stackrel{r+1}{\sum}}a_j C^{\frac{\beta}{2\beta+\varepsilon}}.
\end{eqnarray*}
Finally, (\ref{Prop1ExpAux}) becomes
\begin{eqnarray*}
\int_{A_\sigma \cap E_\sigma} \frac{[f(x)-K_\sigma h_k(x) ]^2}{K_\sigma h_k(x)}dx
&\leq & \frac{\Lambda_{\beta,H}^2}{D} (1+A_{\beta}\sigma^{2\beta})\sigma^{2\beta} \left(\underset{j=1}{\stackrel{r+1}{\sum}}a_j\right)^2 C^{\frac{2 \beta}{2\beta+\varepsilon}}\\
   &+& \frac{\Omega_{\beta,H}^2}{D} (1+A_{\beta}\sigma^{2\beta})\sigma^{2H-8\beta -2} \\
   &+& \frac{2 \Lambda_{\beta,H} \Omega_{\beta,H}}{D} (1+A_{\beta}\sigma^{2\beta})\sigma^{H-3\beta-1} \left(\underset{j=1}{\stackrel{r+1}{\sum}}a_j\right) C^{\frac{\beta}{2\beta+\varepsilon}}.
\end{eqnarray*}
By taking $H =  5\beta+1$, it gives that there exists $\bar \sigma (\beta)>0$ such that for all  $\sigma < \bar \sigma (\beta)$,
$$
    \int_{A_\sigma \cap E_\sigma} \frac{[f(x)-K_\sigma h_k(x) ]^2}{K_\sigma h_k(x)}dx = O_{\beta}(\sigma^{2\beta}).
$$

\noindent $\bullet$ \textbf{Control of (\ref{PropKLINT2})}:\\
According to Lemma~\ref{Lemfk},
\begin{eqnarray*}
h_k(x)
&=& \left(\int_{\mathbb{R}} g_k(x) dx\right)^{-1}\ \left\{f_k(x)\mathds{1}_{J_{\sigma,k}}(x) + \frac 1 2 f(x) \mathds{1}_{J_{\sigma,k}^c}(x)\right\}\\
&\leq& \left\{2 \sum_{i=0}^k (-1)^i \left(\substack{k+1 \\ i+1}\right) K_\sigma^i f(x)\right\} \mathds{1}_{J_{\sigma,k}}(x) + f(x) \mathds{1}_{J_{\sigma,k}^c}(x)
\end{eqnarray*}
thus
$$
    K_\sigma h_k(x) \leq 2 \sum_{j=1}^{k+1}  \left(\substack{k+1\\j}\right) K_\sigma^{j} f(x)  + K_\sigma f(x).
$$
According to (\ref{W23g}) and (\ref{W23d}) in Lemma~\ref{Lemme2} with $H_1=4\beta+1$, there exists $\bar \sigma(\beta)>0$ such that for all $\sigma < \bar \sigma(\beta)$,
\begin{eqnarray*}
\int_{A_\sigma^c \cup E_\sigma^c} \hskip-0.5cm  [K_\sigma h_k(x) -f(x)] dx
&\leq & \int_{A_\sigma^c \cup E_\sigma^c} K_\sigma h_k(x) dx + \int_{A_\sigma^c \cup E_\sigma^c} f(x)dx\\
&\leq & 2 \sum_{j=1}^{k+1}  \left(\substack{k+1\\j}\right) \int_{A_\sigma^c \cup E_\sigma^c}K_\sigma^{j} f(x) dx + \int_{A_\sigma^c \cup E_\sigma^c} K_\sigma f(x) dx + \int_{A_\sigma^c \cup E_\sigma^c} f(x)dx\\
&\leq & 2 \sum_{j=2}^{k+1}  \left(\substack{k+1\\j}\right) \int_{A_\sigma^c} \hskip-0.2cm K_\sigma^{j} f(x) dx + [2(k+1)+1]\int_{A_\sigma^c}  \hskip-0.2cm K_\sigma f(x) dx + \int_{A_\sigma^c} \hskip-0.2cm K_\sigma^0 f(x) dx\\
& & \;+ \; 2 \sum_{j=2}^{k+1} \left(\substack{k+1\\j}\right) \int_{E_\sigma^c} \hskip-0.2cm K_\sigma^{j} f(x) dx + [2(k+1)+1]\int_{E_\sigma^c} \hskip-0.2cm  K_\sigma f(x) dx+\int_{E_\sigma^c}  \hskip-0.2cm  K_\sigma^0 f(x) dx\\
&\leq & 2 \left[2 \sum_{j=2}^{k+1} \left(\substack{k+1\\j}\right) + 2(k+2)\right] c_\beta \sigma^{2\beta} .
\end{eqnarray*}

\noindent $\bullet$ \textbf{Control of (\ref{PropKLINT3})}:\\
According to Lemma~\ref{BoiteOutilBornes0}, for all $\sigma < \bar \sigma(\beta)$, $K_\sigma h_k(x)\geq \frac{D}{1+ A_{\beta} \sigma^{2\beta}} f(x)$ then
\begin{eqnarray*}
\int_{A_\sigma^c \cup E_\sigma^c} f(x) \ln\left(\frac{f(x)}{K_\sigma h_k(x)}\right) dx
&\leq & \ln\left(\frac{1+ A_{\beta} \sigma^{2\beta}}{D}\right) \int_{A_\sigma^c \cup E_\sigma^c} f(x) dx\\
&\leq & \ln\left(\frac{1+ A_{\beta} \sigma^{2\beta}}{D}\right) \left\{\int_{A_\sigma^c} K_\sigma^0 f(x) dx + \int_{E_\sigma^c} K_\sigma^0 f(x) dx\right\}\\
&\leq & \ln\left(\frac{1+ A_{\beta} \sigma^{2\beta}}{D}\right) 2 c_{\beta} \sigma^{2\beta}.
\end{eqnarray*}

In conclusion, there exists $\bar \sigma(\beta)>0$ such that for all $ \sigma <\bar \sigma( \beta)$, $\KL(f,K_\sigma h_k)=O_{\beta}(\sigma^{2\beta})$.

\end{proof}


\subsection{Proof of Theorem~\ref{ThApproximation}}\label{subsect:ProofThApproximation}
\begin{proof}
For the definition of  $E_\sigma$, we choose $H_1=4(\beta+1)$. Let $\tilde{h}_k$ be the restriction of $h_k$ on an interval $[-\mu_\sigma,\mu_\sigma]$, normalized in order to have a density function: 
$$
    \tilde h_k:x\in\R\mapsto \left(\int_{[-\mu_{\sigma},\mu_{\sigma}]} h_k(y) dy\right)^{-1} h_k(x) \mathds{1}_{[-\mu_\sigma,\mu_\sigma]}(x)
$$ 
where $\mu_\sigma$ depends on $\sigma$ and will be chosen further such that
\begin{equation}
\label{aleqsig}
\mu_{\sigma} \geq \sigma .
\end{equation}
Let $\varepsilon \in(0,\pi^{-1/2})$. According to Proposition~\ref{PropMesureDiscrete} in Appendix~\ref{AppMesDiscrete}, there exists a discrete distribution $\tilde F$ on $[-\mu_{\sigma},\mu_{\sigma}]$ with at most $54 \mu_{\sigma} \sigma^{-1} e^2 \left[-\ln\left(\sqrt{\pi}\varepsilon\right) \vee 1\right]$  support points such that
\begin{equation}
\label{applimelangediscret}
\| \tilde h_k \ast\psi_{\sigma} -  \tilde{F}  \ast \psi_\sigma \|_{\infty} \leq \frac{2 \varepsilon}{\sigma}.
\end{equation}
Denoting $\tilde \mel(x) dx= \left( \int_{[-\mu_\sigma,\mu_\sigma]}h_k(x) dx \right) \tilde F \ast \psi_\sigma (x)$, it gives for all $x \in \R$,
\begin{eqnarray*}
 \left|  K_\sigma h_k (x) - \tilde{\mel} (x) \right|
 &= & \left( \int_{[-\mu_\sigma,\mu_\sigma]}h_k(x)  dx \right) \left|  \frac{h_k}{   \int_{[-\mu_\sigma,\mu_\sigma]}h_k(y) dy} \ast \psi_{\sigma}(x) -  \tilde{F}  \ast \psi_\sigma (x) \right|  \\	
 & \leq & \left|  \frac{h_k  \mathds{1}_{[-\mu_\sigma,\mu_\sigma]} }{   \int_{[-\mu_\sigma,\mu_\sigma]}h_k(y) dy } \ast \psi_{\sigma}(x) -  \tilde{F}  \ast \psi_\sigma (x) \right| +  \left(h_k  \mathds{1}_{[-\mu_\sigma,\mu_\sigma]^c} \right) \ast \psi_{\sigma}(x).
\end{eqnarray*}
By applying Lemma~\ref{BoiteOutilBornes} with $p=\frac 1 2$, it gives that for all $\sigma \leq  1 - 2^{-1/k}$ and for all $x \in \R$,
$$
    h_k(x) \leq 4 M \left( \frac{4}{\sqrt 3}\right)^k \psi\left(\frac x 2\right)
$$
and thus $\left(h_k  \mathds{1}_{[-\mu_\sigma,\mu_\sigma]^c} \right) \ast \psi_{\sigma}(x) \leq   4 M \left( \frac{4}{\sqrt 3}\right)^k  \psi(\frac {\mu_\sigma} 2)$. Now, we choose $\mu_\sigma:= 2 \sqrt{    \ln\left(\frac{4 M}{\sqrt{\pi}}\left( \frac{4}{\sqrt 3}\right)^k \frac \sigma  \varepsilon\right)}$  in order to obtain that $\left\| \left(h_k  \mathds{1}_{[-\mu_\sigma,\mu_\sigma]^c} \right) \ast \psi_{\sigma} \right\|_{\infty} \leq \frac \varepsilon \sigma$. This last inequality together with (\ref{applimelangediscret}) yields
\begin{equation*} \label{controlepseudomelange}
 \left\| K_\sigma h_k  -\tilde{\mel} \right\|_{\infty}  \leq  \frac {3 \varepsilon } \sigma  .
\end{equation*}
We also define the function $t:= \tilde{\mel} + \sigma^{6 \beta + 5} \psi_{\sigma}$ and the finite Gaussian mixture with density
$$
   \mel(x) :=   \frac{t(x)}{\int_{\R} t(y) dy}  = \frac{\tilde{\mel}(x) + \sigma^{6 \beta + 5} \psi_{\sigma}(x)}{  \int_{[-\mu_\sigma,\mu_\sigma]} h_k(y) dy + \sigma^{6 \beta + 5}}.
$$
Then we want to upper bound
\begin{eqnarray*}
\mbox{KL}(f,\mel) &=& \int_{\mathbb{R}} f(x) \ln\left(\frac{f(x)}{\mel(x)}\right) dx\\
                  &=& \int_{\mathbb{R}} f(x) \ln\left(\frac{f(x)}{K_\sigma h_k(x)}\right) dx
                 + \int_{\R} f(x) \ln\left(\frac{K_\sigma h_k(x)}{t(x)}\right) dx
                 + \int_{\R} f(x) \ln\left(\frac{t(x)}{\mel(x)}\right) dx \\
                  &=&  \int_{\mathbb{R}} f(x) \ln\left(\frac{f(x)}{K_\sigma h_k(x)}\right) dx\\
               &   &  + \int_{E_\sigma^c} f(x) \ln\left(\frac{K_\sigma h_k(x)}{t(x)}\right) dx
                      + \int_{E_\sigma} f(x) \ln\left(\frac{K_\sigma h_k(x)}{t(x)}\right) dx
                      + \int_{\R} f(x) \ln\left(\frac{t(x)}{\mel(x)}\right) dx \\
               &=& \fbox{I1} + \fbox{I2} + \fbox{I3} + \fbox{I4}.
\end{eqnarray*}

\noindent $\bullet$ \textbf{Control of $\fbox{I1}$}: According to Proposition~\ref{PropKLFKsigmahk}, for all $\sigma < \bar\sigma(\beta)$,
$$
    \int_{\mathbb{R}} f(x) \ln\left(\frac{f(x)}{K_\sigma h_k(x)}\right) dx =O_{\beta}(\sigma^{2\beta}).
$$

\noindent $\bullet$ \textbf{Control of $\fbox{I2}$}: According to Lemma~\ref{BoiteOutilBornes}, $K_\sigma h_k(x)\leq 4 M \left(\frac{4}{\sqrt 3}\right)^k $ for $\sigma$ small enough and since $s(x)\geq \sigma^{6\beta+5}\psi_\sigma(x)$,
\begin{eqnarray*}
\fbox{I2}
& \leq &  \int_{ E_{\sigma}^c} f(x) \ln \left(\frac {4 M \left(\frac{4}{\sqrt 3}\right)^k }{\sigma^{6\beta+5}\psi_\sigma(x)}\right)dx \\
& \leq &  \left( \int_{ E_{\sigma}^c} f(x) dx\right)  \left[ (6\beta +4) | \ln \sigma | + \ln\left(4 M \left(\frac{4}{\sqrt 3}\right)^k\right) \right] + \int_{ E_{\sigma}^c} f(x) \frac {x^2}{\sigma^2} dx .
 \end{eqnarray*}
For the second integral,
\begin{eqnarray*}
\int_{E_\sigma^c} \frac{x^2}{\sigma^2} f(x)dx &\leq & \sigma^{\frac{H_1}{2}-2} \int_{E_\sigma^c} x^2 \sqrt{f(x)}dx\\
&\leq & \sigma^{2\beta}   \int_{\R} x^2 \sqrt{M \psi(x)}dx \\
&\leq & \sigma^{2\beta} \sqrt{2M}\pi^{-\frac 1 4} 2 \nu_{1,2} = 4\pi  \sqrt{M} \sigma^{2\beta}.
\end{eqnarray*}
Similarly, $ \int_{E_\sigma^c} f(x)dx \leq  \sigma^{2\beta+2} \sqrt{2 M}$ and finally
$$
\fbox{I2}\leq \left\{\ln \left(4 M \left(\frac{4}{\sqrt 3}\right)^k\right) + (6\beta+4)|\ln\sigma| \right\} \sqrt{2M} \sigma^{2\beta+2} +   4\pi  \sqrt{M} \sigma^{2\beta}.
$$
Thus $\fbox{I2}=O_\beta(\sigma^{2\beta})$.\\

\noindent $\bullet$ \textbf{Control of $\fbox{I3}$}: On the one hand,
\begin{eqnarray*}
\left| K_{\sigma} h_k(x) - t(x) \right| & \leq  & \left| K_{\sigma} h_k(x) - \tilde \mel (x) \right| + \left| \tilde \mel (x) - t(x) \right| \\
                                        & \leq & 3 \varepsilon \sigma^{-1} + \sigma^{6 \beta + 5} \psi_{\sigma}(x)\\
                                        & \leq & 3 \varepsilon \sigma^{-1} + \sigma^{6 \beta + 4} \pi^{-1/2} .
\end{eqnarray*}
On the other hand, according to Lemma \ref{BoiteOutilBornes0}, for all $x \in \R$ and for all $\sigma< \bar \sigma(\beta)$,
$K_\sigma h_k(x)  \geq \frac{\cinf\sqrt{\pi}}{6 M(1+A_\beta\sigma^{2\beta})}  f(x)$.
Since $x \in E_{\sigma}$ then $K_\sigma h_k(x) \geq \frac{\cinf\sqrt{\pi}}{6 M(1+A_\beta\sigma^{2\beta})}\sigma^{4(\beta+1)}$. Thus, $t(x)\geq  \tilde{\mel}(x)   \geq   K_{\sigma} h_k(x)  - 3 \varepsilon \sigma^{-1}    \geq   \frac{\sigma^{4(\beta+1)}\sqrt{\pi}}{6M (1+ A_\beta \sigma^{2 \beta})}  - 3 \varepsilon \sigma^{-1} $.
Finally,
\begin{eqnarray*}
\fbox{I3}
& \leq& \int_{E_{\sigma}} f(x)  \frac{ K_{\sigma} h_k(x)  -t(x) } { t(x) }  \ dx \\
& \leq&  \frac{ 3 \varepsilon \sigma^{-1} + \sigma^{6 \beta + 4 } \pi^{-1/2} } { \frac{\sigma^{4(\beta+1)}}{2 (1+ A_\beta \sigma^{2 \beta})}  - 3 \varepsilon \sigma^{-1} } \int_{E_{\sigma}} f(x)   \ dx \\
& \leq&  \frac{ 3 \varepsilon \sigma^{-1} + \sigma^{6 \beta + 4 } \pi^{-1/2} } { \frac{\sigma^{4(\beta+1)}}{2 (1+ A_\beta \sigma^{2 \beta})}  - 3 \varepsilon \sigma^{-1} }.
\end{eqnarray*}
Let $\delta'  := 1 + \frac{ \beta}{2(\beta+1)}$ and we set $ \varepsilon :=  \sigma^{\delta' 4(\beta+1) +1}$. It yields
\begin{eqnarray*}
\fbox{I3}  & \leq&  \frac{ (\pi^{-1/2} + 3 ) \sigma^{6 \beta +4}  } { \frac{\sigma^{4(\beta+1)}}{2 (1+ A_\beta \sigma^{2 \beta})}  - 3 \sigma^{6 \beta+4} } = O_{\beta} \left( \sigma^{2 \beta }\right).
\end{eqnarray*}

\noindent $\bullet$ \textbf{Control of $\fbox{I4}$}: Note that $\frac{t(x)}{\mel(x)} = \int_{[-\mu_\sigma,\mu_\sigma]} h_k(y) dy  + \sigma^{6 \beta + 5} \leq 1 + \sigma^{6 \beta +5} $ and thus
\begin{eqnarray*}
\fbox{I4}  & \leq & \int_{\R}   f(x) \ln \left( 1 + \sigma^{6 \beta + 5} \right) \ dx \\
& \leq & \sigma^{6 \beta +5}  \leq  \sigma^{2 \beta} .
\end{eqnarray*}

Finally, we obtain that $\KL(f,\mel)=O_\beta(\sigma^{2\beta})$. Moreover, according to the choice of $\varepsilon$, we have that
\begin{eqnarray*}
\mu_\sigma &=& 2 \sqrt{\ln\left(\frac{4 M}{\sqrt{\pi}} \left( \frac{4}{\sqrt{3}} \right)^k \frac{\sigma}{\varepsilon}\right)}\\
         &=& 2 \sqrt{\ln\left(\frac{4 M}{\sqrt{\pi}} \left( \frac{4}{\sqrt{3}} \right)^k \sigma^{-(6\beta+4)} \right)}\\
         &=&  \tilde G_\beta |\ln \sigma|^{\frac 1 2}
\end{eqnarray*}
where
\begin{equation}
\label{Gtildebeta}
\tilde G_\beta =  2\sqrt{\ln\left(\frac{4 M}{\sqrt{\pi}}\right) + k \ln\left(\frac{4}{\sqrt{3}} \right) + (6\beta+4)}.
\end{equation}
Thus there exists $\bar \sigma(\beta) $ continuous in $\beta$ such that (\ref{aleqsig}) is fulfilled for $\sigma < \bar \sigma(\beta) $. Furthermore, the mixture $\mel$ has $k_\sigma$ components such that
\begin{eqnarray}
k_\sigma &\leq& 54 \mu_\sigma \sigma^{-1} e^2 \left[1\vee\ln\left(\frac{1}{\sqrt{\pi}\varepsilon} \right)\right] +1 \notag \\
         &\leq & \tilde G_\beta |\ln \sigma|^{\frac 1 2} 54 \sigma^{-1} e^2 \left[1\vee\ln\left(\frac{1}{\sqrt{\pi}\sigma^{6\beta+5}} \right)\right] +1  \notag  \\
         & = &  G_\beta \sigma^{-1}  |\ln \sigma|^{\frac 3 2} . \label{Gbeta}
\end{eqnarray}
\end{proof}


\section{Proof of the lower bound}\label{ProofLowB}
\subsection{Proof of Proposition~\ref{prop:exist}}\label{subs:propexist}
Note that for every $j$, $\varphi_{j}$ is supported by
$$
    J_j:=\left[-\frac{\xM}{2} + \frac{\xM}{D}(j-1) + \frac{\xM}{4D}, -\frac{\xM}{2} + \frac{\xM}{D}j - \frac{\xM}{4D}\right]
    \subsetneq I_j=\left[-\frac{\xM}{2} + \frac{\xM}{D}(j-1), -\frac{\xM}{2} + \frac{\xM}{D}j\right]
$$
and thus the supports of the $\varphi_{j}, 1\leq j\leq D$ are disjoint.  We also note that for all $x\in[-\frac{\xM}{2},\frac{\xM}{2}]^c,\ f_\theta(x)=\omega(x)$ and for all $x\in[-\frac{\xM}{2},\frac{\xM}{2}]$, there exists an unique $ j\in\{1,\dots,D\}$ such that $f_\theta(x)=2\cinf + (2\theta_j -1) \varphi_j(x)$ where $ \varphi_j(x)=0$ if $x\in I_j \backslash J_j$. The proof of Proposition~\ref{prop:exist} is decomposed into two lemmas.

\begin{lem}{Density function and monotonicity conditions.}~\label{propMonotonicity}\\
For all $D \in \N^*$ and all $\theta\in\{0,1\}^D$, the function
$f_{\theta}$ defined by (\ref{deffvtheta}) is a positive density
function such that for all $x\in [-\frac{\xM}{2},\frac{\xM}{2}],\
f_{\theta}(x)\in\left[\cinf,3\cinf\right]$. This function fulfills also the following monotonicity conditions:
\begin{enumerate}
    \item $\forall x\in[-\xM,\xM]$, $f_{\theta}(x)\geq \cinf$ and $\forall x\in[-\xM,\xM]^c$, $f_{\theta}(x)\leq \cinf$.
    \item $f_{\theta}$ is nondecreasing on $(-\infty,-\xM)$ and nonincreasing on $(\xM,\infty)$.
    \item $\forall x\in\R,\ f_\theta(x)\leq M\psi(x)$ with $M=\tilde M \vee 3\sqrt{\pi} \cinf \exp(\xM^2/4)$.
\end{enumerate}
\end{lem}

\begin{proof}
For all $x \in [-\frac{\xM}{2},\frac{\xM}{2}]^c,\ f_{\theta}(x) = \omega(x)>0$ since $\omega$ is positive.
Moreover, for all $x \in [-\frac{\xM}{2},\frac{\xM}{2}], \exists ! j\in\{1,\dots,D\}$ such that $x\in I_j$. Then,
$$
    f_{\theta}(x) = \omega(x) + (2\theta_j -1) \varphi_{j}(x)= 2\cinf + (2\theta_j -1) \varphi_{j}(x).
$$
Thus
\begin{eqnarray*}
\left|f_{\theta}(x) - 2\cinf\right| &=& |(2 \theta_j -1)|\ |\varphi_{j}(x)|\\
                                    &=& \left|\frac{\cinf D^{-\beta}}{A}\ \varphi\left(\frac{D}{\xM} \left(x+\frac{\xM}{2}\right) - (j-1)\right)\right|\\
                                    &\leq & \cinf D^{-\beta}\\
                                    &\leq & \cinf
\end{eqnarray*}
since $D^{-\beta}\leq 1$. Thus for all $x\in [-\frac{\xM}{2},\frac{\xM}{2}], f_{\theta}(x) \in \left[\cinf,3\cinf\right]$.
Finally, $f_\theta$ is a positive function on $\R$.
Moreover,
\begin{eqnarray*}
\int_{\R} f_{\theta}(x) dx &=& \int_{\R} \omega(x) dx + \underset{j=1}{\stackrel{D}{\sum}} (2\theta_j -1)\int_{I_j}\varphi_{j}(x) dx\\
                           &=& \int_{\R} \omega(x) dx + \underset{j=1}{\stackrel{D}{\sum}} (2\theta_j -1) \frac{\cinf D^{-\beta}}{A} \frac {\xM}{D} \int_{\R} \varphi(y) dy\\
                           &=& 1
\end{eqnarray*}
because $\int_{\R} \omega(x) dx=1$ and $\int_{\R} \varphi(y) dy=0$. Thus, $f_\theta$ is a density function.

On $(-\infty, -\xM)$, since $f_{\theta}(x)= \omega(x)$ and $\omega$ is a nondecreasing function
on $(-\infty,-\xM)$, the function $f_{\theta}$ is a nondecreasing
function on $(-\infty,-\xM)$. Moreover,
$$
    \forall x < -\xM, f_{\theta}(x) \leq f_{\theta}(-\xM) = \omega(-\xM)=\cinf.
$$
In the same way, the function $f_{\theta}$ is a nonincreasing
function on $(\xM,\infty)$ and
$$
    \forall x > \xM, f_{\theta}(x) \leq f_{\theta}(\xM) = \omega(\xM)=\cinf.
$$
For all $x\in [-\xM,\xM]$,
    \begin{itemize}
    \item if $x\in [-\xM,-\frac{\xM}{2})$, $f_{\theta}(x)= \omega(x) \geq \omega(-\xM) = \cinf$ because $\omega$ non-decreases and $\omega(-\xM)= \cinf$;

    \item if $x\in (\frac{\xM}{2},\xM]$, $f_{\theta}(x)= \omega(x) \geq \omega(\xM) = \cinf$ because $\omega$ non-increases and $\omega(\xM)= \cinf$;

    \item if $x\in[-\frac{\xM}{2},\frac{\xM}{2}]$,
    $f_{\theta}(x)\in\left[\cinf,3\cinf\right]$ thus $f_{\theta}(x)\geq \cinf$.
    \end{itemize}
\noindent For the last point, we have that for all $x\in[-\frac{\xM}{2},\frac{\xM}{2}]^c$, $f_\theta(x)=\omega(x)\leq \tilde M \psi(x)$. Moreover, for all $x\in[-\frac{\xM}{2},\frac{\xM}{2}]$, $f_\theta(x)\leq 3\cinf \leq 3 \cinf \sqrt{\pi}\exp(\xM^2 /4) \psi(x)$. Finally,
for all $x\in\R$, $f_\theta(x)\leq M(\cinf,\xM,\tilde M) \psi(x)$ with $M(\cinf,\xM,\tilde M):=\tilde M \vee 3 \sqrt{\pi}\cinf \exp(\xM^2/4)$.
\end{proof}

\begin{lem}~\label{LemmaSbeta2}
Let  $\beta\in[\underline\beta,\bar\beta]$. For all $\theta\in\{0,1\}^D$, the function
$\ln f_{\theta}$ is locally $\beta$-Hölder: for all $x,y$ such that $|x-y|\leq \frac{\xM}{4}$,
$$
    |(\ln f_\theta)^{(r)}(x) - (\ln f_\theta)^{(r)}(y)|\leq   L(\underline\beta,\bar\beta,\tilde L,\xM) r! |x-y|^{\beta-r}
$$
where $ L(\underline\beta,\bar\beta,\tilde L,\xM)$ does not depend on $D$. Moreover, there exists a constant $C(\underline\beta,\bar\beta,\tilde C,\xM)$, which can be taken identical for every $D$, such that for any integer $j=1,\dots,r$ and for all $D \in \N^*$,
$$
    \int_{\R} |(\ln f_\theta)^{(j)}(x)|^{\frac{2\beta+\tilde\varepsilon}{j}} f_\theta(x) dx \leq C(\underline\beta,\bar\beta,\tilde C,\xM),
$$
and
$$
    \int_{\R} |L(\underline\beta,\bar\beta,\tilde L,\xM)|^{\frac{2\beta+\tilde\varepsilon}{\beta}} f_\theta(x) dx \leq C(\underline\beta,\bar\beta,\tilde C,\xM).
$$
If $D$ is a positive even integer, for any integer $j=0,\dots,r$, $|(\ln f_\theta)^{(j)}(0)|\leq \ln(2\cinf).$
\end{lem}

\begin{proof}

Let  $j\in\{1,\dots,D\}$ and $1\leq t\leq r+1$. We start by upper bounding $ \mbox{sup}_{x \in I_j }|(\ln f_\theta)^{(t)}(x)|$.
According to
Lemma~\ref{LemLogDer}, for all $x \in I_j$,
$$
    (\ln f_\theta)^{(t)}(x) = f_\theta(x)^{-2^{t-1}} \underset{(\eta_0,\dots,\eta_{t})\in\Xi_{t}}{\sum} \rho(\eta_0,\dots,\eta_{t}) \underset{u=0}{\stackrel{t}{\prod}} \left(f_\theta^{(u)}(x)\right)^{\eta_u}
$$
with
$$
    \Xi_{t}=\left\{(\eta_0,\dots,\eta_{t})\in\mathbb{N}^{t+1}; \sum_{u=0}^{t}u\eta_u = t,\ \sum_{u=0}^{t}\eta_u =2^{t-1}\right\}.
$$
For all $u\in\{1,\dots,t\}$,
$$
   |f_\theta^{(u)}|\leq \frac{\cinf D^{-\beta}}{A} \left(\frac{D}{\xM}\right)^u \|\varphi^{(u)}\|_\infty \leq  \frac{\cinf D^{u-\beta}}{\xM^{u}}.
$$
Then, for all $(\eta_0,\dots,\eta_{t})\in\Xi_{t}$,
\begin{eqnarray*}
   \left|\prod_{u=0}^{t} (f_\theta^{(u)})^{\eta_u}\right|
   &\leq& D^{\sum_{u=1}^{t} u\eta_u - \beta\sum_{u=1}^{t}\eta_u} \cinf^{\sum_{u=1}^{t} \eta_u}  \xM^{-\sum_{u=1}^{t}u\eta_u} \times |f_\theta|^{\eta_0}\\
   &\leq& \cinf^{2^{t-1}-\eta_0} D^{t-\beta(2^{t-1} - \eta_0)} \xM^{-t}  \times |f_\theta(x)|^{\eta_0}
\end{eqnarray*}
since $\sum_{u=1}^{t} u\eta_u=t$ and $\sum_{u=1}^{t}\eta_u=2^{t-1} - \eta_0$.  
Since $f_\theta(x)\in\left[\cinf,3\cinf\right]$ and $2^{t-1}-\eta_0\geq 1$,
\begin{eqnarray*}
    |(\ln f_\theta)^{(t)}(x)|
    &\leq & \underset{(\eta_0,\dots,\eta_{t})\in\Xi_{t}}{\sum} |\rho(\eta_0,\dots,\eta_{t})| \left|\underset{u=1}{\stackrel{t}{\prod}}\left(f_\theta^{(u)}(x)\right)^{\eta_u}\right| |f_\theta(x)|^{\eta_0 - 2^{t-1}}\notag\\
    &\leq & \underset{(\eta_0,\dots,\eta_{t})\in\Xi_{t}}{\sum} |\rho(\eta_0,\dots,\eta_{t})| \cinf^{2^{t-1}-\eta_0} D^{t-\beta(2^{t-1} - \eta_0)}\xM^{-t} \cinf^{\eta_0 - 2^{t-1}}\notag\\
    &\leq & \underset{(\eta_0,\dots,\eta_{t})\in\Xi_{t}}{\sum} |\rho(\eta_0,\dots,\eta_{t})| D^{t-\beta(2^{t-1} - \eta_0)} \xM^{-t}.
\end{eqnarray*}
Denoting $\mathcal{B}(t) := \mbox{card}(\Xi_{t})$ and  $B(t) := \mbox{max}_{ (\eta_0,\dots,\eta_{t})\in\Xi_{t}}  |\rho(\eta_0,\dots,\eta_{t})| $, it leads to
\begin{equation}  \label{logftheta(t)}
 \mbox{sup}_{x \in I_j }|(\ln f_\theta)^{(t)}(x)|     \leq \mathcal{B}(t) B(t) D^{t-\beta} \xM^{-t}.
\end{equation}

We now use this preliminary result to prove that $\ln f_\theta$ is locally $\beta$-Holder. Let $(x,y)\in\R^2$ such that $|x-y|\leq \frac{\xM}{4}$.
\begin{itemize}
\item If $x,y\in[-\frac{\xM}{2},\frac{\xM}{2}]^c$,
\begin{eqnarray*}
|(\ln f_\theta)^{(r)}(x) - (\ln f_\theta)^{(r)}(y)|&=& |(\ln \omega)^{(r)}(x) - (\ln \omega)^{(r)}(y)|\\
                                                     &\leq& \tilde L r! |x-y|^{\beta-r}.
\end{eqnarray*}
since $\ln\omega$ is locally $\beta$-Holder with $\gamma_\omega=\frac{\xM}{4}$ and a constant $\tilde L$.

\item If $y\in[-\frac{\xM}{2},\frac{\xM}{2}]^c$ and $x\in I_j$:
    \begin{itemize}
     \item If $|x-y|<\frac{\xM}{4D}$ then $x\in I_j\backslash J_j$. Thus, $\ln f_\theta(x)=\ln \omega(x)$ and
    \begin{eqnarray*}
     |(\ln f_\theta)^{(r)}(x) - (\ln f_\theta)^{(r)}(y)|&=& |(\ln \omega)^{(r)}(x) - (\ln \omega)^{(r)}(y)|\\
                                                            &\leq& \tilde L r! |x-y|^{\beta-r}.
    \end{eqnarray*}

    \item If $\frac{\xM}{4D}\leq|x-y|<\frac{\xM}{4}$, $\ln \omega(y)=\ln\left(2\cinf\right)$ since $x\in [-3\xM/4,-\xM/2]\cup[\xM/2,3\xM/4]$ thus if $r\geq 1$,
    \begin{eqnarray*}
     |(\ln f_\theta)^{(r)}(x) - (\ln f_\theta)^{(r)}(y)|
     &\leq & \|(\ln f_\theta)^{(r)}\|_{\infty,[-\xM/2,\xM/2]} + \|(\ln \omega)^{(r)}\|_{\infty,[\xM/2,3\xM/4]}\\
     &\leq& \mathcal{B}(r) B(r) D^{r-\beta} \xM^{-r} \left(\frac{4D}{\xM}\right)^{\beta-r} |x-y|^{\beta-r} + 0\\
     &\leq& \frac{\mathcal{B}(r) B(r)}{r!} 4^{\beta -r} \xM^{-\beta}\ r!\ |x-y|^{\beta-r}
     \end{eqnarray*}
    and if $r=0$,
     \begin{eqnarray*}
     |(\ln f_\theta)(x) - (\ln f_\theta)(y)|
     &\leq & \left|\ln\left(2\cinf\right) - \ln\left(2\cinf + (2\theta_j-1)\varphi_j(y)\right)\right|\\
     &\leq & \left|- \ln\left(1+ (2\cinf)^{-1}(2\theta_j-1)\varphi_j(y)\right)\right|\\
     &\leq & \left| (2\cinf)^{-1}(2\theta_j-1)\varphi_j(y)\right|\\
     &\leq & (2\cinf)^{-1} \cinf D^{-\beta} (4D)^{\beta} \xM^{-\beta} |x-y|^\beta \\
     &\leq & 4^{\beta} \xM^{-\beta} |x-y|^\beta = \frac{\mathcal{B}(1) B(1)}{0!}4^{\beta} \xM^{-\beta}  0!\ |x-y|^{\beta-r}.
     \end{eqnarray*}
    \end{itemize}

\item For all $x,y\in[-\xM/2,\xM/2]$, $\exists! (j,j')\in\{1,\dots,D\}^2$ such that $x\in I_j$ and $y\in I_{j'}$.
    \begin{itemize}
    \item If $|x-y|\leq \frac{\xM}{4D}$,
            \begin{itemize}
            \item if $j'\neq j$, $x\in I_j\backslash J_j$ and $y\in I_{j'}\backslash J_{j'}$, thus
                $$
                    |(\ln f_\theta)^{(r)}(x) - (\ln f_\theta)^{(r)}(y)|=0.
                $$

            \item if $j'=j$,
           \begin{eqnarray*}
           |(\ln f_\theta)^{(r)}(x) - (\ln f_\theta)^{(r)}(y)|
           &\leq& |x-y|^{\beta-r} |x-y|^{r+1-\beta}\ \|\ln f_\theta^{(r+1)}\|_{\infty,[-\xM/2,\xM/2]}\\
           &\leq & \xM^{-\beta +r+1}(4D)^{\beta - r - 1} \frac{\mathcal{B}(r+1) B(r+1)}{r!} \frac{D^{r+1-\beta}}{\xM^{r+1}} r! |x-y|^{\beta-r}\\
           &\leq & \frac{\mathcal{B}(r+1) B(r+1)}{r!} 4^{\beta-r-1} \xM^{-\beta} r! |x-y|^{\beta-r}
           \end{eqnarray*}
           \end{itemize}
   \item If $\frac{\xM}{4D}<|x-y|<\frac{\xM}{4}$:
   if $r=0$,
   \begin{eqnarray*}
    |(\ln f_\theta)(x) - (\ln f_\theta)(y)|
     &= &  \left|\ln\left(\frac{1+(2\cinf)^{-1}(2\theta_j-1)\varphi_j(x)}{1+(2\cinf)^{-1}(2\theta_j-1)\varphi_j(y)}\right)\right|\\
     &\leq& \left|\frac{(2\cinf)^{-1}(2\theta_j-1) [\varphi_j(x)-\varphi_j(y)]}{1+(2\cinf)^{-1}(2\theta_j-1)\varphi_j(y)}\right|\\
     &\leq& \frac{2 \|\varphi_j\|_\infty}{\cinf}\\
     &\leq& 2 D^{-\beta} (4D)^\beta \xM^{-\beta} |x-y|^\beta\\
     &\leq& 2 4^\beta \xM^{-\beta}  |x-y|^\beta = 2 4^\beta \xM^{-\beta}  \frac{\mathcal{B}(1) B(1)}{0!} 0! |x-y|^\beta
    \end{eqnarray*}
    and if $r\geq 1$
    \begin{eqnarray*}
    |(\ln f_\theta)^{(r)}(x) - (\ln f_\theta)^{(r)}(y)|
     &\leq & 2 \|(\ln f_\theta)^{(r)}\|_{\infty,[-\xM/2,\xM/2]}\\
     &\leq& 2 \mathcal{B}(r) B(r) \frac{D^{r-\beta}}{\xM^r} \left(\frac{4D}{\xM}\right)^{\beta-r}|x-y|^{\beta-r}\\
     &\leq& 2 \frac{\mathcal{B}(r) B(r)}{r!} 4^{\beta} \xM^{-\beta} r! |x-y|^{\beta-r}.
    \end{eqnarray*}
   \end{itemize}
\end{itemize}
Finally, for all $\beta\in[\underline\beta,\bar\beta]$, for all $(x,y)\in\R^2$ such that $|x-y|<\frac{\xM}{4}$,
$$
 |(\ln f_\theta)^{(r)}(x) - (\ln f_\theta)^{(r)}(y)|\leq L(\underline\beta,\bar\beta,\xM) r! |x-y|^{\beta-r}
$$
with
$$
    L(\underline\beta,\bar\beta,\tilde L,\xM):= \tilde L\vee \underset{\beta\in[\underline\beta,\bar\beta]}{\max}\left(2 \frac{\mathcal{B}(\lceil\beta\rceil) B(\lceil\beta\rceil)}{\lfloor\beta\rfloor!} \left(\frac{4}{\xM}\right)^\beta \right).
$$
According to (\ref{logftheta(t)}), for any integer $j\in\{1,\dots,r\}$, $\|(\ln f_\theta)^{(j)}\|_{\infty,[-\xM/2,\xM/2]}\leq B(j) \mathcal{B}(j) \xM^{-j}$ thus it yields
\begin{eqnarray*}
\int_\R |(\ln f_\theta)^{(j)}(x)|^{\frac{2\beta+\tilde\varepsilon}{j}} f_\theta(x) dx
&\leq & \int_{[-\xM/2,\xM/2]^c }  \hskip-0.7cm  |(\ln \omega)^{(j)}(x)|^{\frac{2\beta+\tilde\varepsilon}{j}} \omega(x) dx + \left[B(j) \mathcal{B}(j) \xM^{-j}  \right]^{\frac{2\beta+\tilde\varepsilon}{j}} \int_{[-\xM/2,\xM/2]}  \hskip-0.7cm f_\theta(x) dx\\
&\leq & \tilde C + \left[B(j) \mathcal{B}(j)\xM^{-j}\right]^{\frac{2\beta+\tilde\varepsilon}{j}}.
\end{eqnarray*}
Thus there exists a constant $C(\underline\beta,\bar\beta,\tilde C,\tilde\varepsilon,\xM)$ such that for any integer $j\in\{1,\dots,r\}$,
$$
    \int_\R |(\ln f_\theta)^{(j)}(x)|^{\frac{2\beta+\varepsilon}{j}} f_\theta(x) dx \leq \tilde C + \underset{1\leq j\leq r+1}{\max}\left[B(j) \mathcal{B}(j)\right]^{\frac{2\beta+\tilde\varepsilon}{j}} \leq C(\underline\beta,\bar\beta,\tilde C,\tilde\varepsilon,\xM)
$$
and
$$
    \int_\R |L(\underline\beta,\bar\beta,\tilde L,\xM)|^{2+\frac{\tilde\varepsilon}{\beta}} f_\theta(x) dx = |L(\underline\beta,\bar\beta,\tilde L,\xM)|^{2+\frac{\tilde\varepsilon}{\beta}} \leq C(\underline\beta,\bar\beta,\tilde C,\tilde\varepsilon,\xM).
$$

The last point assumes that $D$ is even, thus $0\in I_{D/2}\backslash J_{D/2}$.  Then, $\ln f_\theta$ is equal to $\ln(2\cinf)$ in a neighborhood of $0$ and for all  $j\in\{1,\dots,r\}$, $|(\ln f_\theta)^{(j)}(0)| = 0$.
\end{proof}

Lemmas \ref{propMonotonicity} and \ref{LemmaSbeta2} show that for any positive even integer $D$ and for all $\beta \in [\underline \beta, \bar \beta]$, $ \mathcal{J}(\beta,D) \subset  \mathcal{H}\left(\beta, \Pc(\underline \beta, \bar \beta)\right) $
where
$$
    \Pc(\underline\beta,\bar\beta) = \left\{ \frac{\xM}{4} , \ln(2\cinf) , L(\underline\beta,\bar\beta,\tilde L,\xM) , \tilde \varepsilon , C(\underline\beta,\bar\beta,\tilde C,\tilde\varepsilon,\xM), \xM,\cinf, M(\cinf,\xM,\tilde M) \right\}.
$$


\subsection{Proof of Theorem~\ref{ThMinoration}}\label{subs:proofThMinor}
\begin{lem}\label{propHellinger}
Let $\theta,\theta'\in\{0,1\}^D$. The Hellinger
distance between two functions $f_{\theta}$ and $f_{\theta'}$ of $\mathcal{J}(\beta,D)$ fulfills
\begin{enumerate}
\item $d_H^2(f_{\theta},f_{\theta'}) \leq \frac{\cinf \xM}{8 A^2}D^{-2\beta}$,
\item $\forall \theta\neq \theta',\ d_H^2(f_{\theta},f_{\theta'}) \geq \cinf \xM (2A)^{-2} \delta(\theta,\theta') D^{-(2\beta + 1)}$ where
$\delta(\theta,\theta')=\underset{j=1}{\stackrel{D}{\sum}}\mathds{1}_{\theta_j
\neq \theta'_j}$ is the Hamming distance between $\theta$ and
$\theta'$.
\end{enumerate}
\end{lem}

\begin{proof}~\\
The Hellinger distance between $f_{\theta}$ and $f_{\theta'}$ can be decomposed as follows:
\begin{eqnarray*}
d_H^2(f_{\theta},f_{\theta'})
                              &=& \frac 1 2 \int_{[-\xM/2,\xM/2]}  \left[\sqrt{f_{\theta}(x)} - \sqrt{f_{\theta'}(x)}\right]^2 dx + \frac 1 2 \int_{[-\xM/2,\xM/2]^c}  \left[\sqrt{\omega(x)} - \sqrt{\omega(x)}\right]^2 dx\\
                              &=& \frac 1 2 \underset{j=1}{\stackrel{D}{\sum}} \int_{I_j} \left[\sqrt{2\cinf + (2\theta_j -1) \varphi_{j}(x)} - \sqrt{2\cinf + (2\theta'_j -1) \varphi_{j}(x)}\right]^2 dx.
\end{eqnarray*}
Since the quantity under the brackets is equal to zero if
$\theta_j= \theta'_j$, it gives
\begin{eqnarray*}
d_H^2(f_{\theta},f_{\theta'}) &=& \frac 1 2 \underset{j=1}{\stackrel{D}{\sum}} \int_{I_j} \left[\sqrt{2\cinf + \varphi_{j}(x)} - \sqrt{2\cinf - \varphi_{j}(x)}\right]^2 dx \ \mathds{1}_{\theta_j \neq \theta'_j}\\
                            &=& \frac 1 2 \underset{j=1}{\stackrel{D}{\sum}} \int_{I_j} \left[4\cinf - 2 \sqrt{\left(2\cinf\right)^2 - \varphi_{j}(x)^2}\right] dx \ \mathds{1}_{\theta_j \neq \theta'_j}.
\end{eqnarray*}
Note that $\left(\frac{\varphi_{j}(x)}{2\cinf}\right)^2 \leq 1$ for all $x\in I_j$ and $\|\varphi_{j}\|_\infty = \frac{\cinf D^{-\beta}}{A} \|\varphi\|_\infty\leq \cinf$.
Then,
$$
    \sqrt{\left(2\cinf\right)^2 - \varphi_{j}(x)^2} = 2\cinf \sqrt{1 - \left(\frac{\varphi_{j}(x)}{2\cinf}\right)^2} \geq \frac{1}{4} \left[ 1 - \left(\frac{\varphi_{j}(x)}{2\cinf}\right)^2\right]
$$
since  $\sqrt{1-y}\geq 1 - y$  for all $y\in [0,1]$. Thus,
\begin{eqnarray*}
\int_{I_j} \left[4\cinf - 2 \sqrt{(2\cinf)^2 - \varphi_{j}(x)^2}\right] dx
&\leq & \int_{I_j} \left[ 4\cinf - 4\cinf + \frac{\varphi_{j}^2(x)}{4\cinf}\right] dx\\
&\leq & (4\cinf)^{-1} \int_{I_j} \left[ \left(\frac{\cinf D^{-\beta}}{A}\right)^2 \varphi^2 \left(\frac{D}{\xM}(x+1) - (j-1)\right)\right] dx\\
&\leq & (4\cinf)^{-1} \left(\frac{\cinf D^{-\beta}}{A}\right)^2 \frac{\xM}{D}
\end{eqnarray*}
since $\int_{\R} \varphi^2(y) dy=1$. Finally,
\begin{eqnarray*}
d_H^2(f_{\theta},f_{\theta'})&\leq & (4\cinf)^{-1} \left(\frac{\cinf D^{-\beta}}{A}\right)^2 \frac{\xM}{D}  \frac 1 2 \delta(\theta,\theta')\\
&\leq & \frac{\cinf \xM}{8 A^2}D^{-2\beta}
\end{eqnarray*}
since $\delta(\theta,\theta')\leq D$.

\noindent For the lower bound, we have
$$
    \sqrt{\left(2\cinf\right)^2 - \varphi_{j}(x)^2} = 2\cinf \sqrt{1 - \left(\frac{\varphi_{j}(x)}{2\cinf}\right)^2} \leq 2\cinf \left[ 1 - \frac 1 2 \left(\frac{\varphi_{j}(x)}{2\cinf}\right)^2\right]
$$
since  $\sqrt{1-y}\leq 1 - \frac 1 2 y$ for all $y\in [0,1]$.
Thus,
\begin{eqnarray*}
\int_{I_j} \left[4\cinf - 2 \sqrt{\left(2\cinf\right)^2 - \varphi_{j}(x)^2}\right] dx
&\geq & \int_{I_j} \left[ 4\cinf - 4\cinf + \frac{ \varphi_{j}^2(x)}{2\cinf}\right] dx\\
&\geq & (2\cinf)^{-1} \left(\frac{\cinf D^{-\beta}}{A}\right)^2  \frac{\xM}{D} \int_{\R} \varphi^2(y) dy\\
&\geq & (2\cinf)^{-1}  \left(\frac{\cinf D^{-\beta}}{A}\right)^2 \frac{\xM}{D}
\end{eqnarray*}
and finally
\begin{eqnarray*}
d_H^2(f_{\theta},f_{\theta'})
&\geq & (2\cinf)^{-1} \left(\cinf\frac{D^{-\beta}}{A}\right)^2 \frac{\xM}{D} \frac 1 2 \underset{j=1}{\stackrel{D}{\sum}}\mathds{1}_{\theta_j \neq \theta'_j}\\
&\geq & \cinf \xM (2A)^{-2} D^{-(2\beta + 1)} \delta(\theta,\theta').
\end{eqnarray*}
\end{proof}

\begin{lem}\label{propKullback}
Let $\theta,\theta'\in\{0,1\}^D$. The
Kullback-Leibler divergence between two functions $f_{\theta}$ and
$f_{\theta'}$ of $\mathcal{J}(\beta,D)$ fulfills
$$
    \KL(f_{\theta},f_{\theta'}) \leq \frac{5 \cinf \xM}{4 A^2} D^{-2\beta}.
$$
\end{lem}

\begin{proof} The Kullback-Leibler divergence between $f_{\theta}$ and
$f_{\theta'}$ is given by
\begin{eqnarray*}
\KL(f_{\theta},f_{\theta'}) &=& \int_{\R} f_{\theta}(x) \ln\left(\frac{f_{\theta}(x)}{f_{\theta'}(x)}\right) dx\\
                                &=& \int_{[-\xM/2,\xM/2]} f_{\theta}(x) \ln\left(\frac{f_{\theta}(x)}{f_{\theta'}(x)}\right) dx + \int_{[-\xM/2,\xM/2]^c} \omega(x) \ln\left(\frac{\omega(x)}{\omega(x)}\right) dx\\
                                &=& \int_{[-\xM/2,\xM/2]} f_{\theta}(x) \ln\left(\frac{f_{\theta}(x)}{f_{\theta'}(x)}\right) dx.
\end{eqnarray*}
Then for all $x\in [-\xM/2,\xM/2]$ and for all $\theta\in\{0,1\}^D$,
$f_{\theta}(x)\in \left[\cinf,3\cinf\right]$
according to Lemma~\ref{propMonotonicity} thus
$\left\|\frac{f_{\theta}}{f_{\theta'}}\right\|_{\infty,[-1,1]}\leq
3$. According to Lemma 7.23 in \citet{Mas07},
$$ \KL(f_{\theta},f_{\theta'})
\leq  2 \left[2 + \ln\left(\left\|\frac{f_{\theta}}{f_{\theta'}}\right\|_{\infty}\right)\right] d_H^2\left(f_{\theta},f_{\theta'}\right). $$
Lemma~\ref{propMonotonicity} gives that  for all $x \in [-\xM/2,\xM/2]$, $f_{\theta}(x)\in \left[\cinf,3\cinf\right]$ and furthermore, $f_{\theta} = f_{\theta'}$ on $[-\xM/2,\xM/2]^c$. Thus,
\begin{eqnarray*}
\KL(f_{\theta},f_{\theta'})
&\leq & 2 \left[2 + \ln\left(  \sup_{[-\xM/2,\xM/2]}  \left|\frac{f_{\theta}(x)}{f_{\theta'}(x)}\right|  \right)   \right] d_H^2\left(f_{\theta},f_{\theta'}\right) \\
&\leq & 10\ d_H^2\left(f_{\theta},f_{\theta'}\right)\\
&\leq & \frac{5 \cinf \xM}{4 A^2} D^{-2\beta}
\end{eqnarray*}
according to Lemma~\ref{propHellinger}.

\end{proof}

\begin{proof}[Proof of Theorem~\ref{ThMinoration}]
The proof consists of applying Corollary~\ref{Cor2.19} given in Appendix~\ref{sect:AppendixMinor} with the space $\mathcal{J}(\beta,D)$, the Hellinger distance $d_H$, $p=2$ and the finite subset $\mathcal{C}=\{f_{\theta},\ \theta\in \Theta\}$ where $\Theta$ is the subset of $\{0,1\}^D$ provided by Lemma~\ref{Lemma4.76}. Then, it has to be checked that
$$
    n \underset{\theta,\theta'\in\Theta}{\max}\ \KL(f_{\theta},f_{\theta'}) \leq \kappa \ln|\Theta|.
$$
According to Lemma~\ref{Lemma4.76}, $\ln|\Theta|> \frac D 8$
and $\kappa \geq \frac 1 2$. Moreover,
$\KL(f_{\theta},f_{\theta'})\leq \frac{5 \cinf \xM}{4 A^2} D^{-2\beta}$
and thus $D$ is chosen such that
$$
    n \frac{5 \cinf \xM}{4 A^2} D^{-2\beta} \leq \frac{D}{16} \Leftrightarrow \frac{20 \cinf \xM n}{A^2}\leq D^{2\beta + 1}.
$$
Since $3\cinf \xM\leq 1$ then  $20 \cinf \xM n A^{-2} \leq \frac{20}{3} n\leq 7n$ and we finally choose $D=\min\{ 2 k ; k \in \N^*, (2k)^{2\beta +1}\geq 7 n\}$. It gives that for any estimator $\tilde s$,
\begin{eqnarray*}
\underset{\theta\in\Theta}{\sup}\ \E_s[d_H^2(f_{\theta},\tilde s)]
&\geq & 2^{-2} (1-\kappa) \left[\underset{\theta,\theta'\in\Theta,\theta\neq\theta'}{\min} d_H(f_{\theta},f_{\theta'})\right]^2\\
&\geq & 2^{-2} (1-\kappa) \cinf \xM (2A)^{-2} D^{-(2\beta+1)} \underset{\theta,\theta'\in\Theta,\theta\neq\theta'}{\min} \delta(\theta,\theta')\\
&\geq & 2^{-2} (1-\kappa) \cinf \xM (2A)^{-2} D^{-(2\beta+1)} \frac{D}{4}\\
&\geq & \frac{(1-\kappa) \cinf \xM}{A^2}\ 2^{-6 - 2\beta}\ (7n)^{-\frac{2\beta}{2\beta +1}}
\end{eqnarray*}
according to Lemma~\ref{Lemma4.76}.
\end{proof}

\section{Proof of Theorem~\ref{ThAdapt}} \label{ProofAdapt}
Under the hypotheses of Section \ref{subs:adaptative}, let
$\Pc(\underline \beta, \bar \beta)$ be the parameter set given in Proposition~\ref{prop:exist}. In order to prove Theorem~\ref{ThAdapt}, we start with the following lemma that makes the connection between the models $\S_m$ and the approximation result given in Theorem~\ref{ThApproximation}.

\begin{lem} \label{lem:adap}
There exists a positive constant $c_{\underline \beta, \bar \beta}$ such that for all $\beta\in[ \underline\beta  , \bar\beta]$ and for all $s \in \mathcal{H}\left(\beta,\mathcal P(\underline \beta, \bar \beta)\right)$,
$$
    \KL(s,\S_m) \leq  c_{\underline \beta, \bar \beta}\ \lambdai(m)^{ \beta}.
$$
\end{lem}

\begin{proof}
According to Theorem~\ref{ThApproximation}, the level $\bar \sigma(\beta)$ under which the approximation (\ref{ctrapproxKL}) is valid is a continuous function of $\beta$. Thus we can define the positive constant $\bar \sigma(\underline \beta, \bar \beta) :=  \underset{\beta  \in [\underline \beta, \bar \beta]}{\inf} \bar \sigma(\beta)$. Next, let
$$
    m_0(\underline \beta, \bar \beta) :=   \inf \left\{m \geq 2;\  \sqrt{\lambdai(m)}  < \bar \sigma(\underline \beta, \bar \beta)   \right\}
$$
and consider $m \geq m_0(\underline \beta, \bar \beta)$. Then Theorem~\ref{ThApproximation} can be applied for $\sigma  =  \sqrt{\lambdai(m)}$: for all $\beta\in[ \underline\beta  , \bar\beta]$ and for all $s \in \mathcal{H}\left(\beta,\mathcal P(\underline \beta, \bar \beta)\right)$, there exists a mixture $\mel$
with less than $ G_\beta \lambdai(m)^{-\frac 1 2}  \left|\ln \sqrt{\lambdai(m)}\right|^{\frac 3 2}$
components, with means belonging to $[-\mus(m) ,\mus(m) ]$
and with the same variance $\lambdai(m)$ for each component such that
$$
    \KL(s,\mel) \leq c_\beta\ \lambdai(m)^\beta .
$$
Since $G_\beta$ is a non decreasing function of $\beta$, the number of components is less than
\begin{eqnarray*}
G_{\bar \beta} \left(\sqrt{\lambdai(m)}\right)^{-1}  \left|\ln \sqrt{\lambdai(m)}\right|^{\frac 3 2}
& \leq  & G_{\bar \beta} \left[\frac{a_{\bar \beta}}{m}  \left( \ln m \right)^{\frac 3 2}\right]^{-1}  \left|\ln  \left\{ \frac{ a_{\bar \beta} }{ m }  \left( \ln m \right)^{\frac 3 2}    \right\}   \right|^{\frac 3 2} \\
& \leq  &  m \frac{G_{\bar \beta}}{a_{\bar \beta}}
  \left[  \frac{\ln  a_{\bar \beta}}{\ln m} + 1+  \frac 3 2 \frac{  | \ln \ln m| }{\ln m}    \right]^{\frac 3 2}   \\
& \leq  &  m
\end{eqnarray*}
according to the definition of $\sqrt{\lambdai(m)}$ and Condition (\ref{Cdttechn}). This shows that $\mel \in \S_m$ and thus  $\KL(s,\S_m) \leq  c_\beta\ \lambdai(m)^\beta$ for all  $m \geq m_0(\underline \beta, \bar \beta)$. Since $c_\beta $ is continuous in $\beta$, there exists $c_ {\underline \beta, \bar \beta}>0 $ such that for all $\beta\in[ \underline\beta  , \bar\beta]$, for all $s \in \mathcal{H}\left(\beta,\mathcal P(\underline \beta, \bar \beta)\right)$, and for all  $m \geq m_0(\underline \beta, \bar \beta)$,
\begin{equation}
    \KL(s,\S_m) \leq  c_ {\underline \beta, \bar \beta} \left[ \lambdai(m) \right] ^\beta .
    \label{KLsmel}
 \end{equation}
It remains to show the same result for $m \leq m_0(\underline \beta, \bar \beta)$ : let $t_m$ be a mixture of $\S_m$, for all $\beta\in[ \underline\beta  , \bar\beta]$ and for all $s \in \mathcal{H}\left(\beta,\mathcal P(\underline \beta, \bar \beta)\right)$,
 \begin{eqnarray*}
 \KL(s, \S_m)  & \leq &  \KL(s, t_m) \\
 & \leq & \int M \psi(x)  \ln \left(\frac{M \psi(x)}{t_m(x)}\right) <  + \infty.
\end{eqnarray*}
Then it can be easily shown that (\ref{KLsmel}) is valid for all $m \geq 1$ by changing the constant $c_ {\underline \beta, \bar \beta}$.
\end{proof}

\begin{proof}[Proof of Theorem~\ref{ThAdapt}]
In order to upper bound the right-hand side of the oracle inequality (\ref{OracleIneq}), we first control the constant $\mathcal{A}$ defined by (\ref{defA}) that depends on the parameters of the Gaussian mixture model $\S_m$ :
$$
    \mathcal{A}^2 \leq 4 \left\{\ln(6\pi e^2) + \pi + \ln\left(\mus(m) \sqrt{\frac{8}{c_1 \lambdai(m)}}\right) + \ln\left(\frac{144 \lambdas(m)}{\lambdai(m)}\right)\right\}.
$$
For the third term, we note that
\begin{eqnarray*}
    \ln\left(\mus(m)\sqrt{\frac{8}{c_1 \lambdai(m)}}\right) &=&  \frac 1 2 \ln\left( \frac{  4 \tilde (G_{\bar\beta})^2 }{ c_1}  \frac{|\ln\lambdai(m)|}{\lambdai(m)}\right) \\
    & \leq & c_{\bar \beta} \ln(m)
\end{eqnarray*}
since $\sqrt{\lambdai  (m)} := \frac{a_{\bar \beta}}{m} (\ln m)^{\frac 3 2}$. For the last term,
$$
    \ln\left(\frac{144 \lambdas(m)}{\lambdai(m)}\right) = \ln\left(\frac{144 \, m ^2 \, \lambdas}{(a_{\bar \beta})^2 (\ln m)^3}\right)\leq c_{\bar\beta} \ln(m)
$$
and thus $\mathcal{A}^2$ is upper bounded by $ c_{\bar \beta}\ln(m)$. For the observation of a $n$ sample, the model collection is indexed by $\mathcal{M}_n = \{2, \dots, n \}$ and then $m \leq n$. Thus for all $m \in \mathcal{M}_n$,
\begin{eqnarray*}
    \pen(m) &= & \kappa\frac{3m-1}{n}\left\{1+2\,\mathcal{A}^2 + \ln\left(\frac{1}{1\wedge\frac{D(m)}{n}\,\mathcal{A}^2}\right)\right\} \\
    &\leq &  c_{\bar\beta} \frac{m}{n}\left[ \ln n + \ln m \right]\\
    &\leq &  c_{\bar\beta} \frac{m}{n} \ln(n).
\end{eqnarray*}
According to Lemma~\ref{lem:adap} and the definition of $\lambdai(m)$, the oracle inequality is upper bounded by
\begin{eqnarray*}
    \E\left[d_H^{\,2}(s,\hat{s}_{\hat{m}})\right]
    &\leq& \mathcal{C}  \inf_{m \in \mathcal{M}_n}  \left[\KL(s,\S_{m})+\pen(m)+\frac{1}{n}\right] \\
    &\leq& c_{\underline\beta,\bar\beta} \inf_{m \in \mathcal{M}_n} \left[ \frac{(\ln m)^{3\beta} }{m^{2\beta}} + m \frac{\ln n}{n}\right].
\end{eqnarray*}
Let $m_n := \inf \left\{ m  \geq 2 \,  ; m \in \N  ;  \,  ;  \frac{(\ln m)^{3\beta} }{m^{2\beta}}  \leq  m \frac{\ln n}{n}   \right\}$.  Note that if $m_n = 2$, then $ \E\left[d_H^{\,2}(s,\hat{s}_{\hat{m}})\right] \leq  4  c_{\underline\beta,\bar\beta}  \frac{\ln n}{n}  $  and this case is completed. Assuming now that $m_n>2$, we want to check that $m_n \leq n$. According to the definition of $m_n$,
$$
    \frac{(m_n - 1) ^{2\beta + 1 }}{\left[ \ln (m_n - 1) \right]^{3\beta} }  <  \frac{n}{\ln n}
$$
thus
$$
    \left[  \frac{(m_n - 1) }{ \ln (m_n - 1)  } \right]^{ \square} <  \frac{n}{\ln n}
$$
where $\square = 3 \beta$ if $\beta >1$ and $\square = 2 \beta +1 $ otherwise. Next,
since $ \frac{(m_n - 1) }{ \ln (m_n - 1)  }  > 1 $,
$$
    \frac{(m_n - 1) }{ \ln (m_n - 1)  } <  \frac{n}{\ln n}
$$
in all cases. Assuming that $n \geq 3$, it leads that $ m_n \leq n$. Since $m_n \in \mathcal{M}_n$,
\begin{eqnarray*}
    \E\left[d_H^{\,2}(s,\hat{s}_{\hat{m}})\right]     &\leq&  2 c_{\underline\beta,\bar\beta}   \frac{\left[ \ln (m_n - 1) \right]^{3\beta} }{(m_n - 1) ^{2\beta  }} \\
     &\leq&  2 c_{\underline\beta,\bar\beta}  2^{2 \beta} \frac{\left[ \ln m_n \right]^{3\beta} }{m_n  ^{2\beta  }} \\
     & \leq&  2^{2 \beta+1}  c_{\underline\beta,\bar\beta}    \left[ \ln m_n \right]^{3\beta}  \left[   \frac{\ln n }{n} (\ln m)^{-3\beta} \right] ^{\frac{2 \beta}{2\beta+1}}  \\
     & \leq&   \tilde{c}_{\underline\beta,\bar\beta}  \,   n^{-\frac{2 \beta}{2\beta+1}}  \left(\ln n \right)^{\frac{5 \beta}{2\beta+1}}   .
\end{eqnarray*}
\end{proof}
\section{Conclusion}

In this paper, the penalized estimator $\hat s _{\hat m}$ defined in  \citet{MaugisMichel09} is shown to be adaptive to the regularity on some density classes $\H_{\beta}$ which elements are univariate densities whose logarithm is locally $\beta$-Hölder. To prove this result, the approximation result given in \citet{KRV} has been adapted to control the bias term between our Gaussian mixture models and the density classes $\H_{\beta}$. A lower bound for the minimax risk on the density classes $\H_{\beta}$ has also been stated to finally prove that our estimator reaches the minimax rate.

In \cite{MaugisMichel09}, a Gaussian mixture estimator, fulfilling an oracle inequality as (\ref{OracleIneq}), is proposed in the context of multivariate data clustering.
In a future work, it would be interesting to extend our adaptive result to this multivariate case. This requires to state an approximation result as Theorem~\ref{ThApproximation} on multivariate density classes which have to be determined, that is obviously a technical task.

\medskip

\section*{Acknowledgements}
The authors are grateful to Willem Kruijer and Judith Rousseau for helpful discussions.

\appendix
\section{Appendices for the approximation result}\label{sect:AppendixApprox}
\subsection{Gaussian kernel properties}

\begin{lem} \label{Lemma10}
Let $\psi_{(p)} (x) = C_p e ^{-|x|^p}$ for all reals $x$ where $C_p$ denotes the normalizing constant $\left[2\Gamma(1+\frac 1 p)\right]^{-1}$. Given a positive integer $u$,  let $\varphi_{u,p}$ be the $u$-fold convolution of $\psi_{(p)}$. Then, for any $t \geq 0$ and for all $H>0$, there exists a number $k'= k'(p,t,u,H)$ such that for all $ \sigma < 1$,
$$
    \int_{|x| > k' |\ln \sigma|^{1/p}} \varphi_{u,p}(x) |x|^{t} dx  = O_{\beta,H}(\sigma^H).
$$
Furthermore, $k'$ is a continuous function of  $H$.
\end{lem}

The reader is referred to Lemma 10 in \citet{KRV} for the proof of Lemma~\ref{Lemma10}.
Next lemma is a technical result used in Lemma~\ref{Lemma1} to prove the general case $\beta \geq 2$.
\begin{lem}~\label{LemmeSomme}
    For all positive integer $u$ and for all integer $k\geq u$,
    $$
        \sum_{j=1}^{k+1} (-1)^j \left(\substack{k+1\\j}\right) \nu_{j,2u} =0
    $$
    where $\nu_{j,h}$ is the $h$-th moment of the $j$-fold convolution of the Gaussian kernel $\psi$.
\end{lem}

\begin{proof}
Let $u=1$ and $k\in\mathbb{N}^\ast$. For all $j\in\{1,\dots,k+1\}$, let $(X_1,\dots,X_j)$ be a sample with density $\psi$. Then
$$
    \nu_{j,2}=\mathbb{E}[(X_1+\dots,X_j)^2]=\underset{q_1+\dots+q_j=2}{\sum} \frac{2!}{q_1!\dots q_j!} \mathbb{E}[X_1^{q_1}]\dots\mathbb{E}[X_j^{q_j}]=j\mathbb{E}[X_1^2]=j\nu_{1,2}
$$
since the odd moments of $\psi$ are equal to zero. Thus,
$$
    \sum_{j=1}^{k+1} (-1)^j \left(\substack{k+1\\j}\right) \nu_{j,2} =\sum_{j=1}^{k+1} (-1)^j \left(\substack{k+1\\j}\right) j \nu_{1,2} = (k+1) \nu_{1,2} \sum_{j=0}^{k} (-1)^j \left(\substack{k\\j}\right) =0.
$$
We assume now that the result is true until rank $u-1$. Let $k\geq u$ and note that
\begin{eqnarray*}
\sum_{j=1}^{k+1} (-1)^j \left(\substack{k+1\\j}\right) \nu_{j,2u}
&=& \left[\sum_{t=1}^{k+1} (-1)^t \left(\substack{k+1\\t}\right) \right] \nu_{1,2u} + \sum_{j=2}^{k+1} \left[\sum_{t=j}^{k+1}(-1)^t \left(\substack{k+1\\t}\right)\right] (\nu_{j,2u} -\nu_{j-1,2u}).
\end{eqnarray*}
Moreover,
$\nu_{j,2u}=\sum_{h=0}^{2u} \left(\substack{2u\\h}\right)\E[(X_1+\dots+X_{j-1})^{2u-h}]\E[X_j^h]=\sum_{p=0}^{u} \left(\substack{2u\\2p}\right) \nu_{j-1,2(u-p)}\nu_{1,2p}$
with the convention that $\nu_{h,0}=1$. Thus,
\begin{eqnarray}
\sum_{j=1}^{k+1} (-1)^j \left(\substack{k+1\\j}\right) \nu_{j,2u}
&=&\left[\sum_{t=1}^{k+1} (-1)^t \left(\substack{k+1\\t}\right) \right] \nu_{1,2u} + \sum_{j=2}^{k+1} \left[\sum_{t=j}^{k+1}(-1)^t \left(\substack{k+1\\t}\right)\right] \left(\sum_{p=1}^{u} \left(\substack{2u\\2p}\right) \nu_{j-1,2(u-p)}\nu_{1,2p}\right) \notag \\
&=&\left[\sum_{t=1}^{k+1} (-1)^t \left(\substack{k+1\\t}\right)  + \sum_{j=2}^{k+1}\sum_{t=j}^{k+1}(-1)^t \left(\substack{k+1\\t}\right)\right] \nu_{1,2u}\label{EqLemSomme1}\\
&+& \sum_{p=1}^{u-1} \left(\substack{2u\\2p}\right) \nu_{1,2p} \left\{\sum_{j=2}^{k+1} \left[\sum_{t=j}^{k+1}(-1)^t \left(\substack{k+1\\t}\right)\right] \nu_{j-1,2(u-p)}\right\}\label{EqLemSomme}.
\end{eqnarray}
It can be checked that the term inside the brackets of  (\ref{EqLemSomme1}) is null.
For (\ref{EqLemSomme}), noting that
\begin{eqnarray*}
\sum_{t=j+1}^{k+1} (-1)^t \left(\substack{k+1\\t}\right)
&=& \sum_{t=j+1}^{k} (-1)^t \left[\left(\substack{k\\t-1}\right)+\left(\substack{k\\t}\right)\right]+(-1)^{k+1}\\
&=& (-1)^k - (-1)^j\left(\substack{k\\j}\right)+(-1)^{k+1}=- (-1)^j\left(\substack{k\\j}\right),
\end{eqnarray*}
we have that for all $1\leq p\leq u-1$,
\begin{eqnarray*}
\sum_{j=2}^{k+1} \left[\sum_{t=j}^{k+1}(-1)^t \left(\substack{k+1\\t}\right)\right] \nu_{j-1,2(u-p)}
&=& -\sum_{j=1}^{k} (-1)^j\left(\substack{k\\j}\right) \nu_{j,2(u-p)}=0
\end{eqnarray*}
according to the induction assumption.
Finally, $\underset{j=1}{\stackrel{k+1}{\sum}} (-1)^j \left(\substack{k+1\\j}\right) \nu_{j,2u} =0$.
\end{proof}

\subsection{Measure discretization}\label{AppMesDiscrete}
The following result is adapted from Lemma 2 in \citet{GhosalVaart07b}.
It allows us to approximate a general Gaussian mixture by a finite Gaussian mixture with a limited number of components.

\begin{prop}~\label{PropMesureDiscrete}~\\
Let $F$ be a probability measure on $[-a,a]$ and $\sigma>0$ such that $\sigma<a$. Let $\varepsilon\in(0,\pi^{-\frac 1 2})$. Then there exists a
discrete distribution $F'$ on $[-a,a]$ with at most $54 a \sigma^{-1} e^2 \left[1\vee\ln\left(\frac{1}{\sqrt{\pi}\varepsilon}\right)\right]$ support points such that
$$
    \|F\ast \psi_{\sigma} - F'\ast\psi_{\sigma}\|_{\infty}\leq \frac{2\varepsilon}{\sigma}.
$$
\end{prop}

\begin{proof}
The interval $[-a,a]$ can be partitioned into $k=\lfloor\frac{2a}{\sigma}\rfloor$ disjoint consecutive subintervals $I_1,\dots,I_k$ of length $\sigma$ and a final subinterval $I_{k+1}$ of length $l\leq\sigma$: $I_i = [a_i,a_i+\sigma[,\ i=1,\dots,k$ and $I_{k+1}=[a_{k+1},a_{k+1}+l]$.\\
We decompose $F$ on this partition $F=\underset{i=1}{\stackrel{k+1}{\sum}} F(I_i) F_i$
where each $F_i$ is a probability measure concentrated on $I_i$. Then, $F\ast \psi_{\sigma}(x) = \underset{i=1}{\stackrel{k+1}{\sum}} F(I_i) (F_i\ast \psi_{\sigma})(x)$.
Let $Z_i$ be a random variable distributed according to $F_i$, and let $G_i$ be the law of $W_i = (Z_i - a_i)/\sigma$. Thus $G_i$ is a probability measure on
$[0,1]$ for $i = 1,\dots,k$ and on $[0, l/\sigma]\subset[0, 1]$ for $i = k + 1$. Lemma~\ref{MesureDiscreteLemme} is applied for each measure $G_i$ and with $D=\ln\left(\frac{1}{\sqrt{\pi}\varepsilon}\right)^{-\frac 1 2}$. We obtain discrete distributions $G'_i$ such that $\|G_i\ast \psi - G'_i\ast \psi\|_{\infty}\leq 2\varepsilon$. Let $F'_i$ be the law of $a_i + \sigma W'_i$ if $W'_i$ has law $G'_i$ and set $F' = \sum_{i=1}^{k+1}F(I_i)F'_i$.
We have
$$
    F_i\ast\psi_{\sigma}(x) = \mathbb{E}\left[\psi_{\sigma}(x-Z_i)\right] = \mathbb{E}\left[\frac 1 \sigma \psi\left(\frac{x-Z_i}{\sigma}\right)\right] = \mathbb{E}\left[\frac 1 \sigma \psi\left(\frac{x-a_i}{\sigma} - W_i\right)\right] = \frac 1 \sigma G_i\ast\psi\left(\frac{x-a_i}{\sigma}\right)
$$
and $F'_i\ast\psi_{\sigma}(x)=\frac 1 \sigma G'_i\ast\psi\left(\frac{x-a_i}{\sigma}\right)$. Thus
$$
    |F_i\ast\psi_{\sigma}(x) - F'_i\ast\psi_{\sigma}(x)| = \frac 1 \sigma \left|G_i\ast\psi\left(\frac{x-a_i}{\sigma}\right) - G'_i\ast\psi\left(\frac{x-a_i}{\sigma}\right)\right|
    \leq  \frac 1 \sigma \|G_i\ast \psi - G'_i\ast \psi\|_{\infty}
    \leq  \frac{2\varepsilon}{\sigma}.
$$
Then
\begin{eqnarray*}
    |F\ast\psi_{\sigma}(x) - F'\ast\psi_{\sigma}(x)| &=& \left|\sum_{i=1}^{k+1} F(I_i) \left[F_i\ast\psi_{\sigma}(x) - F'_i\ast\psi_{\sigma}(x)\right]\right| \\
  &  \leq&  \frac{2\varepsilon}{\sigma} \sum_{i=1}^{k+1} F(I_i) .
\end{eqnarray*}
Thus $\|F\ast\psi_{\sigma} - F'\ast\psi_{\sigma}\|_\infty \leq \frac{2\varepsilon}{\sigma}$
and the number of support points of the discrete distribution $F'$ is upper bounded by
\begin{eqnarray*}
\sum_{i=1}^{k+1} 18\left[1\vee \ln\left(\frac{1}{\sqrt{\pi}\varepsilon}\right)^{-1/2}\right]^2 e^2 \ln\left(\frac{1}{\sqrt{\pi}\varepsilon}\right)
&=& (k+1) 18 \left[1\vee \ln\left(\frac{1}{\sqrt{\pi}\varepsilon}\right)^{-1}\right] e^2 \ln\left(\frac{1}{\sqrt{\pi}\varepsilon}\right)\\
&\leq& 54 a \sigma^{-1} e^2 \left[1\vee \ln\left(\frac{1}{\sqrt{\pi}\varepsilon}\right)\right].
\end{eqnarray*}
\end{proof}

The following lemma is an adaptation of Lemma 3.1 in \citet{GhosalVaart01}. For this lemma, one introduces the inverse function of $\psi_\sigma(.)$ defined by $\psi_\sigma^{-1}(y)=\sigma \sqrt{-\ln(\sqrt{\pi}y)}$ on $(0,\pi^{-\frac 1 2}]$.
\begin{lem}\label{MesureDiscreteLemme}
Let $F$ be a probability measure on $[0,B]$. Let $\varepsilon\in(0,\pi^{-\frac 1 2})$ and let $D$ be a positive constant such that $B\leq D \psi^{-1}(\varepsilon)$. Then there exists a discrete distribution $F'$ on $[0,B]$ with at most $18(1\vee D)^2 e^2 \ln\left(\frac{1}{\sqrt{\pi}\varepsilon}\right)$ support points such that
$$
\|F\ast \psi - F'\ast\psi\|_{\infty}\leq 2\varepsilon.
$$
\end{lem}

\begin{proof}
Let $x_0$ be a positive constant which can be calibrated.
\begin{itemize}
\item[$\bullet$] Case 1: Suppose that $|x|\geq x_0$. Then,
\begin{eqnarray*}
|F\ast \psi(x) - F'\ast \psi(x)|
                               &\leq& \int_0^B \psi(x-u) d(|F-F'|(u))\\
                               &\leq & \frac{1}{\sqrt{\pi}} \int_0^B \exp[-|x-u|^2] d(|F-F'|(u)).
\end{eqnarray*}
If $x_0\geq 2B$ then $x_0-B\geq \frac {x_0} 2$. Thus, for all $|x|>x_0$ and $|u|\leq B$, $|x-u|\geq (x_0-B)^2\geq \frac{x_0^2}{4}$ and
$$
    |F\ast \psi(x) - F'\ast \psi(x)| \leq \frac{1}{\sqrt{\pi}} \int_0^B \exp[-|x-u|^2] d(|F-F'|(u)) \leq \frac{2}{\sqrt{\pi}} \exp\left[-\frac{x_0^2}{4}\right].
$$
If $\varepsilon \leq \pi^{-\frac 1 2}$, we choose $x_0$ such that $\exp\left[-\frac{x_0^2}{4}\right] \leq \varepsilon \sqrt{\pi} \Leftrightarrow x_0\geq 2 \sqrt{-\ln\left(\sqrt{\pi}\varepsilon\right)} = 2 \psi^{-1}(\varepsilon)$.
Finally if $x_0=2\max(B,\psi^{-1}(\varepsilon))$, then $\|F\ast \psi - F'\ast \psi\|_{\infty, [-x_0,x_0]^c}\leq 2 \varepsilon$.

\item[$\bullet$] Case 2: Suppose that $|x|\leq x_0$. By Taylor's expansion of $e^y$ and $k!\geq k^k e^{-k}$, we have for any $y<0$, $k>1$,
$$
    \left|e^y - \sum_{j=0}^{k-1}\frac{y^j}{j!}\right|\leq \left|\frac{ y ^k}{k!}\right|\leq \left(\frac{e|y|}{k}\right)^k.
$$
We use this inequality with $y=-x^2$ thus
$$
    \left|\psi(x) - \frac{1}{\sqrt{\pi}}\sum_{j=0}^{k-1}\frac{x^{2j}(-1)^j}{j!}\right|\leq \frac{1}{\sqrt{\pi}}\left(\frac{e x^2}{k}\right)^k.
$$
Then, it leads to
\begin{eqnarray}
\left|F\ast \psi(x) - F'\ast \psi(x)\right|
&\leq & \left|\int_0^B \psi(x-u) - \frac{1}{\sqrt{\pi}}\sum_{j=0}^{k-1}\frac{(x-u)^{2j}(-1)^j}{j!} d(F-F')(u)\right|\label{eq1lem10}\\
& & + \left|\int_0^B  \frac{1}{\sqrt{\pi}}\sum_{j=0}^{k-1}\frac{(x-u)^{2j}(-1)^j}{j!} d(F-F')(u)\right|.\label{eq2lem10}
\end{eqnarray}
The term (\ref{eq2lem10}) can be written
$$
\left|\int_0^B  \frac{1}{\sqrt{\pi}}\sum_{j=0}^{k-1}\frac{(x-u)^{2j}(-1)^j}{j!} d(F-F')(u)\right|
         = \left| \sum_{j=0}^{k-1}\sum_{t=0}^{2j}\frac{1}{\sqrt{\pi}j!} \left(\substack{ 2j\\t}\right) (-1)^j (-x)^{2j-t} \int_0^B u^t d(F-F')(u)\right|.
$$
According to Lemma A.1 in \citet{GhosalVaart01}, there exists a discrete distribution $F'$ with at most $2k-1$ support points such that $\int_0^B u^t dF'(u)=\int_0^B u^t dF(u)$ for all $1\leq t\leq 2k-2$. Finally, considering this discrete distribution $F'$, we obtain that (\ref{eq2lem10}) is null. For the term (\ref{eq1lem10}),
$$
    \left|\int_0^B \psi(x-u) -\frac{1}{\sqrt{\pi}}\sum_{j=0}^{k-1}\frac{(x-u)^{2j}(-1)^j}{j!} d(F-F')(u)\right|  \leq \int_0^B \frac{1}{\sqrt{\pi}}\left(\frac{e (x-u)^2}{k}\right)^k d(F+F')(u).
$$
Since $|x|<x_0$ and $0\leq u\leq B$, $|x-u|\leq |x|+|u|\leq x_0+B\leq \frac{3x_0}{2}$, we obtain that
$$
    \int_0^B \frac{1}{\sqrt{\pi}}\left(\frac{e (x-u)^2}{k}\right)^k d(F+F')(u) \leq \frac{2}{\sqrt{\pi}} \left(\frac{ 9\,e\, x_0^2}{4\,k}\right)^k.
$$
Moreover, since $x_0=2\max(B,\psi^{-1}(\varepsilon))$ and $B\leq D \psi^{-1}(\varepsilon)$, $x_0\leq 2 (1\vee D) \psi^{-1}(\varepsilon) = 2 (1\vee D) \sqrt{\ln\left(\frac{1}{\sqrt{\pi}\varepsilon}\right)}$. Thus
\begin{eqnarray*}
\frac{2}{\sqrt{\pi}} \left(\frac{9\,e\, x_0^2}{4\, k}\right)^k
        &\leq& \frac{2}{\sqrt{\pi}} \left[\frac{9\,e\, (1\vee D)^2}{k} \ln\left(\frac{1}{\sqrt{\pi}\varepsilon}\right) \right]^k\\
        &\leq& \frac{2}{\sqrt{\pi}} \exp\left[-k \left\{\ln\left(\frac{9\,e\,(1\vee D)^2}{k}\right)- \ln\left[\ln\left(\frac{1}{\sqrt{\pi}\varepsilon}\right)\right] \right\}\right].
\end{eqnarray*}
We have that $\varepsilon\leq \pi^{-\frac 1 2}\leq 1$ and we choose $k$ such that $k\geq \ln\left(\frac{1}{\varepsilon}\right)$ and
$
\ln\left(\frac{9e(1\vee D)^2}{k}\right)- \ln\left[\ln\left(\frac{1}{\sqrt{\pi}\varepsilon}\right)\right] \geq 1.
$
This is the case if $k=9(1\vee D)^2 e^2 \ln\left(\frac{1}{\sqrt{\pi}\varepsilon}\right)$.
Finally, the term (\ref{eq1lem10}) is upper bounded by $2 \pi^{-\frac 1 2}\varepsilon \leq 2 \varepsilon$ and $\|F\ast \psi - F'\ast \psi\|_{\infty, [-x_0,x_0]}\leq 2 \varepsilon$.
\end{itemize}
\end{proof}


\subsection{Technical results for $f$, $f_k$, $g_k$, $h_k$ and their convolutions}
The following lemma allows to bound the derivative functions of $\ln f$. It is based on the smoothness assumptions (\ref{Cond1}) and (\ref{Cond1b}) and is used in the proof of Lemma~\ref{Lemma1}.

\begin{lem}\label{LemMajlu}~\\
    For all $j \in \{0,\dots,r \}$ and for all $n \in \Z$, there exists a constant $0 < l^+_{j,n} < \infty $ such that for all $f \in \mathcal{H}\left( \beta,\Pc\right)$,
        \begin{equation*}
        \sup_{y \in [n \gamma, (n+1) \gamma]} \left| (\ln f)^{(j)}(y) \right| \leq l^+_{j,n} .
        \label{Controlel_u}
        \end{equation*}
\end{lem}

\begin{proof}
We first prove Lemma \ref{LemMajlu} on $[-\gamma,\gamma]$. For all $j \in \{1,\dots,r\}$,  all $f \in \mathcal{H}\left( \beta,\Pc \right)$ and all $y \in [-\gamma,\gamma]$, there exists $\tilde{y} \in [-|y|,|y|]$ such that
$$
    (\ln f)^{(j)}(y)=\sum_{u= 0} ^{r-j} \frac{(\ln f)^{(j+u)}(0)}{u !} y ^{u}  +   \frac{y^{r-j}}{(r-j)!} \left[ (\ln f)^{(r)}(\tilde{y})  -  (\ln f)^{(r)}(0)  \right].
$$
Thus,
\begin{eqnarray*}
| (\ln f)^{(j)}(y) |& \leq & \sum_{u= 0} ^{r-j}  \frac{ \left| (\ln f)^{(j+u)} (0)\right|}{u !} |y| ^{u} + \frac{r!}{(r-j)!}  |y| ^{r-j}  L(0) |\tilde{y}|^{\beta-r}  \\
                & \leq & l^+  \sum_{u= 0} ^{r-j}   \frac{\gamma^u}{u !}+ \frac{r!}{(r-j)!}  L(0) \gamma^{\beta-j}  \\
                & \leq & l^+  \exp(\gamma)+ \frac{r!}{(r-j)!}  L(0) \gamma^{\beta-j} := l^+_{j,-1} = l^+_{j,0}
\end{eqnarray*}
and  Lemma \ref{LemMajlu} is proved for $n=-1$ and $n =0$. Now, assume that Lemma~\ref{LemMajlu} is valid for  $n-1\geq 0$. Then, proceeding as before, for all $j \in \{1,\dots,r\}$ and all $y \in [n\gamma,(n+1)\gamma]$, there exists  $\tilde{y} \in [n \gamma,(n+1) \gamma]$ such that
\begin{eqnarray*}
(\ln f)^{(j)}(y)& = & \sum_{u= 0} ^{r-j-1} \frac{(\ln f)^{(j+u)}( n \gamma)}{u !} (y-n\gamma) ^{u}  +   \frac{(y-n\gamma)^{r-j}}{(r-j)!} \left[ (\ln f)^{(r)}(\tilde{y})  -  (\ln f)^{(r)}(n\gamma)  \right]
\end{eqnarray*}
and thus
\begin{eqnarray*}
| (\ln f)^{(j)}(y) |& \leq & \sum_{u= 0} ^{r-j}\frac{ \left| (\ln f)^{(j+u)} (n \gamma) \right|}{u !} \gamma ^{u} +  \frac{r!}{(r-j)!}  \gamma ^{\beta-j}  L(n\gamma) \\
          & \leq &  \sum_{u= 0} ^{r-j}  l_{j+u,n-1}^+  \frac{\gamma^u}{u !}+ \frac{r!}{(r-j)!}  L(n \gamma) \gamma^{\beta-j}   := l^+_{j,n}.
\end{eqnarray*}
Finally, Lemma \ref{LemMajlu} is proved for all $n \in \N$ and a similar proof gives this result for all $n \in \Z\backslash\N$.
\end{proof}

\begin{lem}~\label{Lemfk}
    Let $f_0=f$ and $\forall k\in\mathbb{N}^\ast,\ f_{k+1}=f -\Delta_\sigma f_k$ with $\Delta_\sigma f_k=K_\sigma f_k - f_k$.
    \begin{enumerate}
    \item For all $x\in\R$, $f_k(x)=\underset{i=0}{\stackrel{k}{\sum}} \left(\substack{k+1\\i+1}\right) (-1)^i K_\sigma^i f(x).$
    \item For all $k\in\mathbb{N}$, $\int_\mathbb{R} f_k(x) dx =1$.
    \item For all $i\in\mathbb{N}$ and for all $x\in\R$, $K_\sigma^i f(x)\leq \frac{M}{\sqrt{\pi}}$ and thus $|f_k(x)|\leq (2^{k+1}-1)\frac{M}{\sqrt{\pi}}.$
    \end{enumerate}
\end{lem}

\begin{proof}
The first result is trivial for $k=0$. For $k=1$, we remark that $f_1(x)=f(x)-\Delta_\sigma f(x)=2f(x) - K_\sigma f(x)$.
Then recursively, we have
\begin{eqnarray*}
f_{k+1}(x)
           &=& f(x) - K_\sigma f_k(x) + f_k(x)\\
           &=& f(x) - \underset{i=0}{\stackrel{k}{\sum}} \left(\substack{k+1\\i+1}\right) (-1)^i K_\sigma^{i+1} f(x) + \underset{i=0}{\stackrel{k}{\sum}} \left(\substack{k+1\\i+1}\right) (-1)^i K_\sigma^i f(x)\\
           &=& f(x) + \underset{j=1}{\stackrel{k+1}{\sum}} \left(\substack{k+1\\j}\right) (-1)^j K_\sigma^{j} f(x) + \underset{i=0}{\stackrel{k}{\sum}} \left(\substack{k+1\\i+1}\right) (-1)^i K_\sigma^i f(x)\\
           &=& \left(\substack{k+1\\k+1}\right) (-1)^{k+1} K_\sigma^{k+1} f(x) + \underset{j=0}{\stackrel{k}{\sum}} \left[\left(\substack{k+1\\j}\right)+ \left(\substack{k+1\\j+1}\right) \right] (-1)^j K_\sigma^j f(x).
\end{eqnarray*}
Since $\left(\substack{k+1\\j}\right)+ \left(\substack{k+1\\j+1}\right) =\left(\substack{k+2\\j+1}\right)$
and $\left(\substack{k+1\\k+1}\right)=\left(\substack{k+2\\k+2}\right)$, we have
$$
f_{k+1}(x)=\underset{i=0}{\stackrel{k+1}{\sum}} \left(\substack{k+2\\i+1}\right) (-1)^i K_\sigma^i f(x).
$$
Consequently, for all $k\in\mathbb{N}$, $\int_\mathbb{R} f_k(x) dx = \underset{i=0}{\stackrel{k}{\sum}} \left(\substack{k+1\\i+1}\right) (-1)^i \int_\mathbb{R} K_\sigma^i f(x) dx.$
Moreover, it can be easily proved by induction that for all nonnegative integer $i$, $\int_\mathbb{R} K_\sigma^i f(x) dx=1$.
Thus,
$$
\int_\mathbb{R} f_k(x) dx = \underset{i=0}{\stackrel{k}{\sum}} \left(\substack{k+1\\i+1}\right) (-1)^i = - \underset{i=0}{\stackrel{k+1}{\sum}} \left(\substack{k+1\\i}\right) (-1)^i + \left(\substack{k+1\\0}\right) = (1-1)^{k+1} +1=1.
$$
For the third result, according to Condition (\ref{Cond2b}), $f(x)\leq M \psi(x)\leq \frac{M}{\sqrt{\pi}}$. And by induction,
$$
    K_\sigma^i f(x) =\int_\R K_\sigma^{i-1}f(u) \psi_\sigma(x-u) du \leq \frac{M}{\sqrt{\pi}}\int_\R \psi_\sigma(x-u) du \leq \frac{M}{\sqrt{\pi}}.
$$
Finally, $|f_k(x)|\leq \underset{i=0}{\stackrel{k}{\sum}} \left(\substack{k+1\\i+1}\right)  K_\sigma^i f(x)\leq (2^{k+1}-1)\frac{M}{\sqrt{\pi}}.$
\end{proof}

\begin{lem} \label{BoiteOutilBornes0}
Let $\beta>0$ and $k\in\mathbb{N}$ such that $\beta\in(2k,2k+2]$.
Let $f$ be a density function belonging to $\H(\beta,\Pc)$ where $\Pc = \left\{\gamma,l^+,L,\varepsilon,C,\xM,\cinf,M\right\}$.
\begin{enumerate}
\item Let $\bar{\sigma}>0$ such that if $Y$ is distributed from a centered Gaussian density with variance $\bar\sigma^2$, then $P(0<Y<2\xM)=\frac 1 3$. For all $\sigma<\bar{\sigma}$,
    \begin{equation}
    K_\sigma f(x) \geq \frac{\cinf \sqrt{\pi}}{3 M} f(x).
    \label{Ksigmafgeqf}
    \end{equation}
\item There exists $\bar \sigma(\beta)>0$ and $A_\beta>0$ such that for all $\sigma < \bar \sigma (\beta) $,
    \begin{equation*}
    K_\sigma h_k(x) \geq \frac{\cinf \sqrt{\pi}}{6M(1+A_{\beta} \sigma^{2\beta})}  f(x).
    \label{Ksigmahkgeqf}
    \end{equation*}
    Furthermore, $\bar \sigma (\beta)$ can be chosen as a continuous function of $\beta$.
\end{enumerate}
\end{lem}

\begin{Rem}
The first result of Lemma~\ref{BoiteOutilBornes0} is based on the monotonicity assumption on $f$. It comes from Remark 3 of \citet{GGR99}.
In the second result, the constants $\bar \sigma(\beta)$ and $A_{\beta}$ are due to the result (\ref{W24}) in Lemma~\ref{Lemme2}.
\end{Rem}

\begin{proof} For the first point, let $\sigma<\bar{\sigma}$ and $Z$ be a standard centered Gaussian random variable.
\begin{itemize}
\item[$\bullet$] If $x\in[-\xM,\xM]$,
    \begin{eqnarray*}
    K_\sigma f(x) &\geq & \int_{-\xM}^{\xM} f(u) \psi_{\sigma}(x-u) du\\
                          &\geq & \cinf \int_{\frac{x-\xM}{\sigma}}^{\frac{x+\xM}{\sigma}} \psi\left(z\right) dz\\
                          &\geq & \cinf \left\{P\left(\frac{x-\xM}{\sigma}\leq Z \leq 0\right) + P\left(0\leq Z \leq \frac{2\xM}{\sigma}\right) - P\left(\frac{x+\xM}{\sigma}\leq Z \leq \frac{2\xM}{\sigma}\right)\right\}\\
                          &\geq & \cinf \left\{\frac 1 3 + P\left(0\leq Z \leq \frac{\xM-x}{\sigma}\right) - P\left(\frac{x+\xM}{\sigma}\leq Z \leq \frac{2\xM}{\sigma}\right)\right\} \\
			  &\geq& \frac{ \cinf}{3}.
    \end{eqnarray*}
According to Condition (\ref{Cond2b}), $f(x)\leq M\psi(x)\leq \pi^{-1/2} M$ for all $x\in\R$. Then for all $x\in[-\xM,\xM]$,
$$
    K_\sigma f(x) \geq \frac{\cinf\sqrt{\pi}}{3 M} f(x).
$$

\item[$\bullet$] If $x>\xM$,
\begin{eqnarray*}
K_\sigma f(x)
                      &\geq & f(x) \int_0^{\frac{x+\xM}{\sigma}} \psi(u) du\\
                      &\geq & f(x) \left\{P\left(0\leq Z \leq \frac{2\xM}{\sigma}\right) + P\left(\frac{2\xM}{\sigma}\leq Z \leq \frac{x+\xM}{\sigma}\right)\right\}\\
                      &\geq & f(x) P\left(0\leq Z \leq \frac{2\xM}{\sigma}\right)\\
                   &\geq & \frac 1 3 f(x)  .
\end{eqnarray*}
In the same way, for all $x<-\xM$, $K_\sigma f(x) \geq  \frac 1 3 f(x)$.
\end{itemize}
Finally, since $M \pi^{-1/2}\geq \cinf$, $ K_\sigma f(x) \geq  \frac{\cinf\sqrt{\pi}}{3 M} f(x) $  for all $x \in \R$.

For the second point, we take $H_1 = 4 \beta +1 $ as in Remark~\ref{Rem:W24bis}; let  $\sigma < \bar{\sigma} (\beta)$ in order to have (\ref{W24bis}) and (\ref{Ksigmafgeqf}). Then, for all $\sigma < \bar \sigma (\beta) $,  $\int_\R g_k(u) du \leq 1+A_{\beta} \sigma^{2\beta}$ and since $g_k(x)\geq \frac 1 2 f(x)$ for all $x\in\R$,
\begin{eqnarray*}
K_\sigma h_k(x) &=& \int_\R h_k(u) \psi_\sigma(x-u) du\\
                &\geq & \int_\R \frac{f(u)}{2 \int g_k} \psi_\sigma(x-u) du\\
                &\geq & \frac{\cinf\sqrt{\pi}}{6M(1+A_{\beta} \sigma^{2\beta})} f(x) .
\end{eqnarray*}
\end{proof}

\begin{lem}\label{BoiteOutilBornes}
Let $p \in (0,1)$. For all $x\in \R$, we have that
\begin{itemize}
\item[$\bullet$] for all $i \in\mathbb{N}$ and for all $\sigma < 1 - p^{1/i}$,
        $ K_{\sigma}^i f(x) \leq M \left(\frac{2}{\sqrt{3}}\right)^{i}  \psi(p  x)$.
\item[$\bullet$] for all $\sigma < 1 - p^{1/k}$,
$$
    \max\left(f_k (x),g_k (x),\frac 1 2 h_k(x)\right)\leq   2 M \left(\frac{4}{\sqrt{3}}\right)^{k} \psi(p  x).
$$
\end{itemize}
\end{lem}

\begin{proof}
The control of $K^i_{\sigma} f$ can be proved by applying successively Lemma~\ref{reccontrole} to $f$, $K_{\sigma} f$, ..., $K_{\sigma}^{i-1} f$ with $q_1=p^{k/i}$ and $q_2=p^{1/i}$ for each step $k$. It finally gives that $K_{\sigma}^{i} f (x) \leq M \left(\frac{4}{\sqrt{3}}\right)^{i} \psi(p x)$ for all $x$ in $\R$ and for all $\sigma < 1 - p^{1/i}$. This control on $K^i_{\sigma} f$ together with Lemma \ref{Lemfk} give the control on $f_k$. According to the definition of $g_k$ and previous results,
$$
    g_k(x)\leq 2M\left(\frac{4}{\sqrt{3}}\right)^{k} \psi(px) \mathds{1}_{J_{\sigma,k}}(x) + M\psi(x)\mathds{1}_{J_{\sigma,k}^c}(x)\leq 2M\left(\frac{4}{\sqrt{3}}\right)^{k} \psi(px).
$$
Finally, $h_k(x) = g_k(x)/ \int_\R g_k(y) dy \leq 2 g_k(x) \leq 4M\left(\frac{4}{\sqrt{3}}\right)^{k} \psi(px)$.
\end{proof}

\begin{lem} \label{reccontrole}
Let $f$ be a positive application on $\R$ such that there exists $M> 0 $ and $q_1 \in (0,1]$ such that for all $x \in \R$, $ f(x) \leq M \psi(q_1 x)$. Then for all $q_2 \in (0,1)$ and all $\sigma\in (0,1- q_2^2)$,
$$
    \forall x \in \R, \  K_{\sigma} f(x) \leq \frac{2}{\sqrt{3}} M \psi(q_1 q_2 x).
$$
\end{lem}

\begin{proof}
Let $\sigma\in (0,1- q_2^2)$. For all $x\in\R$,
\begin{eqnarray*}
K_{\sigma} f (x) & =    & \int_{\R} f(u) \psi_{\sigma}(x-u) du \\
                 & \leq & \frac{M}{\pi}  \int_{\R} \exp\left[- q_1^2 ( \sigma y - x )^2\right] \exp(-y^2) dy \\
                 & \leq & \frac{M}{\pi} \exp\left[ - q_1^2 x^2 ( 1 - \sigma) \right]  \int_{\R} \exp(-y^2+ q_1^2 y^2\sigma (1-\sigma)) dy \\
                 & \leq & \frac{M}{\pi} \exp\left[ - q_1^2 q_2^2 x^2  \right]  \int_{\R} \exp\left(-y^2+\frac 1 4 q_1^2 y^2\right)dy \\
                 &\leq &  \frac{2M}{\sqrt{3}} \psi(q_1 q_2 x).
\end{eqnarray*}
\end{proof}

\section{Appendices for the lower bound result}\label{sect:AppendixMinor}
The two following results are crucial for establishing the lower bound: 
The first one is the so-called Varshamov-Gilbert's lemma and the second one is a corollary of a lemma given in \citet{Birge05}.
They correspond to Lemma 4.7 and Corollary 2.19 in \citet{Mas07} respectively.

\begin{lem}\label{Lemma4.76}
Let $\{0, 1\}^D$ be equipped with Hamming distance $\delta$. Given
$\alpha\in (0, 1)$, there exists some subset $\Theta$ of $\{0,
1\}^D$ with the following properties
$$\left\{\begin{array}{l}
\delta(\theta,\theta') > \frac{(1-\alpha) D}{2} \textrm{ for every } (\theta,\theta')\in\Theta, \theta\neq \theta'\\
\ln|\Theta| > \frac{\rho D}{2}
\end{array}\right.
$$
where $\rho = (1 + \alpha)\ln(1+\alpha) + (1 - \alpha) \ln (1 -
\alpha)$. In particular $\rho>\frac 1 4$ when $\alpha=\frac 1 2$.
\end{lem}

\begin{cor}\label{Cor2.19}
Let $(S, d)$ be some pseudo-metric space, $\{\mathbb{P}_s, s \in
S\}$ be some statistical model. Let $\kappa$ denote an absolute
constant \citep[given in Corollary 2.18 of ][]{Mas07}. Then for any estimator $\tilde s$ and
any finite subset $\mathcal{C}$ of $S$ such that $\underset{s,t\in
\mathcal{C}}{\max} \KL(\mathbb{P}_s,\mathbb{P}_t) \leq \kappa
\ln|\mathcal{C}|$, the following lower bound holds for every $p>1$

$$
    \underset{s\in \mathcal{C}}{\sup}\ \mathbb{E}_s[d^p(s,\tilde s)] \geq 2^{-p} (1-\kappa) \left[\underset{s,t\in \mathcal{C}, s\neq t}{\min}d(s,t)\right]^p.
$$
\end{cor}

The following lemma, used to prove Proposition~\ref{prop:exist}, gives an expression of the derivatives of the logarithm of a function. 

\begin{lem}~\label{LemLogDer}
Let $i\in\N^*$ and let $t$ be a strictly positive function, $t\in\mathcal{C}^i$. Then
$$
    (\ln t)^{(i)}(x) = \frac{P_i(x)}{t(x)^{2^{i-1}}}
$$
where
$$
  P_i(x)=\sum_{(\eta_0,\ldots,\eta_i)\in \Xi_i} \rho(\eta_0,\ldots,\eta_i) \prod_{j=0}^i \left[t^{(j)}(x)\right]^{\eta_j}
$$
with
$$
  \Xi_i=\left\{(\eta_0,\ldots,\eta_i)\in \mathbb{N}^{i+1}; \sum_{j=0}^i \eta_j = 2^{i-1},\ \sum_{j=0}^i j\eta_j=i\right\}
$$
and $\rho(\eta_0,\ldots,\eta_i)$'s are the polynomial coefficients.
\end{lem}

\begin{proof}
The result is trivial for $i=1$. Assume that $t$ is $\mathcal{C}^{i+1}$ and that the result is valid for the $i$-th derivative. Then
$$
    (\ln t)^{(i+1)}(x) = \frac{t(x)^{2^{i-1}}P_i(x)' - 2^{i-1}t(x)'t(x)^{2^{i-1}-1} P_i(x)}{t(x)^{2^{i}}}=\frac{\diamondsuit}{t(x)^{2^{i}}}
$$
with
\begin{eqnarray*}
\diamondsuit&=& t(x)^{2^{i-1}}\left\{\underset{(\eta_0,\ldots,\eta_i)\in \Xi_i}{\sum} \rho(\eta_0,\ldots,\eta_i) \left[\underset{j=0}{\stackrel{i}{\sum}} \eta_j \frac{t(x)^{(j+1)}}{t(x)^{(j)}}\underset{u=0}{\stackrel{i}{\prod}} \left(t(x)^{(u)}\right)^{\eta_u}\right]\right\}\\
& & \hspace*{2cm} - 2^{i-1}t(x)'t(x)^{2^{i-1}-1}\underset{(\eta_0,\ldots,\eta_i)\in \Xi_i}{\sum} \rho(\eta_0,\ldots,\eta_i) \underset{j=0}{\stackrel{i}{\prod}}\left(t(x)^{(j)}\right)^{\eta_j}\\
\end{eqnarray*}
\begin{eqnarray*}
\diamondsuit &=& \underset{(\eta_0,\ldots,\eta_i)\in \Xi_i}{\sum} \rho(\eta_0,\ldots,\eta_i) \left[\underset{j=0}{\stackrel{i}{\sum}} \eta_j t(x)^{2^{i-1}}\frac{t(x)^{(j+1)}}{t(x)^{(j)}}\underset{u=0}{\stackrel{i}{\prod}} \left(t(x)^{(u)}\right)^{\eta_u}\right]\\
& & \hspace*{2cm}- \underset{(\eta_0,\ldots,\eta_i)\in \Xi_i}{\sum} \rho(\eta_0,\ldots,\eta_i) 2^{i-1}t(x)'t(x)^{2^{i-1}-1}\underset{j=0}{\stackrel{i}{\prod}}\left(t(x)^{(j)}\right)^{\eta_j}.
\end{eqnarray*}
Let $\tilde \eta_j$ denotes the new power of the $j$-th derivative for $j=0\ldots,i+1$. In the second sum, we have that
$\tilde \eta_0=2^{i-1}-1+\eta_0$, $\tilde\eta_1=\eta_1+1$, $\tilde\eta_j=\eta_j$ for all $j=2,\ldots,i$ and $\tilde\eta_{i+1}=0$
thus $\underset{j=0}{\stackrel{i+1}{\sum}}\tilde \eta_j =\underset{j=0}{\stackrel{i}{\sum}}\eta_j +1+2^{i-1}-1 = 2^i$ and $\underset{j=0}{\stackrel{i+1}{\sum}}j \tilde \eta_j =\underset{j=0}{\stackrel{i}{\sum}}j \eta_j + 1 = i+1$.
In the first sum,
\begin{itemize}
        \item if $j<i$: $\tilde \eta_0=2^{i-1}+\eta_0$, $\tilde\eta_j=\eta_j-1$, $\tilde\eta_{j+1}=\eta_{j+1}+1$, $\forall u\in\{1,\ldots,i\}\backslash\{j,j+1\}, \tilde\eta_u=\eta_u$ and $\tilde\eta_{i+1}=0$ thus $\underset{u=0}{\stackrel{i+1}{\sum}}\tilde \eta_u =\underset{u=0}{\stackrel{i}{\sum}}\eta_u +1+2^{i-1}-1 = 2^i$ and $\underset{u=0}{\stackrel{i+1}{\sum}}u \tilde \eta_u =\underset{u=0}{\stackrel{i}{\sum}}u \eta_u + j+1-j = i+1$.

        \item if $j=i$: $\tilde \eta_0=2^{i-1}+\eta_0$, $\tilde\eta_i=\eta_i-1$, $\tilde\eta_{i+1}=1$, $\forall u\in\{1,\ldots,i-1\}, \tilde\eta_u=\eta_u$ thus $\underset{u=0}{\stackrel{i+1}{\sum}}\tilde \eta_u =\underset{u=0}{\stackrel{i}{\sum}}\eta_u +1+2^{i-1}-1 = 2^i$ and $\underset{u=0}{\stackrel{i+1}{\sum}}u \tilde \eta_u =\underset{u=0}{\stackrel{i}{\sum}}u \eta_u + i+1-i = i+1$.
\end{itemize}
\end{proof}

\bibliographystyle{apalike}
\bibliography{BiblioESAIM}

\begin{thebibliography}{}

\bibitem[Birg{\'e}, 2005]{Birge05}
Birg{\'e}, L. (2005).
\newblock A new lower bound for multiple hypothesis testing.
\newblock {\em IEEE Trans. Inform. Theory.}, 51:1611--1615.

\bibitem[Cheney and Light, 2009]{ChL09}
Cheney, W. and Light, W. (2009).
\newblock {\em A course in approximation theory}, volume 101 of {\em Graduate
  Studies in Mathematics}.
\newblock American Mathematical Society, Providence, RI.

\bibitem[Ghosal et~al., 1999]{GGR99}
Ghosal, S., Ghosh, J.~K., and Ramamoorthi, R.~V. (1999).
\newblock Posterior consistency of {D}irichlet mixtures in density estimation.
\newblock {\em Annals of Statistics}, 27:143--158.

\bibitem[Ghosal and van~der Vaart, 2001]{GhosalVaart01}
Ghosal, S. and van~der Vaart, A. (2001).
\newblock Entropy and rates of convergence for maximum likelihood and {B}ayes
  estimation for mixtures of normal densities.
\newblock {\em Annals of Statistics}, 29:1233--1263.

\bibitem[Ghosal and van~der Vaart, 2007]{GhosalVaart07b}
Ghosal, S. and van~der Vaart, A. (2007).
\newblock Posterior convergence rates of {D}irichlet mixtures at smooth
  densities.
\newblock {\em Annals of Statistics}, 35:697--723.

\bibitem[Grenander, 1981]{Grenander}
Grenander, U. (1981).
\newblock {\em Abstract inference}.
\newblock John Wiley and Sons Inc., New York.

\bibitem[Hangelbroek and Ron, 2010]{Hang10}
Hangelbroek, T. and Ron, A. (2010).
\newblock Nonlinear approximation using {G}aussian kernels.
\newblock {\em Journal of Functional Analysis}, 259(1):203--219.

\bibitem[Hartigan, 1975]{Hartigan75}
Hartigan, J.~A. (1975).
\newblock {\em Clustering algorithms}.
\newblock John Wiley \& Sons, New York-London-Sydney.
\newblock Wiley Series in Probability and Mathematical Statistics.

\bibitem[Hastie et~al., 2009]{Hastie09}
Hastie, T., Tibshirani, R., and Friedman, J. (2009).
\newblock {\em The elements of statistical learning}.
\newblock Springer Series in Statistics. Springer, New York, second edition.
\newblock Data mining, inference, and prediction.

\bibitem[Kruijer et~al., 2010]{KRV}
Kruijer, W., Rousseau, J., and van~der Vaart, A. (2010).
\newblock Adaptive {B}ayesian {D}ensity {E}stimation with {L}ocation-{S}cale
  {M}ixtures.
\newblock {\em Electronic Journal of Statistics}, 4:1225--1257.

\bibitem[Lindsay, 1995]{Lindsay95}
Lindsay, B. (1995).
\newblock {\em Mixtures {M}odels: {T}heory, {G}eometry and {A}pplications}.
\newblock IMS, Hayward, CA.

\bibitem[Massart, 2007]{Mas07}
Massart, P. (2007).
\newblock {\em Concentration Inequalities and Model Selection. {\'E}cole
  d'{\'e}t{\'e} de Probabilit{\'e}s de Saint-Flour 2003. Lecture Notes in
  Mathematics}.
\newblock Springer.

\bibitem[Maugis and Michel, 2009]{MaugisMichel09}
Maugis, C. and Michel, B. (2009).
\newblock A non asymptotic penalized criterion for {G}aussian mixture model
  selection.
\newblock {\em ESAIM: Probability and Statistics}.
\newblock doi:10.1051/ps/2009004.

\bibitem[Maugis and Michel, 2010]{MauMic10}
Maugis, C. and Michel, B. (2010).
\newblock Data-driven penalty calibration: a case study for {G}aussian mixture
  model selection.
\newblock {\em ESAIM: Probability and Statistics}.
\newblock doi:10.1051/ps/2010002.

\bibitem[McLachlan and Peel, 2000]{McLachlan:Peel:2000}
McLachlan, G. and Peel, D. (2000).
\newblock {\em Finite {M}ixture {M}odels}.
\newblock Wiley.

\bibitem[Tsybakov, 2009]{Tsybakov09}
Tsybakov, A.~B. (2009).
\newblock {\em Introduction to nonparametric estimation}.
\newblock Springer Series in Statistics. Springer, New York.

\bibitem[Wolfowitz, 1950]{Wolfowitz50}
Wolfowitz, J. (1950).
\newblock Minimax estimation of the mean of a normal distribution with known
  variance.
\newblock {\em Annals of Mathematical Statistics}, 21:218--230.

\end{thebibliography}
\end{document}